\documentclass[11pt,fleqn]{article}

\usepackage{graphicx}
\usepackage{amsthm,amsmath}
\usepackage{amssymb,dsfont,mathrsfs}

\numberwithin{equation}{section}
\newtheorem{theorem}{Theorem}[section]

\newtheorem{proposition}{Proposition}[section]
\newtheorem{corollary}{Corollary}[section]
\newtheorem{remark}{Remark}[section]
\newtheorem{assumption}{Assumption}[section]
\newtheorem{definition}{Definition}[section]
\newtheorem{example}{Example}[section]
\def\ba{\boldsymbol{a}}

\def\be{\boldsymbol{e}}
\def\bg{\boldsymbol{g}}
\def\bj{\boldsymbol{j}}

\def\bw{\boldsymbol{w}}
\def\bx{\boldsymbol{x}}
\def\by{\boldsymbol{y}}
\def\bz{\boldsymbol{z}}
\def\bE{\boldsymbol{E}}
\def\bJ{\boldsymbol{J}}
\def\bO{\boldsymbol{O}}
\def\bX{\boldsymbol{X}}
\def\bY{\boldsymbol{Y}}
\def\bZ{\boldsymbol{Z}}
\def\balpha{\boldsymbol{\alpha}}
\def\bbeta{\boldsymbol{\beta}}
\def\blambda{\boldsymbol{\lambda}}
\def\bmu{\boldsymbol{\mu}}
\def\bphi{\boldsymbol{\phi}}
\def\bpi{\boldsymbol{\pi}}

\def\bzero{\mathbf{0}}
\def\bone{\mathbf{1}}
\def\calA{{\cal A}}
\def\calB{{\cal B}}

\def\calS{{\cal S}}

\def\calV{{\cal V}}

\def\calL{{\cal L}}

\def\scrD{\mathscr{D}}
\def\scrI{\mathscr{I}}
\def\scrP{\mathscr{P}}
\title{Stability of  multidimensional skip-free Markov modulated reflecting random walks%
: Revisit to Malyshev and Menshikov's results and application to queueing networks}
\author{Toshihisa Ozawa \\ 
Faculty of Business Administration, Komazawa University \\
1-23-1 Komazawa, Setagaya-ku, Tokyo 154-8525, Japan \\
E-mail: toshi@komazawa-u.ac.jp }
\date{\today}

\begin{document}

\maketitle

\begin{abstract} 
Let $\{\bX_n\}$ be a discrete-time $d$-dimensional process on $\mathbb{Z}_+^d$ with a supplemental (background)  process $\{J_n\}$ on a finite set and assume the joint process $\{\bY_n\}=\{(\bX_n,J_n)\}$ to be Markovian. Then, the process $\{\bX_n\}$ can be regarded as a kind of reflecting random walk (RRW for short) in which the transition probabilities of the RRW are modulated according to the state of the background process $\{J_n\}$; we assume this modulation is space-homogeneous inside $\mathbb{Z}_+^d$ and on each boundary face of $\mathbb{Z}_+^d$. 
Further we assume the process $\{\bX_n\}$ is skip free in all coordinates and call the joint process $\{\bY_n\}$ a $d$-dimensional skip-free Markov modulated reflecting random walk (MMRRW for short). The MMRRW is an extension of an ordinary RRW and stability of ordinary RRWs have been studied by Malyshev and Menshikov \cite{Malyshev81}. 
Following their results, we obtain stability and instability conditions for MMRRWs and apply our results to stability analysis of a two-station network. 

\smallskip
{\it Keywards}: Multidimensional reflecting random walk, Markov chain, stability, positive recurrent, Foster's criterion, multiclass queueing network

\end{abstract}

%{\small 
%\tableofcontents
%}

%%%%%%%%%%%%%%%%%%%%%%%
%
% Section 1
%
%%%%%%%%%%%%%%%%%%%%%%%
%
\section{Introduction} \label{sec:introduction}

We consider an extension of a multidimensional skip-free reflecting random walk and investigate its stability; then, we apply our results to a multiclass queuing network with two stations.

An ordinary $d$-dimensional reflecting random walk (RRW for short), denoted by $\{\bX_n\}$, is a Markov chain on the state space $\mathbb{Z}_+^d$, where $\mathbb{Z}_+$ is the set of all nonnegative integers and each $\bX_n$ a $d$-dimensional stochastic vector taking values in $\mathbb{Z}_+^d$. RRWs are used as a basic model in many fields; for instance, in the area of queueing theory, a $d$-dimensional RRW $\{\bX_n\}$ is used for representing the behavior of a discrete-time queueing network, where $d$ is the number of queues and $X_n(l)$, the $l$-th element of $\bX_n$, is the number of customers in queue $l$ at time $n$ (see, for example, Miyazawa \cite{Miyazawa11}). 
However, objects that can be represented as ordinary RRWs are rather restricted. In order to see this point, we here consider a multiclass queueing network whose behavior is represented as a continuous-time homogeneous Markov chain (CTMC for short) on the state space $\mathbb{Z}_+^d$; we denote the CTMC by $\{\bar{\bX}_t\}$, where $\bar{X}_t(l)$, the $l$-th element of $\bar{\bX}_t$,  is the number of customers in queue $l$ at time $t$. 
Assume single arrivals and single departures, then the CTMC $\{\bar{\bX}_t\}$ becomes a continuous-time version of a multidimensional skip-free RRW and, by uniformization, we can obtain a (discrete-time) Markov chain $\{\bX_n\}$ that has the same stationary distribution as that of the CTMC, if it exists. Furthermore, the Markov chain $\{\bX_n\}$ becomes a multidimensional skip-free RRW and it can also be regarded as a discrete-time queueing network. 
Therefore, in order to investigate stability of the original queueing network, we can utilize the RRW $\{\bX_n\}$ instead of the CTMC $\{\bar{\bX}_t\}$; but, in this case, the external arrival processes of the original queueing network are restricted to Poisson processes and the service time distributions are restricted to exponential distributions. Furthermore, service disciplines are also restricted to those having memoryless properties; for example, a preemptive-resume priority service is included but a non-preemptive priority service is not. 
Introducing a background process is a way to remove such restrictions and it has often been used, especially in the area of queueing theory (see, for example, Latouche and Ramaswami \cite{Latouche99} and Ozawa \cite{Ozawa04}). 

Let $\{\bX_n\}$ be a $d$-dimensional process on $\mathbb{Z}_+^d$ with a supplemental (background) process $\{J_n\}$ on a finite set and assume the joint process $\{\bY_n\}=\{(\bX_n,J_n)\}$ to be Markovian. Then, the process $\{\bX_n\}$ can be regarded as a kind of RRW in which the transition probabilities of the RRW are modulated according to the state of the background process $\{J_n\}$; we assume this modulation is space-homogeneous inside $\mathbb{Z}_+^d$ and on each boundary face of $\mathbb{Z}_+^d$. 
Further we assume that the process $\{\bX_n\}$ is skip free in all coordinates, which means that the increment of each element of the process takes values in $\{-1,0,1\}$. We call the joint process $\{\bY_n\}$ a $d$-dimensional skip-free Markov modulated reflecting random walk (MMRRW for short) and investigate stability of the MMRRW. 
In queueing theory, this extension of ordinary RRW enables us to deal with multiclass queueing networks having Markovian arrival processes and service times subject to phase-type distributions; furthermore, non-memoryless-type service disciplines including a non-preemptive priority service also become available. 
In order to widen the range of applications, we further assume that the Markov chain $\{\bY_n\}$ does not have to be irreducible but it has just one irreducible class whose states are accessible from every state of the Markov chain; we call such a Markov chain semi-irreducible. For example, a CTMC arising from a famous reentrant line with a preemptive-resume priority service and its corresponding discrete-time Markov chain are semi-irreducible (see Bramson \cite{Bramson08} and Dumas \cite{Dumas97}). In Dumas \cite{Dumas97}, states in the unique irreducible class are called essential.

%%%%%%%%%%%%%%%%%%%%%%%%%
%
Stability of ordinary multidimensional RRWs have been studied in a lot of literature (see, for example, Fayolle et al.\ \cite{Fayolle95} and references therein) and notable results have been obtained by Malyshev and Menshikov \cite{Malyshev81} (also see Fayolle et al.\ \cite{Fayolle95}). In this paper, following their results, we study stability of semi-irreducible multidimensional skip-free MMRRWs. 
Key notions we use are ``induced Markov chain" and ``mean drift vector" introduced by Malyshev and Menshikov \cite{Malyshev81} (also see Fayolle et al.\ \cite{Fayolle95}); an induced Markov chain is a subprocess generated from a multidimensional RRW and the mean drift vector corresponding to the induced Markov chain is the expected increments of the RRW, evaluated by the stationary distribution of the induced Markov chain. 
Malyshev and Menshikov \cite{Malyshev81} dealt with multidimensional RRWs, which were also irreducible Markov chains, and obtained the condition for the RRWs to be positive recurrent (ergodic) and that for them to be transient, where existence of a certain test function was assumed. To construct such a test function seems to be crucial for applying their results. 
Malyshev and Menshikov introduced the notion of second vector field, which was the nonnegative orthant every point of which was assigned one of the mean drift vectors, and they proposed a way to construct a desired test function by using flows on the second vector field. 
Furthermore, applying their results, they classified two and three-dimensional RRWs with respect to stability. 

The notion of induced Markov chain can easily be applied to multidimensional MMRRWs; we define it in a way little bit different from that used in Malyshev and Menshikov \cite{Malyshev81}. First, we consider Markov modulated partially-reflecting random walks obtained from a multidimensional MMRRW $\{(\bX_n,J_n)\}$ by removing some boundaries. We call those random walks expanded Markov chains. Then, the expanded Markov chains are expressed as Markov additive processes and each induced Markov chain is given by the background process of the corresponding Markov additive process. Furthermore, each mean drift vector is given by the expectation of the time-averaged increments of the process $\{\bX_n\}$ of the corresponding expanded Markov chain. Hence, the expanded Markov chains play an important role in our analysis. 

As mentioned above, in applying Malyshev and Menshikov's results (Theorem 2.1 of Ref.\ \cite{Malyshev81}; also see Theorem 4.3.4 of Fayolle et al.\ \cite{Fayolle95}), it is crucial to select a suitable test function and the method using flows on the second vector field is effective for constructing such a test function. However, we consider such a method relying on geometric perception is sometimes too intuitive; hence we adopt linear and piecewise-liner functions as specific test functions.  Using such test functions, we obtain stability and instability conditions for multidimensional skip-free semi-irreducible MMRRWs, where ``stability'' corresponds to ``positive recurrence" in the case of irreducible Markov chains and ``instability" to ``transience" (see Subsection \ref{sec:semi_irreducible} for details); those results are proved by using a kind of Foster's criterion. 
We here note that possibility of linear test functions (Lyapounov functions) has already been suggested in Fayolle \cite{Fayolle89}. 
The conditions we obtain are represented in terms of the mean drift vectors and, in order to derive them, we need approximation formulas for the expected increments of a certain embedded Markov chain of the original MMRRW; the approximation formulas are also represented in terms of the mean drift vectors. 
Applying our results, we classify one through three-dimensional MMRRWs with respect to stability; the classifications we obtain are compatible with those for ordinary RRWs, obtained in Malyshev and Menshikov \cite{Malyshev81}. 
In order for the paper to be self-contained, we give proofs to all our propositions and theorems.

%%%%%%%%%%%%%%%%%%%%%%%%%%%%%%%%%
%
Multiclass queueing networks are typical application of our results. 
Stability of multiclass queueing networks have intensively been studied for the last few decades, especially, by the method using relations between stability of queueing networks and that of the corresponding fluid limits and fluid models (see, for example, Bramson \cite{Bramson08}, Chen \cite{Chen95}, Dai \cite{Dai95,Dai96}, Dai and Vande Vate \cite{Dai96c,Dai00}, Dai and Weiss \cite{Dai96b}, Down and Meyn \cite{Down97},  Gamarnik and Hasenbein \cite{Gamarnik05} and Meyn \cite{Meyn95}). 
A main result in those studies is that the stability of a fluid model implies the stability of the corresponding queueing network (see, for example, Chen \cite{Chen95} and Dai \cite{Dai95}). Several versions of its converse have also been obtained (see, for example, Dai \cite{Dai96} and Gamarnik and Hasenbein \cite{Gamarnik05}). 
It is not doubted that fluid limits and fluid models are the most promising tools to investigate stability of multiclass queueing networks. However, at the same time, we have some difficulty in applying them in some cases; for instance, in order to represent a particular service discipline in a fluid model, we need some additional equation and to describe such an equation is sometimes difficult. A non-preemptive priority service is a typical example; it cannot be represented in an ordinary framework of fluid model. Hence, we consider the MMRRW can be another tool to investigate stability of multiclass queueing networks. 

As discussed in Dumas \cite{Dumas97} and Tezcan \cite{Tezcan13}, each induced Markov chain of a MMRRW arising from a queueing network corresponds to the queueing network in which some queues are saturated with customers; the model of Dumas \cite{Dumas97} is a 6-dimensional RRW and that of Tezcan \cite{Tezcan13} is a two-dimensional MMRRW. Therefore, it can be seen from the queueing network equation that the mean drift vectors can be represented in terms of the external arrival rates of the corresponding queueing network and the departure rates of the saturated queues. This makes the mean drift vectors very easy to get; in other words, we do not need the stationary distributions of the induced Markov chains to get the mean drift vectors if the departure rates of the saturated queues are obtained in a different way. We demonstrate this point in a two-station network.

%%%%%%%%%%%%%%%%%%%%%%%%%%%%
%
The rest of the paper is organized as follows.  
In Section \ref{sec:model}, we describe the $d$-dimensional skip-free MMRRW in detail. Then, we define semi-irreducible Markov chains and discuss stability of them, where a model of Lu-Kumar network is presented as an example of 4-dimensional semi-irreducible MMRRW. 
Section \ref{sec:stability_instability_MMRRW} is a main section of the paper. We give general criteria for stability of semi-irreducible MMRRWs and define expanded Markov chains and induced Markov chains as well as the mean drift vectors. Then, we obtain approximation formulas for the expected increments of a certain embedded Markov chain and give stability and instability conditions for semi-irreducible MMRRWs. After that, we classify low-dimensional MMRRWs with respect to stability. 
In Section \ref{sec:application_QNW}, first we discuss relation between multiclass queueing networks and the corresponding MMRRWs. Then, we consider a two-station network with Markovian arrival processes and service times subject to phase-type distributions. We deal with two examples: one is a two-station network with a non-preemptive priority service and the other that with a $(1,K)$-limited service. The stability and instability regions of both the examples are given and we see that the region of the latter example depends on the value of $K$.

%%%%%%%%%%%%%%%%%%%%%%%
%
% Section 2
%
%%%%%%%%%%%%%%%%%%%%%%%
%
\section{Markov modulated reflecting random walks} \label{sec:model}

%%%%%%%%%%%%%%%%%%%%%%%
%
\subsection{Model description} \label{sec:modeldescription}
Hereafter, we use the following notation. Let $d$ be the dimension of the random walk we consider and let a set $D$ be defined as $D=\{1,2,...,d\}$. For a set $A$, we denote by $\scrP(A)$ the set of all subsets of $A$, including the empty set. We use $\scrP(D)$ as the index set for indicating boundary faces, induced Markov chains and so on. 
For a $d$-dimensional vector $\bx$ on $\mathbb{R}^d$, we denote by $x(l)$ the $l$-th element of $\bx$ and by $\bx(A)$ a part of $\bx$ specified by $A\in\scrP(D)$, i.e., $\bx(A)=(x(l),l\in A)$; for example, when $d=5$ and $A=\{2,4,5\}$, we have $\bx(A)=(x(2),x(4),x(5))$, $\bx(D\setminus A)=(x(1),x(3))$ and $\bx=(\bx(A),\bx(D\setminus A))$. Note that, if $A=\emptyset$, $\bx(A)$ means nothing.
For a set A, we denote by $|A|$ the cardinality of $A$. We have $|\scrP(D)|=2^d$. 
Furthermore, for a real vector $\bx$ and a real number $c$, if every element of $\bx$ is equal to $c$, we express it as $\bx=c$; we analogously define $\bx<c$, $\bx>c$ and so on. 

Consider $d$-dimensional nonnegative orthant $\mathbb{R}_+^d$ and divide $\mathbb{R}_+^d$ into $2^d$ subsets defined by 
\[
\calB^A=\{\bx\in\mathbb{R}_+^d: \bx(A)>0,\ \bx(D\setminus A)=0 \},\ A\in\scrP(D). 
\]
Since we have $\calB^A\cap\calB^B=\emptyset$ if $A\ne B$ and $\mathbb{R}_+^d=\bigcup_{A\in\scrP(D)} \calB^A$, $\{\calB^A: A\in\scrP(D)\}$ is a partition of $\mathbb{R}_+^d$; $\calB^\emptyset$ is the set containing only the origin and $\calB^D$ is the interior of $\mathbb{R}_+^d$. 
For nonempty $A\in\scrP(D)$ such that $A\ne D$, the closure of $\calB^A$, $\{\bx\in\mathbb{R}_+^d: \bx(D\setminus A)=0 \}$% denoted by $\bar{\calB}^A$
, is usually called a boundary face of $\mathbb{R}_+^d$; hence we call $\calB^A$ a sub-boundary face, in the paper.  
Next, we consider $2^d$ finite sets indexed in $\scrP(D)$, say $S^A=\{1,2,...,s^A\}$ for $A\in\scrP(D)$, where $s^A$ is the number of elements of $S^A$. For $A\in\scrP(D)$, assigning the finite set $S^A$ to each point of $(\calB^A\cap\mathbb{Z}^d)$, we define a state space $\calS$ ({\it not} $S$) as 
\begin{equation}
\calS = \bigcup_{A\in\scrP(D)} (\calB^A\cap\mathbb{Z}^d)\times S^A. \label{eq:calS}
\end{equation}
The stochastic process we consider is a discrete-time homogeneous Markov chain $\{\bY_n\}=\{(\bX_n,J_n)\}$ defined on the state space $\calS$, where $\bX_n$ is a $d$-dimensional stochastic vector taking values in $\mathbb{Z}_+^d$ and $J_n$ is the state of the background process taking values in $\bigcup_{A\in\scrP(D)}S^A$. 
We assume that the Markov chain $\{\bY_n\}$ has the properties that the transition probabilities are space-homogeneous in each subset $(\calB^A\cap\mathbb{Z}^d)\times S^A$ in a certain sense and that the process $\{\bX_n\}$ is skip free in all coordinates. 
Then, the Markov chain is a $d$-dimensional reflecting random walk whose transition probabilities are modulated by the background process; we call it a $d$-dimensional skip-free Markov modulated reflecting random walk ($d$-dimensional MMRRW for short). We also denote the Markov chain $\{\bY_n\}$ by $\calL$. 

Here, we explain the homogeneity of the transition probabilities precisely. 
For $n\ge 0$, let $\bZ_{n+1}$ be an increment of the vector process, defined as $\bZ_{n+1}=\bX_{n+1}-\bX_n$, and let $A$ be an element of $\scrP(D)$. If $\bX_n\in(\calB^A\cap\mathbb{Z}^d)$, then $\bZ_{n+1}(A)\in\{-1,0,1\}^{|A|}$ and $\bZ_{n+1}(D\setminus A)\in\{0,1\}^{d-|A|}$ since we have $\bX_n(A)\ge 1$ and $\bX_n(D\setminus A)=0$. 
For $\bx,\,\bz\in\mathbb{R}_+^d$, suppose that $\bx\in(\calB^A\cap\mathbb{Z}^d)$, $\bz(A)\in\{-1,0,1\}^{|A|}$ and $\bz(D\setminus A)\in\{0,1\}^{d-|A|}$. Further suppose that $\bx+\bz\in(\calB^B\cap\mathbb{Z}^d)$ for some $B\in\scrP(D)$ (the case where $B=A$ is included). 
In our model, the homogeneity of the transition probabilities means that if $A$, $B$ and $\bz$ are given, then the transition probability $\mathbb{P}(\bY_{n+1}=(\bx+\bz,j)\,|\,\bY_n=(\bx,i))$, where $i\in S^A$ and $j\in S^B$, does not depend on $\bx$. 
In other words, it is given in the form of 
\begin{equation}
\mathbb{P}(\bY_{n+1}=(\bx+\bz,j)\,|\,\bY_n=(\bx,i)) = p_{\bz}^{A,B}(i,j), \label{eq:tp_calL}
\end{equation}
where $p_{\bz}^{A,B}(i,j)$ is a nonnegative number determined by the values of $A$, $B$, $\bz$, $i$ and $j$; if $B=A$, then $p_{\bz}^{A,B}(i,j)$ is a transition probability that $\bY_n$ changes in $(\calB^A\cap\mathbb{Z}^d)\times S^A$; otherwise, it is a transition probability that $\bY_n$ changes from a state in $(\calB^A\cap\mathbb{Z}^d)\times S^A$ to another state in $(\calB^B\cap\mathbb{Z}^d)\times S^B$. We denote by $P_{\bz}^{A,B}$ the $s^A\times s^B$ matrix whose $(i,j)$-element is $p_{\bz}^{A,B}(i,j)$, i.e., $P_{\bz}^{A,B}=(p_{\bz}^{A,B}(i,j),i\in S^A,j\in S^B)$.

For $\by\in\calS$, we denote by $\balpha_{\by}$ the conditional mean increment vector of $\calL$ given that the state is in $\by$, i.e., $\balpha_{\by}=\mathbb{E}(\bX_{n+1}-\bX_{n}\,|\,\bY_n=\by)$, where we take the expectation element-wise. By the space homogeneity of the transition probabilities, if $\bx,\,\bx'\in (\calB^A\cap \mathbb{Z}^d)$ for some $A\in\scrP(D)$, then we have, for $j\in S^A$, $\balpha_{(\bx,j)}=\balpha_{(\bx',j)}$, which we denote by $\balpha^A_j$.

\begin{example}[QBD process] \label{ex:QBD1}
{\rm 
One-dimensional MMRRW $\{\bY_n\}=\{(X_n,J_n)\}$ is a quasi-birth-and-death process (QBD process for short) on the state space $\calS=(\{0\}\times S^\emptyset)\cup(\mathbb{N}\times S^{\{1\}})$, where $\mathbb{N}$ is the set of positive integers; $X_n$ is called the level and $J_n$ the phase (see, for example, Latouche and Ramaswami \cite{Latouche99}). We have $D=\{1\}$ and $\scrP(D)=\{\emptyset, \{1\}\}$, and the transition matrix $P$ of $\{\bY_n\}$ is given by the block tri-diagonal matrix of 
\begin{equation}
P=\begin{pmatrix} 
P_0^{\emptyset,\emptyset} & P_1^{\emptyset,\{1\}} & & & \cr
P_{-1}^{\{1\},\emptyset} & P_0^{\{1\},\{1\}} & P_1^{\{1\},\{1\}} & & \cr
 & P_{-1}^{\{1\},\{1\}} & P_0^{\{1\},\{1\}} & P_1^{\{1\},\{1\}}  &\cr
 & \ddots & \ddots & \ddots &
\end{pmatrix}. 
\label{eq:QBD_P}
\end{equation}
\hfill$\Box$
}
\end{example}

\begin{example}[Two-dimensional QBD process] \label{ex:2D_QBD1} 
{\rm 
Two-dimensional MMRRW $\{\bY_n\}=\{(\bX_n,J_n)\}$, where $\bX_n=(X_n(1),X_n(2))$, is a Markov chain on the state space
\[
\calS = (\{0\}\times\{0\}\times S^\emptyset)\cup(\mathbb{N}\times\{0\}\times S^{\{1\}})\cup(\{0\}\times\mathbb{N}\times S^{\{2\}})\cup(\mathbb{N}^2\times S^D),
\]
which is called a two-dimensional QBD process in Ozawa \cite{Ozawa12}. We have $D=\{1,2\}$ and $\scrP(D)=\{\emptyset,\{1\},\{2\},D\}$. Transition probabilities are given by $36$ matrices, which are 
\begin{align*}
&P_{\bz}^{D,D},\,\bz\in\{-1,0,1\}^2, \quad 
P_{(z,-1)}^{D,\{1\}},\,z\in\{-1,0,1\}, \quad
P_{(-1,z)}^{D,\{2\}},\,z\in\{-1,0,1\}, \quad
P_{(-1,-1)}^{D,\emptyset}, \cr
&P_{(z,1)}^{\{1\},D},\,z\in\{-1,0,1\}, \quad
P_{(z,0)}^{\{1\},\{1\}},\,z\in\{-1,0,1\}, \quad
P_{(-1,1)}^{\{1\},\{2\}}, \quad
P_{(-1,0)}^{\{1\},\emptyset}, \cr
&P_{(1,z)}^{\{2\},D},\,z\in\{-1,0,1\}, \quad
P_{(1,-1)}^{\{2\},\{1\}}, \quad
P_{(0,z)}^{\{2\},\{2\}},\,z\in\{-1,0,1\}, \quad
P_{(0,-1)}^{\{2\},\emptyset}, \cr
&P_{(1,1)}^{\emptyset,D},\quad
P_{(1,0)}^{\emptyset,\{1\}}, \quad
P_{(0,1)}^{\emptyset,\{2\}}, \quad
P_{(0,0)}^{\emptyset,\emptyset}.
\end{align*} 
For example, $P_{\bz}^{D,D}$ governs transitions in $(\calB^D\cap\mathbb{Z}^2)\times S^D=\mathbb{N}^2\times S^D$ and $P_{(z,-1)}^{D,\{1\}}$ those from states in $\mathbb{N}\times\{1\}\times S^D=\{((x,1),i): x\in\mathbb{N},\,i\in S^D\}$ to other states in $\mathbb{N}\times\{0\}\times S^{\{1\}}=\{((x',0),j): x'\in\mathbb{N},\,j\in S^{\{1\}}\}$; $P_{(-1,-1)}^{D,\emptyset}$ governs transitions from states in $\{1\}\times\{1\}\times S^D=\{((1,1),i): i\in S^D\}$ to other states in $\{0\}\times\{0\}\times S^\emptyset=\{((0,0),j): j\in S^\emptyset\}$. The roles of other matrices are analogously given. 
\hfill$\Box$
}
\end{example}

%%%%%%%%%%%%%%%%%%%%%%%%
% figure: Re-entrant line
%%%%%%%%%%%%%%%%%%%%%%%%
\begin{figure}[bht]
\begin{center}
\setlength{\unitlength}{0.8mm}
\begin{picture}(100,50)(0,0)
%
% Station 1
\thicklines
\put(2,15){\makebox(0,0){\normalsize $\bar{\lambda}$}}
\put(5,15){\vector(1,0){7}}

\put(14,6){\makebox(0,0){\normalsize Q$_1$}}
\put(27,11){\makebox(0,0){\normalsize $\bar{h}_1$}}
\put(12,20){\line(1,0){6}}
\put(12,10){\line(1,0){6}}
\put(18,10){\line(0,1){10}}
\multiput(22,15)(2,0){6}{\line(1,0){1}}
\put(36,15){\vector(1,0){15}}

\put(42,40){\makebox(0,0){\normalsize Q$_4$}}
\put(22,34){\makebox(0,0){\small High priority}}
\put(27,26){\makebox(0,0){\normalsize $\bar{h}_4$}}
\put(37,25){\line(1,0){6}}
\put(37,35){\line(1,0){6}}
\put(37,25){\line(0,1){10}}
\multiput(22,30)(2,0){6}{\line(1,0){1}}
\put(19,30){\vector(-1,0){15}}

\put(28,0){\makebox(0,0){\normalsize Station 1}}
\put(20,5){\line(0,1){35}}
\put(20,5){\line(1,0){15}}
\put(20,40){\line(1,0){15}}
\put(35,5){\line(0,1){35}}
%
% Station 2
\thicklines
\put(54,6){\makebox(0,0){\normalsize Q$_2$}}
\put(75,11){\makebox(0,0){\small High priority}}
\put(67,18){\makebox(0,0){\normalsize $\bar{h}_2$}}
\put(52,20){\line(1,0){6}}
\put(52,10){\line(1,0){6}}
\put(58,10){\line(0,1){10}}
\multiput(62,15)(2,0){6}{\line(1,0){1}}
\put(76,15){\line(1,0){15}}
\put(91,15){\line(0,1){15}}
\put(91,30){\vector(-1,0){8}}

%\put(95,32){\vector(-1,0){12}}

\put(82,40){\makebox(0,0){\normalsize Q$_3$}}
\put(67,33){\makebox(0,0){\normalsize $\bar{h}_3$}}
\put(77,25){\line(1,0){6}}
\put(77,35){\line(1,0){6}}
\put(77,25){\line(0,1){10}}
\multiput(62,30)(2,0){6}{\line(1,0){1}}
\put(59,30){\vector(-1,0){15}}

\put(68,0){\makebox(0,0){\normalsize Station 2}}
\put(60,5){\line(0,1){35}}
\put(60,5){\line(1,0){15}}
\put(60,40){\line(1,0){15}}
\put(75,5){\line(0,1){35}}

\end{picture}
\caption{Lu-Kumar network.}
\label{fig:reentrant}
\end{center}
\end{figure}
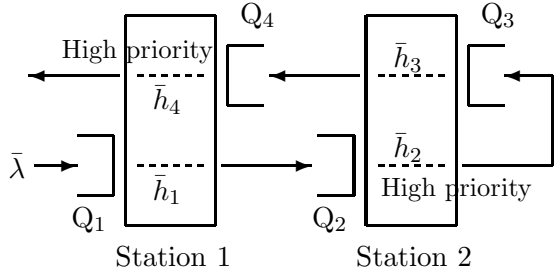
%%%%%%%%%%%%%%%%%%%%%%%

\begin{example}[Lu-Kumar network with a non-preemptive priority service] \label{ex:reentrant1} 
{\rm 
We consider a re-entrant line depicted in Fig.\ \ref{fig:reentrant}, called a Lu-Kumar network (see Lu and Kumar \cite{Lu91}; also see Bramson \cite{Bramson08}). In the model, exogenous customers arrive at station 1 according to a Poisson process with intensity $\bar{\lambda}$ and join queue 1 (Q$_1$); after completing service in station 1, they move to station 2, join queue 2 (Q$_2$) and receive service there; then, they reenter station 2, join queue 3 (Q$_3$) and receive service there again. Finally, they move to station 1, join queue 4 (Q$_4$) and, after completing service there, they depart from the system. 
We assume that, for $l\in\{1,2,3,4\}$, service times for customers in Q$_l$ are subject to an exponential distribution with mean $\bar{h}_l$ and the arrival process and service times are mutually independent. 
In each station, there is just one server which serves customers in two queues according to a non-preemptive priority service; we assume customers in Q$_4$ (resp.\ Q$_2$) have non-preemptive priority over those in Q$_1$ (resp.\ Q$_3$). 

Let $D$ be defined as $D=\{1,2,3,4\}$. For $l\in D$, let $\bar{X}_t(l)$ be the number of customers in Q$_l$ at time $t$, including one being served, and denote by $\bar{\bX}_t$ the vector of $\bar{X}_t(l)$'s, i.e., $\bar{\bX}_t=(\bar{X}_t(l),l\in D)$. Let $\bar{\bJ}_t=(\bar{J}_t(1),\bar{J}_t(2))\in\{0,1,2\}^2$ be the state of the servers at time $t$, given as follow: for $i\in\{1,2\}$, if $\bar{J}_t(i)=0$, then the server in station $i$ is idle; if $\bar{J}_t(i)=1$, then it is serving a high-priority customer; otherwise ($\bar{J}_t(i)=2$), it is serving a low-priority customer. 
The joint process $\{\bar{\bY}_t\}=\{(\bar{\bX}_t,\bar{\bJ}_t)\}$ becomes a CTMC, which is a continuous-time version of a 4-dimensional MMRRW. We have
\begin{align*}
&S^\emptyset=\{(0,0)\},\ 
S^{\{1\}}=\{(2,0)\},\ 
S^{\{2\}}=\{(0,1)\},\ 
S^{\{3\}}=\{(0,2)\},\ 
S^{\{4\}}=\{(1,0)\},\cr
&S^{\{1,2\}}=\{(2,1)\},\ 
S^{\{1,3\}}=\{(2,2)\},\ 
S^{\{1,4\}}=\{(1,0),(2,0)\},\cr 
&S^{\{2,3\}}=\{(0,1),(0,2)\},\ 
S^{\{2,4\}}=\{(1,1)\},\ 
S^{\{3,4\}}=\{(1,2)\},\cr
&S^{\{1,2,3\}}=\{(2,1),(2,2)\},\ 
S^{\{2,3,4\}}=\{(1,1),(1,2)\},\ 
S^{\{1,3,4\}}=\{(1,2),(2,2)\},\cr 
&S^{\{1,2,4\}}=\{(1,1),(2,1)\},\ 
S^D=\{(1,1),(1,2),(2,1),(2,2)\}, 
\end{align*}
and the state space of $\{\bar{\bY}_t\}$ is given by $\calS=\bigcup_{A\in\scrP(D)} (\calB^A\cap\mathbb{Z}^4)\times S^A$, where $\scrP(D)$ is the set of all subsets of $D$. 
For $\by,\by'\in\calS$ such that $\by\ne\by'$, we denote by $q(\by,\by')$ the transition rate that the CTMC changes from state $\by$ to state $\by'$; for $\by\in\calS$, we define $q(\by,\by)$ as $q(\by,\by)=-\sum_{\by'\ne\by}q(\by,\by')$. We omit describing explicit formulas for $q(\by,\by')$'s since they are complicated and we do not use them here. 
Since $|q(\by,\by)|\le \bar{\lambda}+\sum_{l\in D}1/\bar{h}_l$ for any $\by\in\calS$, we can obtain, by uniformization, a discrete-time Markov chain $\{\bY_n\}=\{(\bX_n,\bJ_n)\}$ on the state space $\calS$ that has the same stationary distribution as that of the original CTMC, if it exists. 
The transition probabilities of $\{\bY_n\}$ are given by $\mathbb{P}(\bY_1=\by'|\bY_0=\by)=\delta_{\by,\by'}+q(\by,\by')/\nu$ for $\by,\by'\in\calS$, where $\delta_{\cdot,\cdot}$ is the delta function (if $\by=\by'$, then $\delta_{\by,\by'}=1$; otherwise $\delta_{\by,\by'}=0)$ and $\nu$ is a positive number satisfying $\nu\ge\sup_{\by\in\calS}|q(\by,\by)|$. 
The Markov chain $\{\bY_n\}$ is a 4-dimensional MMRRW on the state space $\calS$ and it is also the discrete-time queueing network corresponding to the Lu-Kumar network we consider; in the discrete-time queueing network, the external arrival rate is given by $\bar{\lambda}/\nu$ and, for $l\in D$, the mean service time of customers in Q$_l$ is given by $\nu\bar{h}_l$. 
\hfill$\Box$
}
\end{example}

%%%%%%%%%%%%%%%%%%%%%%%
%
\subsection{Semi-irreducibility and Stability} \label{sec:semi_irreducible}
The assumption of irreducibility for Markov chains is sometimes too restrictive for analyzing stability of queueing models; for example, the 4-dimensional MMRRW arising from the Lu-Kumar network with a non-preemptive priority service is a reducible Markov chain, which has just one irreducible class (closed communication class); we will explain this point in detail later.  Hence, we consider a class of Markov chain including such reducible Markov chains and call ones in the class semi-irreducible Markov chains. 

\begin{definition}[Semi-irreducible Markov chain] 
A continuous-time or discrete-time homogenous Markov chain on a countable set is said to be semi-irreducible if the chain has just one irreducible class and every state in the irreducible class is accessible from any state of the chain. 
\end{definition}

We note that an irreducible Markov chain is also semi-irreducible. States of the irreducible class are called {\it essential states} (see, for example, Dumas \cite{Dumas97}). 

\begin{remark} \label{re:semiirreducible}
Suppose that a Markov chain has a state that is accessible from any state of the chain, then the communicating class including the state is closed and every state in the class is accessible from any state of the chain. Furthermore, the closed communicating class is unique since if it is not the case, then there exists a state from which any state in the closed communicating class is not accessible and this is a contradiction to the assumption. As a result, the Markov chain is semi-irreducible. 
Therefore, it is not so difficult to check up on semi-irreducibility of a given Markov chain arising from an open queueing network, since a state in which the system is empty is usually accessible from any other states of the model, especially, in the case where the open queueing network is work-conserving. 
\end{remark}

\begin{example}[Semi-irreducibility of the Lu-Kumar network with a non-preemptive priority service] \label{ex:reentrant2}
{\rm 
Consider the CTMC $\{\bar{\bY}_t\}=\{(\bar{\bX}_t,\bar{\bJ}_t)\}$ on the state space $\calS$ arising from the Lu-Kumar network described in Example \ref{ex:reentrant1}. State $\by_0=((0,0,0,0),(0,0))$, in which the system is empty and both the servers are idle, is obviously accessible from any other state in $\calS$ since the Lu-Kumar network is work-conserving. Thus, the CTMC $\{\bar{\bY}_t\}$ is semi-irreducible and the 4-dimensional discrete-time MMRRW derived from the CTMC by uniformization, $\{\bY_n\}$, is also semi-irreducible. 
In addition, they are not irreducible; the reason is as follows.
Consider a state in which the number of customers in Q$_4$ is 2, say $\by_1=((x_1(1),x_1(2),x_1(3),2),(j_1(1),j_1(2)))$. In order for the CTMC $\{\bar{\bY}_t\}$ to reach $\by_1$ from $\by_0$, a customer must enter Q$_4$ when there is just one customer in Q$_4$; we call the former customer (arriving customer) the second customer and the latter the first customer. 
Let $t_1$ be the time when the first customer enters Q$_4$, $t_2$ the time when service for the second customer in Q$_3$ begins, $t_3$ the time when the second customer enters Q$_4$ and $t_4$ the time when the number of customers in Q$_4$ becomes 1 for the first time after $t_3$. We have $t_1\le t_2<t_3<t_4$ and also have $\bar{X}_{t_3-}(4)=1$, $\bar{X}_{t_4}(4)=1$ and, for every $t\in[t_3,t_4)$, $\bar{X}_t(4)\ge 2$. 
Since customers in Q$_2$ have non-preemptive priority over those in Q$_3$, the number of customers in Q$_2$ at $t_2$ must be 0, i.e., $\bar{X}_{t_2}(2)=0$. Furthermore, since customers in Q$_4$ have non-preemptive priority over those in Q$_1$, the number of customers whose service in Q$_1$ is completed in $[t_1,t_4)$ is at most 1. Thus, we must have, for $t\in[t_2,t_4)$, $\bar{X}_t(2)\le \bar{X}_{t_2}(2)+1 =1$. 
This implies that when the number of customers in Q$_4$ is greater than or equal to 2, the number of customers in Q$_2$ must be less than or equal to 1. Therefore, states in $\{(\bx,\bj)\in\calS: x(2)\ge 2\ \mbox{and}\ x(4)\ge 2\}$ are not accessible from $\by_0$ and the CTMC $\{\bar{\bY}_t\}$ is reducible; thus, the 4-dimensional MMRRW $\{\bY_n\}$ is also reducible.
\hfill$\Box$
}
\end{example}

Hereafter, we assume the following condition.
\begin{assumption}
The $d$-dimensional MMRRW $\calL=\{\bY_n\}$ is semi-irreducible.
\end{assumption}

%%%%%%%%%%%%%%%%%%%%%
Next, we classify semi-irreducible MMRRWs with respect to stability; the classification can also be applied to general semi-irreducible Markov chains. 
Let $\{\bY_n\}=\{(\bX_n,J_n)\}$ be a $d$-dimensional semi-irreducible MMRRW on the state space $\calS$ with the unique irreducible class $\calS_0$. 
The MMRRW can have more than one communication class and the number of communication classes may be countable. From the definition of semi-irreducibility, all communication classes except for $\calS_0$ are open and every state in those open communication classes are transient in usual sense. 
However, it is possible that there exists a open communication class, say $\calS_{open}\subset\calS\setminus\calS_0$, for which we have $\lim_{n\to\infty}\mathbb{P}(\bY_n\in\calS_{open}\,|\,\bY_0=\by_0)>0$ for any state $\by_0\in\calS_{open}$; in this case, we should say that the MMRRW $\{\bY_n\}$ is unstable. 
We, therefore, consider the following three conditions for MMRRWs:
\begin{description}
\item[(C1)] For every state in $\calS\setminus\calS_0$, the first passage time from the state to the irreducible class is finite with probability one and has the finite expectation. 
\item[(C2)] For every state in $\calS\setminus\calS_0$, the first passage time from the state  to the irreducible class is finite with probability one, but there exists a state in $\calS\setminus\calS_0$ such that the first passage time from the state to the irreducible class has no finite expectation. 
\item[(C3)] There exists a state in $\calS\setminus\calS_0$ such that, with a positive probability, the MMRRW starting from the state will never reach the irreducible class. 
\end{description}
With respect to the irreducible class, we have the following three conditions: {\bf (C4)} the irreducible class is positive recurrent; {\bf (C5)} it is null recurrent; {\bf (C6)} it is transient, in usual sense. 
The semi-irreducible MMRRW $\{\bY_n\}=\{(\bX_n,J_n)\}$ satisfying conditions C1 and C4 is stable in the sense that the random walk $\{\bX_n\}$ starting from any state will never diverge; in this case, the MMRRW is also positive Harris recurrent (see, for example, Bramson \cite{Bramson08}). 
The semi-irreducible MMRRW $\{\bY_n\}=\{(\bX_n,J_n)\}$ satisfying condition C6 is unstable in the sense that the random walk $\{\bX_n\}$ starting from any state will diverge; in this case, all the states of the MMRRW are transient.
The semi-irreducible MMRRW $\{\bY_n\}=\{(\bX_n,J_n)\}$ satisfying conditions C3 and C4 or conditions C3 and C5 is unstable in the sense that the random walk $\{\bX_n\}$ starting from some state will diverge; however, in this case, the random walk starting from any state in the irreducible class never diverge. 
Considering those points, we will use the following definition of stability for semi-irreducible MMRRWs.
\begin{definition}[Stability of semi-irreducible MMRRWs] \label{def:stability}
The semi-irreducible MMRRW $\{\bY_n\}$ is said to be stable if it satisfies conditions C1 and C4, and it is said to be unstable if it satisfies condition C3 or condition C6. 
\end{definition}

We also apply this definition to continuous-time models of semi-irreducible MMRRW. 
\begin{remark}
In Definition \ref{def:stability}, we exclude three cases: the first is that the MMRRW satisfies conditions C1 and C5, the second is that it satisfies conditions C2 and C4, and the third is that it satisfies conditions C2 and C5. In these cases, the MMRRW is Harris recurrent (see, for example, Bramson \cite{Bramson08}), but is not positive Harris recurrent, and it is difficult to analyze stability of the MMRRW by the methods we will use in the paper; this is a reason why we exclude the cases. 
\end{remark}

\begin{example}[Stability of a combination of two random walks] \label{ex:2D_RW1}
{\rm 
Let $\{X_n(1)\}$ be a one-dimensional skip-free random work on $\mathbb{Z}_+$ having the reflecting barrier at zero and $\{X_n(2)\}$ another one-dimensional skip-free random walk on $\mathbb{Z}_+$ having the absorbing barrier at zero. 
The transition probabilities of $\{X_n(1)\}$ are given as, for $i\ge 0$, $\mathbb{P}(X_1(1)=i+1\,|\,X_0(1)=i)=p_1$ and $\mathbb{P}(X_1(1)=i\,|\,X_0(1)=i+1)=1-p_1$, where we have $0\le p_1\le 1$ and $\mathbb{P}(X_1(1)=0\,|\,X_0(1)=0)=1-p_1$; those of $\{X_n(2)\}$ are given as, for $i\ge 1$, $\mathbb{P}(X_1(2)=i+1\,|\,X_0(2)=i)=p_2$ and $\mathbb{P}(X_1(2)=i-1\,|\,X_0(2)=i)=1-p_2$, where we have $0\le p_2\le 1$, $\mathbb{P}(X_1(2)=0\,|\,X_0(2)=0)=1$ and $\mathbb{P}(X_1(2)=1\,|\,X_0(2)=0)=0$. 
Assuming that $\{X_n(1)\}$ and $\{X_n(2)\}$ are mutually independent, we consider the process combining those two random walks and denoted it by $\{\bX_n\}=\{(X_n(1),X_n(2))\}$. Then, the process $\{\bX_n\}$ is a two-dimensional semi-irreducible random walk on $\mathbb{Z}_+^2$, whose unique irreducible class is given by $\calS_0=\{\bx\in\mathbb{Z}_+^2: x(2)=0\}$, and it is also a semi-irreducible MMRRW. 
If $p_1<1/2$ and $p_2<1/2$, then the process $\{\bX_n\}$ satisfies conditions C1 and C4, and it is stable; if $p_1>1/2$, then the process satisfies condition C6 and it is unstable; if $p_2>1/2$, then for any $i\ge 0$ and any $j\ge 1$, the process starting from state $(i,j)$ will never reach the irreducible class with probability $1-((1-p_2)/p_2)^j$, and hence the process is unstable due to condition C3. 
\hfill$\Box$
}
\end{example}

Since it is not always easy to separate the irreducible class from the whole state space, we will seek a way to analyze stability of the semi-irreducible MMRRW {\it without specifying the irreducible class} in the next section. 
Besides, introducing semi-irreducibility has the following merit.

\begin{remark} \label{re:semiirreduciblemodel}
For a $d$-dimensional semi-irreducible MMRRW $\{\bY_n\}=\{(\bX_n,J_n)\}$ on the state space $\calS$, adding dummy states to $\calS$, we can define a new state space with a simpler structure and obtain a semi-irreducible MMRRW on the new state space corresponding to the original MMRRW. 
For example, defining an extended state space $\calS^{ex}$ as 
\[
\calS^{ex} = \bigcup_{A\in\scrP(D)} (\calB^A\cap\mathbb{Z}^d)\times S^{\scrP(D)} = \mathbb{Z}_+^d\times S^{\scrP(D)}, 
\]
where $S^{\scrP(D)}=\bigcup_{A\in\scrP(D)} S^A$, we can obtain a $d$-dimensional MMRRW $\{\bY^{ex}_n\}=\{(\bX^{ex}_n,J^{ex}_n)\}$ on $\calS^{ex}$ whose transition probabilities are given by 
\[
\mathbb{P}(\bY^{ex}_1=\by'\,|\,\bY^{ex}_0=\by)
=\left\{ \begin{array}{ll}
\mathbb{P}(\bY_1=\by'\,|\,\bY_0=\by), & \by,\by'\in\calS, \cr
1, & \by\in\calS^{ex}\setminus \calS,\,\by'=\by_0, \cr
0, & \by\in\calS^{ex}\setminus \calS,\,\by'\in\calS^{ex}\setminus\{\by_0\}, 
\end{array} \right.
\]
where $\by_0$ is a state belonging to the unique irreducible class of the original MMRRW. 
The new MMRRW $\{\bY^{ex}_n\}$ is also semi-irreducible and has the same irreducible class as that of the original MMRRW. Furthermore, the stability of $\{\bY^{ex}_n\}$ coincides with that of the original MMRRW. 
\end{remark}

%%%%%%%%%%%%%%%%%%%%%%%
%
% Section 3
%
%%%%%%%%%%%%%%%%%%%%%%%
%
\section{Stability of semi-irreducible MMRRWs} \label{sec:stability_instability_MMRRW}

Let $\calL=\{\bY_n\}=\{(\bX_n,J_n)\}$ be a $d$-dimensional semi-irreducible MMRRW on the state space $\calS$ with the unique irreducible class $\calS_0$. In this section, we derive stability and instability conditions for the MMRRW $\calL$. 
Key notions we use are induced Markov chains and the mean drift vectors introduced by Malyshev and Menshikov \cite{Malyshev81} (also see Fayolle et al.\ \cite{Fayolle95}), in which stability of multidimensional reflecting random walks (RRWs) was discussed. We apply those notions to our multidimensional MMRRWs. 
We also reconsider the definition of induced Markov chains and that of the mean drift vectors; briefly speaking, an induced Markov chain is given as the background process of a multidimensional Markov additive process obtained from $\calL$ by removing several reflecting barriers and the corresponding mean drift vector is given as the vector of the mean increments (drifts) of the Markov additive process. 
The stability and instability conditions we derive are represented in terms of the mean drift vectors.

%%%%%%%%%%%%%%%%%%%%%%%%
%
\subsection{Preliminaries}

Here we present several properties of semi-irreducible MMRRWs, which will be used for deriving stability and instability conditions for the MMRRWs in the next subsection. 
If the semi-irreducible MMRRW $\calL=\{\bY_n\}$ is stable in our sense, it has the stationary distribution, denoted by $\bpi=(\pi(\by),\by\in\calS)$, and the ergodic theorem of Markov chains holds, i.e., for a real function $f$ on $\calS$ satisfying 
\[
\sum_{\by\in\calS} |f(\by)|\,\pi(\by) < \infty, 
\]
we have, for any initial distribution, 
\[
\lim_{n\to\infty} \frac{1}{n}\sum_{k=1}^n f(\bY_n) = \sum_{\by\in\calS} f(\by) \pi(\by),\quad a.s.
\]
Note that, since the MMRRW is semi-irreducible with the irreducible class $\calS_0$, we have $\pi(\by)=0$ for all $\by\in\calS\setminus\calS_0$.

%%%%%%%%%
Exclusively divide the state space $\calS$ into a finite number of nonempty subsets, $\calV_k,\,k=1,2,...,m$, where $\calV_k\cap\calV_l=\emptyset$ for $k\ne l$ and $\bigcup_{k=1}^m \calV_k=\calS$, and define a function $u$ on $\calS$ as 
\[
u(\by) = u_k\quad \mbox{if $\by\in\calV_k$ for some $k\in\{1,2,...,m\}$}, 
\]
where $u_k$ is a positive integer. 
Furthermore, define a strictly increasing random sequence $\{\sigma_n\}$ by 
\[
\sigma_0=0,\quad \sigma_{n+1} = \sigma_n+u(\bY_{\sigma_n})\ \mbox{for}\ n\ge 0,
\]
and consider a stochastic process $\{\tilde{\bY}_n\}$ defined by $\tilde{\bY}_n=\bY_{\sigma_n}$ for $n\ge 0$. 
This process $\{\tilde{\bY}_n\}$ is a kind of embedded Markov chain of $\{\bY_n\}$. In order to derive a stability condition for semi-irreducible MMRRWs, we will use the following proposition, which is a modification of Theorem 2.2.4 of Fayolle et al.\ \cite{Fayolle95}; also see Theorem 1.4 of Malyshev and Menshikov \cite{Malyshev81} and Proposition 4.5 of Bramson \cite{Bramson08}.

\begin{proposition} \label{pr:Foster2} 
The semi-irreducible MMRRW $\{\bY_n\}$ is stable in our sense if there exist a positive number $\varepsilon$, a finite subset $\calV\subset\calS$ and a lower bounded real function $f$ on $\calS$ such that
\begin{align}
&\mathbb{E}(f(\tilde{\bY}_{n+1})-f(\tilde{\bY}_n)\,|\,\tilde{\bY}_n=\by) \le -\varepsilon,\quad \by\in\calS\setminus\calV, \label{eq:Foster2a} \\
&\mathbb{E}(f(\tilde{\bY}_{n+1})\,|\,\tilde{\bY}_n=\by) < \infty,\quad \by\in\calV. \label{eq:Foster2b}
\end{align}
\end{proposition}

Usually, when the MMRRW $\{\bY_n\}$ is irreducible and $u(\by)=1$ for all $\by\in\calS$, this proposition is called {\it Foster's criterion}. 

\begin{proof}[Proof of Proposition \ref{pr:Foster2}]
This proposition can be proved in a manner similar to that used for irreducible Markov chains. 
Let $T_{\calV}$ be the time of first entry (return) of $\{\bY_n\}$ into $\calV$ and $\tilde{T}_{\calV}$ that of $\{\tilde{\bY}_n\}$. 
Condition (\ref{eq:Foster2a}) implies that the time of first entry of $\{\tilde{\bY}_n\}$ starting from any state in $\calS\setminus\calV$ into $\calV$ has the finite expectation. This and condition (\ref{eq:Foster2b}) imply that the return time of $\{\tilde{\bY}_n\}$ from any state in $\calV$ to $\calV$ also has the finite expectation (see, for example, the proof of Theorem 1.1 in Chapter 5 of Br\'emaud \cite{Bremaud99}). 
Hence, we have $\mathbb{E}(\tilde{T}_{\calV}\,|\,\tilde{\bY}_0=\by)<\infty$ for all $\by\in\calS$, and this leads us to 
\begin{align}
\mathbb{E}(T_{\calV}\,|\,\bY_0=\by) \le \left( \max_{1\le k\le m} u_k \right) \mathbb{E}(\tilde{T}_{\calV}\,|\,\tilde{\bY}_0=\by)<\infty\ \mbox{for all}\ \by\in\calS. \label{eq:ETy}
\end{align}
If $\calV$ contains no elements in the irreducible class $\calS_0$, the process starting from a state $\by_0\in\calS_0\subset\calS\setminus\calV$ will never reach the set $\calV$. This contradicts expression (\ref{eq:ETy}), the finiteness of the expected time of first entry of $\{\bY_n\}$ starting from the state $\by_0$ into $\calV$, and hence we have $\calV\cap\calS_0\ne\emptyset$. 

Let $\tau_1,\,\tau_2,\,...$ be the successive return times of $\{\bY_n\}$ to $\calV$, and define a process $\{\hat{\bY}_n\}$ by $\hat{\bY}_0=\bY_0$ and $\hat{\bY}_n=\bY_{\tau_n}$ for $n\ge 1$. Then, by the strong Markov property, $\{\hat{\bY}_n\}$ becomes a Markov chain whose state space is $\calV$. 
Select $\by\in\calV\cap\calS_0$; since the original process $\{\bY_n\}$ is semi-irreducible, the state $\by$ in the new process $\{\hat{\bY}_n\}$ is also accessible from any state in $\calV$ and, by Remark \ref{re:semiirreducible}, $\{\hat{\bY}_n\}$ is semi-irreducible. 
For $\by\in\calV$, let $\hat{T}_{\by}$ be the first passage time (return time) of $\{\hat{\bY}_n\}$ to $\by$. Because of the fact that $\calV$ is finite, we have, for any $\by_0\in\calV$ and $\by_1\in\calV\cap\calS_0$, $\mathbb{E}(\hat{T}_{\by_1}\,|\,\hat{\bY}_0=\by_0) < \infty$.
Since the return times $\tau_i,\,i=1,2,...$, are stopping times of $\{\bY_n\}$, we obtain, for any $\by_0\in\calV$ and $\by_1\in\calV\cap\calS_0$, (see, for example, the proof of Lemma 1.1 in Chapter 5 of Br\'emaud \cite{Bremaud99})
\begin{align}
&\mathbb{E}(T_{\by_1}\,|\,\bY_0=\by_0) \le \left(\max_{\by\in\calV} \mathbb{E}(T_{\calV}\,|\,\bY_0=\by) \right) \mathbb{E}(\hat{T}_{\by_1}\,|\,\hat{\bY}_0=\by_0) < \infty.  \label{eq:rtime}
\end{align}
When $\by_0=\by_1\in\calV\cap\calS_0$, expression (\ref{eq:rtime}) implies that the state $\by_1$ as well as the irreducible class $\calS_0$ of $\{\bY_n\}$ is positive recurrent. Furthermore, for any $\by_0\in\calS\setminus\calS_0$ and $\by_1\in\calV\cap\calS_0$, we have 
\[
\mathbb{E}(T_{\by_1}\,|\,\bY_0=\by_0) \le \mathbb{E}(T_{\calV}\,|\,\bY_0=\by_0) + \max_{\by\in\calV} \mathbb{E}(T_{\by_1}\,|\,\bY_0=\by) < \infty, 
\]
and this implies that the time of first entry of $\{\bY_n\}$ starting from any state in $\calS\setminus\calS_0$ into $\calS_0$ has the finite expectation. As a result, $\{\bY_n\}$ is stable in our sense. 
\end{proof}

\begin{remark}
In applying Proposition \ref{pr:Foster2} to a semi-irreducible MMRRW, it is not necessary to specify the irreducible class of the MMRRW; this point is very convenient since to specify the irreducible class of a given semi-irreducible MMRRW is not usually so easy. 
\end{remark}

%%%%%%%%%
In order to derive an instability condition for  semi-irreducible MMRRWs, we will use the following proposition, which is a version of Theorem 2.2.7 of Fayolle et al.\ \cite{Fayolle95} (also see Theorem 1.6 of Malyshev and Menshikov \cite{Malyshev81}). 
\begin{proposition} \label{pr:Markov_unstable2} 
The semi-irreducible MMRRW $\{\bY_n\}$ is unstable in our sense if there exist a positive function $f$ on $\calS$ and positive numbers $\varepsilon$, $c$ and $b$ such that, for $\calA=\{\by\in\calS : f(\by)>c \}$,   
\begin{itemize}
\item[(i)] $\calA\ne\emptyset,\ \calA^C\cap\calS_0\ne\emptyset$, 
\item[(ii)] $\mathbb{E}(f(\tilde{\bY}_{n+1})-f(\tilde{\bY}_n)\,|\,\tilde{\bY}_n=\by) \ge \varepsilon$ \ for all $\by\in\calA$, and
\item[(iii)] the inequality $|f(\by_1)-f(\by_0)|>b$ implies $\mathbb{P}(\bY_1=\by_1\,|\,\bY_0=\by_0)=0$. 
\end{itemize}
\end{proposition}

\begin{proof}
This proposition can be proved in a manner similar to that used for proving Theorem 2.2.7 of Fayolle et al.\ \cite{Fayolle95} since Theorem 2.1.9 of the same book can be applied to semi-irreducible MMRRWs. 
Define process $\{Z_n\}$ as $Z_n=f(\bY_n)$ for $n\ge 0$ and let $T_c$ be the time of first entry of $\{Z_n\}$ into $[0,c]$; $T_c$ is also the time of first entry of $\{\bY_n\}$ into $\calA^C$. 
Since $\tilde{\bY}_n=\bY_{\sigma_n}$ for $n\ge 0$, condition (ii) implies that, for all $n\ge 0$, if $Z_{\sigma_n}>c$, then
\begin{equation}
\mathbb{E}(Z_{\sigma_n+1}\,|\,Z_{\sigma_n}) \ge Z_{\sigma_n} + \varepsilon,\quad a.s. 
\label{eq:Zsigma_epsilon}
\end{equation}
Furthermore, since the process $\{\bX_n\}$ is skip free in all coordinates, we see from condition (ii) that there exists $\by\in\calA$ such that $f(\by)> c+b\,\max_{1\le k\le m}u_k$. 
Therefore, letting $\by_0\in\calA$ be such a state, we obtain, by Theorem 2.1.9 of Fayolle et al.\ \cite{Fayolle95}, $\mathbb{P}(T_c<\infty\,|\,\bY_0=\by_0)<1$. 
If the state $\by_0\in\calS_0$, it implies that the irreducible class $\calS_0$ is transient since we have $\calA^C\cap\calS_0\ne\emptyset$; if $\by_0\notin\calS_0$ and $\calS_0$ is recurrent, it implies that, with a positive probability, the MMRRW $\{\bY_n\}$ starting from $\by_0\in\calS\setminus\calS_0$ will never reach the irreducible class $S_0$. As a result, $\{\bY_n\}$ is unstable in our sense.
\end{proof}

\begin{remark}
Theorem 2.1.9 of Fayolle et al.\ \cite{Fayolle95} requires that expression (\ref{eq:Zsigma_epsilon}) unconditionally holds for all $n\ge 0$. However, we are concerning the first entry time $T_c$ and that point does not become a problem. 
For example, consider a modified process of $\{\bY_n\}$ in which if $\bY_n$ enters $\calA^C$, then $\bY_{n+1}$ is forced to be the initial state $\by_0$ appearing in the proof above; in other words, the modified process restarts from the initial state after entering $\calA^C$. 
Denote the modified process by $\{\bY'_n\}$ and define $\{Z'_n\}$ by $Z'_n=f(\bY'_n)$. The time of first entry of $\{\bY'_n\}$ into $\calA^C$ is obviously equivalent to that of the original process and the process $\{Z'_n\}$ unconditionally satisfies expression (\ref{eq:Zsigma_epsilon}) for all $n\ge 0$. Applying Theorem 2.1.9 of Fayolle et al.\ \cite{Fayolle95} to $\{Z'_n\}$, we obtain the desired results. 
\end{remark}

%%%%%%%%%%%%%%%%%%%%%%%%
%
\subsection{Markov chains generated from $\calL$} \label{sec:generatedMC}

As mentioned above, one of the key notions we use is {\it induced Markov chain}. We define it in several steps. First, we consider Markov modulated partially-reflecting random walks obtained from the original MMRRW $\calL=\{(\bX_n,J_n)\}$ by removing some boundaries. We call those random walks expanded Markov chains. Next, the expanded Markov chains are expressed as Markov additive processes and each induced Markov chain is given by the background process of the corresponding Markov additive process. 
We also introduce a special embedded Markov chain, which is used for deriving an important property of the mean drift vectors in the next subsection.

\subsubsection{Expanded Markov chains}

For $l\in D$, let $\mathbb{B}^{x_l}$ be a boundary accompanied by $x_l$-axis and define it as 
\[
\mathbb{B}^{x_l} = \{\bx=(x(1),...,x(d))\in\mathbb{R}^d: x(l)=0\}, 
\]
which is the superplane perpendicular to $x_l$-axis and cuts $x_l$-axis at the origin. 
Roughly speaking, for nonempty $A\in\scrP(D)$, a Markov modulated partially-reflecting random walk with index $A$, denoted by $\hat{\calL}^A$, is a random walk obtained from $\calL$ by removing the reflecting barriers of $\calL$ on the boundaries $\mathbb{B}^{x_l},\,l\in A$.
We denote the process $\hat{\calL}^A$ by $\{\bY^A_n\}=\{(\bX^A_n,J^A_n)\}$, where $\{\bX^A_n\}$ is the random walk taking values in $\mathbb{Z}^d\cap\{\bx\in\mathbb{R}^d: \bx(D\setminus A)\ge 0\}$ and $\{J^A_n\}$ is the background process. Precisely, $\hat{\calL}^A$ is defined as follows.
Divide $\{\bx\in\mathbb{R}^d: \bx(D\setminus A)\ge 0\}$ into $2^{|D\setminus A|}$ subsets defined by 
\[
\hat{\calB}^{A,B}=\{\bx\in\mathbb{R}^d: \bx(B)>0,\ \bx(D\setminus(A\cup B))=0\},\ B\in\scrP(D\setminus A), 
\]
where we always have $A\cap B= \emptyset$. The state space of $\hat{\calL}^A$, denoted by $\hat{\calS}^A$ ({\it not} $\hat{S}^A$), is given by
\begin{equation}
\hat{\calS}^A = \bigcup_{B\in\scrP(D\setminus A)} (\hat{\calB}^{A,B}\cap\mathbb{Z}^d)\times S^{A\cup B}, 
\label{eq:hatcalSA}
\end{equation}
where $S^{A\cup B}$ is the set of background states that the original MMRRW $\calL=\{(\bX_n,J_n)\}$ takes when $\bX_n$ is in $\calB^{A\cup B}\cap\mathbb{Z}^d$ (see expression (\ref{eq:calS}) in Subsection \ref{sec:modeldescription}). 
For $n\ge 0$, let $\bZ^A_{n+1}$ be defined as $\bZ^A_{n+1}=\bX^A_{n+1}-\bX^A_n$. Let $B$ be an element of $\scrP(D\setminus A)$. If $\bX^A_n\in(\hat{\calB}^{A,B}\cap\mathbb{Z}^d)$, then $\bZ^A_{n+1}(A\cup B)\in\{-1,0,1\}^{|A\cup B|}$ and $\hat{\bZ}_{n+1}(D\setminus (A\cup B))\in\{0,1\}^{d-|A\cup B|}$ since we have $\bX^A_n(A)\in\mathbb{Z}^{|A|}$, $\bX^A_n(B)\ge 1$ and $\bX^A_n(D\setminus(A\cup B))=0$. 
For $\bx,\,\bz\in\mathbb{R}^d$, suppose that $\bx\in(\hat{\calB}^{A,B}\cap\mathbb{Z}^d)$, $\bz(A\cup B)\in\{-1,0,1\}^{|A\cup B|}$ and $\bz(D\setminus(A\cup B))\in\{0,1\}^{d-|A\cup B|}$. Further suppose that $\bx+\bz\in(\hat{\calB}^{A,C}\cap\mathbb{Z}^d)$ for some $C\in\scrP(D\setminus A)$ (the case where $C=B$ is included). 
Then, the transition probabilities of $\hat{\calL}^A$ are given by
\begin{equation}
\mathbb{P}(\bY^A_{n+1}=(\bx+\bz,j)\,|\,\bY^A_n=(\bx,i)) = p_{\bz}^{A\cup B,A\cup C}(i,j),\quad i\in S^{A\cup B},\ j\in S^{A\cup C}, \label{eq:tp_expand}
\end{equation}
where $p_{\bz}^{A\cup B,A\cup C}(i,j)$ is a transition probability of the original MMRRW $\calL$  (see expression (\ref{eq:tp_calL}) in Subsection \ref{sec:modeldescription}). 
We call $\hat{\calL}^A=\{(\bX^A_n,J^A_n)\}$ the expanded Markov chain of $\calL$ with index $A$. 

For nonempty $A\in\scrP(D)$, let a $d$-dimensional vector $\ba^A$ be defined as 
\begin{equation}
\ba^A = \lim_{n\to\infty} \frac{1}{n} \sum_{k=1}^n \bZ^A_k = \lim_{n\to\infty} \frac{1}{n} \sum_{k=1}^n \left( \bX^A_k - \bX^A_{k-1} \right), 
\label{eq:drift_vector}
\end{equation}
if the limit exists with probability one. We call $\ba^A$ the mean drift vector of the expanded Markov chain $\hat{\calL}^A$.

\subsubsection{Induced Markov chains}

Let $A$ be a nonempty element of $\scrP(D)$ and consider expanded Markov chain $\hat{\calL}^A=\{(\bX^A_n,J^A_n)\}$. Dividing $\bX^A_n$ into $\bX^A_n(A)$ and $\bX^A_n(D\setminus A)$ and putting the latter into the background state, then we obtain a $|A|$-dimensional Markov additive process, $\{(\bX^A_n(A),(\bX^A_n(D\setminus A),J^A_n))\}$, where $\bX^A_n(A)$ is the additive part taking values in $\mathbb{Z}^{|A|}$ and $(\bX^A_n(D\setminus A),J^A_n)$ is the background state. 
Obviously, the background process $\{(\bX^A_n(D\setminus A),J^A_n)\}$ is a Markov chain by itself, and the induced Markov chain of $\calL$ with index $A$, denoted by $\calL^A$, is given by the background process, i.e., $\calL^A = \{(\bX^A_n(D\setminus A),J^A_n)\}$. 
$\bX^A_n(D\setminus A)$ takes values in $\mathbb{Z}_+^{d-|A|}$, which is represented as $\mathbb{Z}_+^{d-|A|}=\{\bx(D\setminus A): \bx\in\mathbb{R}^d,\ \bx(D\setminus A)\ge 0\}\cap\mathbb{Z}^{d-|A|}$.
Divide $\{\bx(D\setminus A): \bx\in\mathbb{R}^d,\ \bx(D\setminus A)\ge 0\}$ into $2^{|D\setminus A|}$ subsets defined by (c.f.\ the definition of $\hat{\calB}^{A,B}$) 
\[
\calB^{A,B}=\{\bx(D\setminus A): \bx\in\mathbb{R}^d,\ \bx(B)>0,\ \bx(D\setminus(A\cup B))=0\},\ B\in\scrP(D\setminus A), 
\]
then the state space of $\calL^A$, denoted by $\calS^A$ ({\it not} $S^A$), is given by (c.f.\ the definition of $\hat{\calS}^A$)
\begin{equation}
\calS^A = \bigcup_{B\in\scrP(D\setminus A)} (\calB^{A,B}\cap\mathbb{Z}^{d-|A|})\times S^{A\cup B}. 
\label{eq:calSA}
\end{equation}
Let $B$ be an element of $\scrP(D\setminus A)$. For $\bx,\,\bz\in\mathbb{R}^d$, suppose that $\bx(D\setminus A)\in(\calB^{A,B}\cap\mathbb{Z}^{d-|A|})$, $\bz(B)\in\{-1,0,1\}^{|B|}$ and $\bz(D\setminus(A\cup B))\in\{0,1\}^{d-|A\cup B|}$. Further suppose that $\bx(D\setminus A)+\bz(D\setminus A)\in({\calB}^{A,C}\cap\mathbb{Z}^{d-|A|})$ for some $C\in\scrP(D\setminus A)$ (the case where $C=B$ is included). 
Then, the transition probabilities of $\calL^A$ are given by (c.f.\ expression (\ref{eq:tp_expand}))
\begin{align}
&\mathbb{P}((\bX^A_{n+1}(D\setminus A),J^A_{n+1})=(\bx(D\setminus A)+\bz(D\setminus A),j)\,|\,(\bX^A_n(D\setminus A),J^A_n)=(\bx(D\setminus A),i)) \cr
&\quad = \sum_{\bz(A)\in\{-1,0,1\}^{|A|}} p_{\bz}^{A\cup B,A\cup C}(i,j),\quad i\in S^{A\cup B},\ j\in S^{A\cup C}, 
 \label{eq:tp_induced}
\end{align}
where $\bz=(\bz(A),\bz(D\setminus A))$. 
Hence, we see that the induced Markov chain $\calL^A$ is also a $(d-|A|)$-dimensional MMRRW. 

\begin{remark} \label{re:inducedMC}
For nonempty $A\in\scrP(D)$, since induced Markov chain $\calL^A$ is a MMRRW, we can consider expanded Markov chains and induced Markov chains of $\calL^A$. 
Let $B$ be a nonempty element of $\scrP(D\setminus A)$. We denote by $\hat{\calL}^{A,B}$ the expanded Markov chain of $\calL^A$ with index $B$ and by $\calL^{A,B}$ the induced Markov chain of $\calL^A$ with index $B$. 
$\hat{\calL}^A$ is obtained from $\calL$ by removing the reflecting barriers on $\mathbb{B}^{x_l},\,l\in A$, and $\hat{\calL}^{A,B}$ is obtained from $\calL^A$ by removing the reflecting barriers on $\mathbb{B}^{x_l},\,l\in B$. Hence, we see that $\hat{\calL}^{A,B}$ is a subprocess of $\hat{\calL}^{A\cup B}$ and it is given as $\hat{\calL}^{A,B}=\{(\bX^{A\cup B}_n(D\setminus A),J^{A\cup B}_n)\}$. 
Since the induced Markov chain $\calL^{A,B}$ is the background process of the Markov additive process $\hat{\calL}^{A,B}=\{(\bX^{A\cup B}_n(B),(\bX^{A\cup B}_n(D\setminus (A\cup B)),J^{A\cup B}_n))\}$, it is given as $\calL^{A,B}=\{(\bX^{A\cup B}_n(D\setminus (A\cup B)),J^{A\cup B}_n)\}$. 
This implies that $\calL^{A,B}=\calL^{A\cup B}$; in other words, the induced Markov chain of $\calL^A$ with index $B$ is the induced Markov chain of $\calL$ with index $A\cup B$.
\end{remark}

Throughout the paper, we assume the following.
\begin{assumption} \label{as:inducedMC}
For all nonempty $A\in\scrP(D)$, induced Markov chain $\calL^A$ satisfies one of the following conditions. 
\begin{itemize}
\item[(i)] $\calL^A$ is semi-irreducible and stable in our sense.
\item[(ii)] $\calL^A$ is transient; this means that all the states of $\calL^A$ are transient but it is not necessary to be semi-irreducible.
\end{itemize}
\end{assumption}

Note that we will use statement (i) of Assumption \ref{as:inducedMC} for applying the ergodic theorem of Markov chains to our model. 
Let $\scrD_{stable}$ be the index set of stable induced Markov chains, i.e., 
\[
\scrD_{stable}=\{A\in\scrP(D): A\ne\emptyset,\ \mbox{$\calL^A$ is semi-irreducible and stable} \}. 
\]
Since $\calL^D$ is a finite Markov chain, it is semi-irreducible and stable in our sense, by Assumption \ref{as:inducedMC}; hence we always have $D\in\scrD_{stable}$.
For $A\in\scrD_{stable}$, $\calL^A$ has the unique stationary distribution, denoted by $\bpi^A=(\pi^A(\by),\,\by\in\calS^A)$, and the mean drift vector of $\hat{\calL}^A$, $\ba^A$, is given as 
\begin{equation}
\ba^A = \sum_{B\in\scrP(D\setminus A)}\ \ \sum_{(\bx,j)\in(\calB^{A,B}\cap\,\mathbb{Z}^{d-|A|})\times S^{A\cup B}} \balpha^{A\cup B}_j\,\pi^A((\bx,j)),  \label{eq:aA}
\end{equation}
where $\balpha^{A\cup B}_j$ is the conditional mean increment vector of $\calL$ defined in subsection \ref{sec:modeldescription}. Since $\bpi^A$ is the stationary distribution of $\calL^A$, we have $\ba^A(D\setminus A)=0$. 

\begin{remark}
Let $A$ be an element of $\scrP(D)$. By the definition of the mean drift vectors (expression (\ref{eq:drift_vector})), mean drift vector $\ba^A$ may exist even though the corresponding induced Markov chain $\calL^A$ is unstable in our sense. However, in our analysis, we only use the mean drift vectors whose corresponding induced Markov chains are stable in our sense. 
\end{remark}

\begin{example}[Two-dimensional QBD process] \label{ex:2D_QBD2} {\rm 
Consider a two-dimensional QBD process (two-dimensional MMRRW), $\calL=\{(\bX_n,J_n)\}$, described in Example \ref{ex:2D_QBD1}. 
$\calL$ has three induced Markov chains: $\calL^{\{1\}}$, $\calL^{\{2\}}$ and $\calL^D$, where $D=\{1,2\}$. 
$\calL^{\{1\}}=\{(X^{\{1\}}_n(2),J^{\{1\}}_n)\}$ and $\calL^{\{2\}}=\{(X^{\{2\}}_n(1),J^{\{2\}}_n)\}$ are QBD processes (one-dimensional MMRRWs), whose state spaces are respectively given by 
\[
\calS^{\{1\}}=(\{0\}\times S^{\{1\}}) \cup (\mathbb{N}\times S^D),\quad \calS^{\{2\}}=(\{0\}\times S^{\{2\}}) \cup (\mathbb{N}\times S^D). 
\]
We denote by $P^{\{1\}}$ the transition probability matrix of $\calL^{\{1\}}$ and by $P^{\{2\}}$ that of $\calL^{\{2\}}$; they are given as
\begin{align}
P^{\{1\}} = \begin{pmatrix} 
P_{(*,0)}^{\{1\},\{1\}} & P_{(*,1)}^{\{1\},D} & & & \cr
P_{(*,-1)}^{D,{\{1\}}} & P_{(*,0)}^{D,D} & P_{(*,1)}^{D,D} & & \cr
 & P_{(*,-1)}^{D,D} & P_{(*,0)}^{D,D} & P_{(*,1)}^{D,D}  &\cr
 & \ddots & \ddots & \ddots &
\label{eq:2DQBD_P1}
\end{pmatrix},
\\
P^{\{2\}} = \begin{pmatrix} 
P_{(0,*)}^{\{2\},\{2\}} & P_{(1,*)}^{\{2\},D} & & & \cr
P_{(-1,*)}^{D,{\{2\}}} & P_{(0,*)}^{D,D} & P_{(1,*)}^{D,D} & & \cr
 & P_{(-1,*)}^{D,D} & P_{(0,*)}^{D,D} & P_{(1,*)}^{D,D}  &\cr
 & \ddots & \ddots & \ddots &
\end{pmatrix},
\label{eq:2DQBD_P2}
\end{align}
where $P_{(*,z(2))}^{A,B}=\sum_{z(1)\in\{-1,0,1\}} P_{(z(1),z(2))}^{A,B}$ and $P_{(z(1),*)}^{A,B}=\sum_{z(2)\in\{-1,0,1\}} P_{(z(1),z(2))}^{A,B}$. 
$\calL^D=\{J^D_n\}$ is a finite Markov chain with state space $\calS^D=S^D$ and its transition probability matrix, denoted by $P^D$, is given as 
\begin{equation}
P^D = P_{(*,*)}^{D,D} = \sum_{\bz\in\{-1,0,1\}^2} P_{\bz}^{D,D}.
\label{eq:2DQBD_PD}
\end{equation}
\hfill$\Box$}
\end{example}

\subsubsection{Embedded Markov chain of the MMRRW $\calL$} \label{sec:embeddedMC} 

Another Markov chain crucial for our analysis is an embedded Markov chain of $\calL$ defined in this subsection.
For nonempty $A\in\scrP(D)$, let $K_A$ be a positive integer and assume that, for nonempty $A, B\in\scrP(D)$, if $|A|>|B|$ then $K_A < K_B$. 
We divide the state space $\calS$ into $2^d$ exclusive subsets, $\calV_A\subset\calS,\ A\in\scrP(D)$, which are inductively given as follows. 
\begin{description}
\item[Step 1] Set $k := d$ and $\calV_\emptyset := \calS$.
\item[Step 2] For $A\in\scrP(D)$ such that $|A|=k$, set $\calV_A := \calV_\emptyset \cap \{ (\bx,i)\in\calS : \bx(A)\ge K_A\}$. 
Set $\calV_\emptyset := \calV_\emptyset\setminus \left(\bigcup_{A\in\scrP(D),\,|A|=k} \calV_A \right) $ and $k := k-1$.
\item[Step 3] If $k>0$ then go to Step 2, otherwise we obtain $\calV_A$ for all $A\in\scrP(D)$.
\end{description}

\begin{figure}[htbp]
\begin{center}
\setlength{\unitlength}{0.7mm}
\begin{picture}(80,70)(0,0)
\thicklines
\put(0,10){\vector(1,0){70}}
\put(10,0){\vector(0,1){70}}
\put(70,5){\makebox(0,0){\normalsize $x_1$}}
\put(5,65){\makebox(0,0){\normalsize $x_2$}}
\thinlines
\put(30,30){\line(1,0){40}}
\put(30,30){\line(0,1){40}}
\multiput(10,30)(2,0){10}{\line(1,0){1}}
\multiput(30,10)(0,2){10}{\line(0,1){1}}
\put(2,30){\makebox(0,0){\normalsize $K_{\{1,2\}}$}}
\put(30,5){\makebox(0,0){\normalsize $K_{\{1,2\}}$}}
\put(50,50){\makebox(0,0){\normalsize $\calV_{\{1,2\}}$}}
\put(10,40){\line(1,0){20}}
\put(50,10){\line(0,1){20}}
\put(2,40){\makebox(0,0){\normalsize $K_{\{2\}}$}}
\put(50,5){\makebox(0,0){\normalsize $K_{\{1\}}$}}
\put(20,50){\makebox(0,0){\normalsize $\calV_{\{2\}}$}}
\put(60,20){\makebox(0,0){\normalsize $\calV_{\{1\}}$}}
\put(20,20){\makebox(0,0){\normalsize $\calV_\emptyset$}}
\end{picture}
\caption{Partition of the state space $\calS$ ($d=2$).}
\label{fig:partition2}
\end{center}
\end{figure}
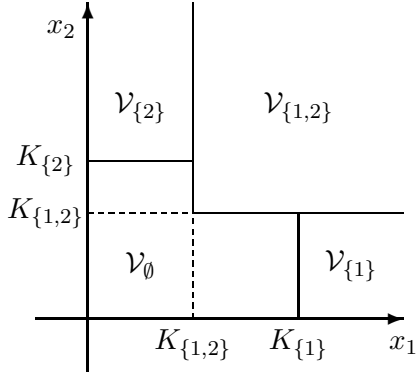

It is obvious that $\bigcup_{A\in\scrP(D)} \calV_A = \calS$ and that $\calV_A\cap\calV_B=\emptyset$ for all $A, B\in\scrP(D)$ such that $A\ne B$. Thus, $\{\calV_A : A\in\scrP(D)\}$ is a partition of $\calS$. We also see that $\calV_\emptyset$ is finite. 
An example when $d=2$ is described in Fig.\ \ref{fig:partition2}, where $\calS$ is divided into 4 subsets: $\calV_\emptyset$, $\calV_{\{1\}}$, $\calV_{\{2\}}$ and $\calV_{\{1,2\}}$. 
Note that, in Malyshev and Menshikov \cite{Malyshev81}, they picked up each part of the state space of a target RRW, where the part was denoted by $B^\Lambda_{ct}$, and used it for proving lemmas; on the other hand, we consider all the parts of the state space $\calS$ at the same time. 

Let a function $u$ on $\calS$ be defined as 
\[
u(\by) = \left\{ \begin{array}{ll}
u_A & \mbox{if $\by\in\calV_A$ for some {\it nonempty} $A\in\scrP(D)$}, \cr
1 & \mbox{otherwise},
\end{array} \right.
\]
where, for nonempty $A\in\scrP(D)$, $u_A$ is a positive integer; we assume that $u_A<K_A$. 
Furthermore, define a random sequence $\{\sigma_n\}$ as $\sigma_{n+1}=\sigma_n+u(\bY_{\sigma_n}),\, \sigma_0=0$, and define a Markov chain $\{\tilde{\bY}_n\}=\{(\tilde{\bX}_n,\tilde{J}_n)\}$ as $\tilde{\bY}_n = \bY_{\sigma_n},\,n\ge 0$. We call $\{\tilde{\bY}_n\}$ the embedded Markov chain of $\calL$ and denote it by $\tilde{\calL}$. 
By Proposition \ref{pr:Foster2}, if there exists a set of parameters $\{(K_A,u_A): A\in\scrP(D), A\ne\emptyset\}$ such that the condition of the proposition holds, then the MMRRW $\calL$ is stable in our sense. Furthermore, by Proposition \ref{pr:Markov_unstable2}, if there exists a set of parameters $\{(K_A,u_A): A\in\scrP(D), A\ne\emptyset\}$ such that the condition of Proposition \ref{pr:Markov_unstable2} holds, then $\calL$ is unstable in our sense. 
\begin{remark} \label{re:emexMC} 
Let $A$ be a nonempty element of $\scrP(D)$. If $\tilde{\bY}_k=\bY_{\sigma_k}=(\bX_{\sigma_k},J_{\sigma_k})\in\calV_A$, then we have $u(\bY_{\sigma_k})=u_A<K_A$ and $\bX_{\sigma_k}(A)\ge K_A$.
Since the process $\{\bX_n\}$ is skip free in all coordinates, we therefore see that $\{\bY_n\}=\{(\bX_n,J_n)\}$ does not touch the reflecting barriers of $\calL$ on the boundaries $\mathbb{B}^{x_l},\,l\in A$, during time interval $[\sigma_k,\sigma_{k+1}]$. Thus, in stochastic sense, the process $\{\bY_n\}$ behaves just like the expanded Markov chain $\hat{\calL}^A$ during that time interval.  
\end{remark}

\subsubsection{Embedded Markov chain of each expanded Markov chain} \label{sec:embeddedMC_expand} 

We also define an embedded Markov chain of each expanded Markov chain in the same manner. Let $A$ be a nonempty element of $\scrP(D)$ and consider expanded Markov chain $\hat{\calL}^A=\{\bY^A_n\}=\{(\bX^A_n,J^A_n)\}$. Using the same parameters, $K_{A\cup B},\,B\in\scrP(D\setminus A)$, as those used in the definition of $\{\tilde{\bY}_n\}$, the embedded Markov chain of $\calL$, we divide the state space $\hat{\calS}^A$ into $2^{|D\setminus A|}$ exclusive subsets, $\hat{\calV}^A_B\subset\hat{\calS}^A,\ B\in\scrP(D\setminus A)$, which are inductively given as follows. 
\begin{description}
\item[Step 1] Set $k := |D\setminus A|$ and $\hat{\calV}^A_\emptyset := \hat{\calS}^A$.
\item[Step 2] For $B\in\scrP(D\setminus A)$ such that $|B|=k$, set $\hat{\calV}^A_B := \hat{\calV}^A_\emptyset \cap \{ (\bx,i)\in\hat{\calS}^A : \bx(B)\ge K_{A\cup B}\}$. 
Set $\hat{\calV}^A_\emptyset := \hat{\calV}^A_\emptyset\setminus \left(\bigcup_{B\in\scrP(D\setminus A),\,|B|=k} \hat{\calV}^A_B \right) $ and $k := k-1$.
\item[Step 3] If $k>0$ then go to Step 2, otherwise we obtain $\hat{\calV}^A_B$ for all $B\in\scrP(D\setminus A)$.
\end{description}

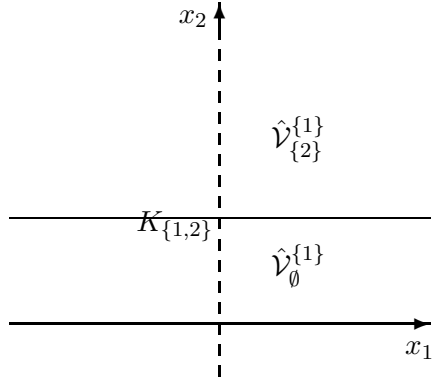
\begin{figure}[htbp]
\begin{center}
\setlength{\unitlength}{0.7mm}
\begin{picture}(85,70)(0,0)
\thicklines
\put(0,10){\vector(1,0){80}}
\put(78,5){\makebox(0,0){\normalsize $x_1$}}
\thicklines
\put(40,66){\vector(0,1){5}}
\multiput(40,0)(0,4){17}{\line(0,1){2}}
\put(35,68){\makebox(0,0){\normalsize $x_2$}}
\thinlines
\put(0,30){\line(1,0){80}}
\put(31.5,28.5){\makebox(0,0){\normalsize $K_{\{1,2\}}$}}
\put(55,45){\makebox(0,0){\normalsize $\hat{\calV}^{\{1\}}_{\{2\}}$}}
\put(55,21){\makebox(0,0){\normalsize $\hat{\calV}^{\{1\}}_\emptyset$}}
\end{picture}
\caption{Partition of the state space $\hat{\calS}^{\{1\}}$ ($d=2$).}
\label{fig:partition3}
\end{center}
\end{figure}

We have that $\bigcup_{B\in\scrP(D\setminus A)} \hat{\calV}^A_B = \hat{\calS}^A$ and $\hat{\calV}^A_B\cap\hat{\calV}^A_C=\emptyset$ for all $B, C\in\scrP(D\setminus A)$ such that $B\ne C$; hence, $\{\hat{\calV}^A_B : B\in\scrP(D\setminus A)\}$ is a partition of $\hat{\calS}^A$. 
From the definition of the partition of $\calS$ and that of the partition of $\hat{\calS}^A$, we see that $\calV_{A\cup B}\subset \hat{\calV}^A_B$ for $B\in\scrP(D\setminus A)$. 
An example when $d=2$ is described in Fig.\ \ref{fig:partition3}, where $\hat{\calS}^{\{1\}}$ is divided into 2 subsets: $\hat{\calV}^{\{1\}}_\emptyset$ and $\hat{\calV}^{\{1\}}_{\{2\}}$ (c.f.\ Fig.\ \ref{fig:partition2}). 

Using the same parameters, $u_{A\cup B},\,B\in\scrP(D\setminus A),\,B\ne\emptyset$, as those used in the definition of $\{\tilde{\bY}_n\}$, we define a function $u^A$ on $\hat{\calS}^A$ as 
\[
u^A(\by) = \left\{\begin{array}{ll}
u_{A\cup B} & \mbox{if $\by\in\hat{\calV}^A_B$ for some {\it nonempty} $B\in\scrP(D\setminus A)$}, \cr
1 & \mbox{otherwise}.
\end{array} \right.
\]
Note that, for nonempty $B\in\scrP(D\setminus A)$, if $\by\in\calV_{A\cup B}\subset\hat{\calV}^A_B$, then we have $u^A(\by)=u(\by)$, where $u$ is the function used in the definition of the embedded Markov chain of $\calL$.  
Further define a random sequence $\{\sigma^A_n\}$ as $\sigma^A_{n+1}=\sigma^A_n+u^A(\bY^A_{\sigma^A_n}),\, \sigma^A_0=0$, and define a Markov chain $\{\tilde{\bY}^A_n\}=\{(\tilde{\bX}^A_n,\tilde{J}^A_n)\}$ as $\tilde{\bY}^A_n = \bY^A_{\sigma^A_n},\,n\ge 0$. We call $\{\tilde{\bY}^A_n\}$ the embedded Markov chain of $\hat{\calL}^A$ and denote it by $\tilde{\calL}^A$. 
\begin{remark} \label{re:parameter_expand}
For nonempty $A\in\scrP(D)$, we only use parameters in $\{(K_B,u_B): B\in\scrP(D),\, |B|>|A|\}$ for defining the partition $\{\hat{\calV}^A_B: B\in\scrP(D\setminus A)\}$ and the random sequence $\{\sigma^A_n,\,n\ge 0\}$. 
\end{remark}

%%%%%%%%%%%%%%%%%%%%%%%
%
% Section 3, continuation
%
%%%%%%%%%%%%%%%%%%%%%%%
%
\subsection{Conditional mean increment vector of $\tilde{\calL}$} \label{sec:mean_increments}

Here we consider the conditional mean increment vectors of the embedded Markov chain $\tilde{\calL}=\{(\tilde{\bX}_n,\tilde{J}_n)\}$ defined in Subsection \ref{sec:embeddedMC} and obtain their approximation formulas represented in terms of $\ba^A,\,A\in\scrD_{stable}$, where $\ba^A$ is the mean drift vector of the expanded Markov chain $\hat{\calL}^A$ whose corresponding induced Markov chain $\calL^A$ is stable in our sense. The approximation formulas will be used for proving our main theorems.

For $\by=(\bx,j)\in\calS$, let $\tilde{\balpha}_{\by}$ be the conditional mean increment vector of $\tilde{\calL}$ given that the state is in $\by$, i.e., 
\[
\tilde{\balpha}_{\by} 
= \mathbb{E}(\tilde{\bX}_1-\tilde{\bX}_0\,|\,\tilde{\bY}_0=\by)
= \mathbb{E}(\bX_{u(\by)}-\bX_0\,|\,\bY_0=\by).
\]
For nonempty $A\in\scrP(D)$, let $\scrD^A_{stable}$ be defined as $\scrD^A_{stable} = \{B\in\scrP(D\setminus A): A\cup B\in\scrD_{stable} \}$. By Remark \ref{re:inducedMC}, the induced Markov chains of $\calL^A$ are induced Markov chains $\calL^{A\cup B},\,B\in\scrP(D\setminus A),\,B\ne\emptyset$; hence, $\scrD^A_{stable}$ is the index set of the stable induced Markov chains of $\calL^A$, including $\calL^A$ itself if it is stable in our sense. We have $D\setminus A\in\scrD^A_{stable}$ since $\calL^D$ is always stable. 
Approximation formulas for the conditional mean increment vectors of $\tilde{\calL}$ are given by the following proposition, which corresponds to Lemmas 2.1, 2.2 and 2.3 of Malyshev and Menshikov \cite{Malyshev81} (also see Lemmas 4.3.1, 4.3.2 and 4.3.3 of Fayolle et al.\ \cite{Fayolle95}). 
\begin{proposition} \label{pr:app_talpha}
For any $\varepsilon>0$ and any $\delta\in(0,1)$, there exist the set of positive integers $\{(K_A,u_A) : A\in\scrP(D),\,A\ne\emptyset\}$ satisfying the conditions mentioned in Subsection \ref{sec:embeddedMC} and the set of nonnegative vectors $\{(p^{A,B}_{\by},\,B\in\scrD^A_{stable}): A\in\scrP(D),\,A\ne\emptyset,\,\by\in\calV_A\}$ such that, for all nonempty $A\in\scrP(D)$ and for all $\by\in\calV_A$, 
\begin{align}
&\delta< \sum_{B\in\scrD^A_{stable}} p^{A,B}_{\by} \le 1,
\label{eq:app_pAB} \\
&\bigg| \tilde{\alpha}_{\by}(l)/u_A - \sum_{B\in\scrD^A_{stable}} p^{A,B}_{\by}\,a^{A\cup B}(l) \bigg| < \varepsilon\quad \mbox{for all $l\in D$}, 
\label{eq:app_talpha}
\end{align}
where if the induced Markov chain $\calL^A$ is stable, expression (\ref{eq:app_talpha}) is simplified and given by
\begin{equation}
\big| \tilde{\alpha}_{\by}(l)/u_A - a^A(l) \big| < \varepsilon\quad \mbox{for all $l\in D$}. 
\label{eq:app_talpha2}
\end{equation}
\end{proposition}

\medskip
In order to prove the proposition above, we propose several other propositions. 
For nonempty $A\in\scrP(D)$, consider expanded Markov chain $\hat{\calL}^A=\{\bY^A_n\}=\{(\bX^A_n,J^A_n)\}$ on the state space $\hat{\calS}_A$ and, for $\by\in\hat{\calS}^A$, define the expectation of the time-averaged increment vector of $\hat{\calL}^A$, denoted by $\bg^A_{\by}$, as 
\[
\bg^A_{\by}
= \mathbb{E}\bigg( \frac{1}{u_A} \sum_{n=1}^{u_A} (\bX^A_n-\bX^A_{n-1})\,\Big|\,\bY^A_0=\by \bigg)
= \frac{1}{u_A} \mathbb{E}(\bX^A_{u_A}-\bX^A_0\,|\,\bY^A_0=\by). 
\]
By Remark \ref{re:emexMC}, if $\bY_0=\bY^A_0=\by\in\calV_A\subset\hat{\calS}_A$, then we have $u(\by)=u_A$ and $\{\bY_n\}$ behaves just like $\{\bY^A_n\}$ during $[0,u_A]$, in stochastic sense; hence we immediately obtain the following proposition.
\begin{proposition} \label{pr:talpha_gA}
For $\by=(\bx,j)\in\calS$, if $\by\in\calV_A$ for some nonempty $A\in\scrP(D)$, then we have $\tilde{\balpha}_{\by}/u_A=\bg^A_{\by}$. 
\end{proposition}

An approximation formula of $\bg^A_{\by}$ is given by the following proposition.
\begin{proposition} \label{pr:app_gA_original}
For any $\varepsilon>0$ and any $\delta\in(0,1)$, there exist the set of positive integers $\{(K_A,u_A) : A\in\scrP(D),\,A\ne\emptyset\}$ satisfying the conditions mentioned in Subsection \ref{sec:embeddedMC} and the set of nonnegative vectors $\{(p^{A,B}_{\by},\,B\in\scrD^A_{stable}): A\in\scrP(D),\,A\ne\emptyset,\,\by\in\hat{\calV}^A_\emptyset\}$ such that, for all nonempty $A\in\scrP(D)$ and for all $\by\in\hat{\calV}^A_\emptyset$, 
\begin{align}
&\delta\le \delta^{(d-|A|+1)/d} < \sum_{B\in\scrD^A_{stable}} p^{A,B}_{\by} \le 1,
\label{eq:app_pAB_original} \\
&\bigg| g^A_{\by}(l) - \sum_{B\in\scrD^A_{stable}} p^{A,B}_{\by}\,a^{A\cup B}(l) \bigg| < \varepsilon(d-|A|+1)/d \le \varepsilon\quad \mbox{for all $l\in D$}, 
\label{eq:app_gA_original}
\end{align}
where if the induced Markov chain $\calL^A$ is stable, expression (\ref{eq:app_gA_original}) is simplified and given by
\begin{equation}
\big| g^A_{\by}(l) - a^A(l) \big| < \varepsilon(d-|A|+1)/d \le \varepsilon\quad \mbox{for all $l\in D$}. 
\label{eq:app_gA2_original}
\end{equation}
\end{proposition}

Proposition \ref{pr:app_talpha} is proved by Propositions \ref{pr:talpha_gA} and \ref{pr:app_gA_original}, as follows.

\begin{proof}[Proof of Proposition \ref{pr:app_talpha}] 
For nonempty $A\in\scrP(D)$, we have $\calV_A\subset\hat{\calV}^A_\emptyset$. Thus, by Propositions \ref{pr:talpha_gA}, $\tilde{\balpha}_{\by}/u_A=\bg^A_{\by}$ for all $\by\in\calV_A$. Thus, by expressions (\ref{eq:app_gA_original}) and (\ref{eq:app_gA2_original}) in Propositions \ref{pr:app_gA_original}, we obtain expressions (\ref{eq:app_talpha}) and (\ref{eq:app_talpha2}) in Proposition \ref{pr:app_talpha}. 
\end{proof}

%%%%%%%%%%%%%%%%%%%%%%%%%%%%%%%%%%%%%%%%%%%%%%%%%%%
%
Therefore, we focus on $\bg^A_{\by}$ and prove Proposition \ref{pr:app_gA_original}. By expressions (\ref{eq:hatcalSA}) and (\ref{eq:calSA}),  $\bg^A_{\by}$ is represented as
\begin{align}
\bg^A_{\by}
&= \frac{1}{u_A} \sum_{n=1}^{u_A}\ \sum_{\by'\in\hat{\calS}^A} \mathbb{E}(\bX^A_n-\bX^A_{n-1}\,|\,\bY^A_{n-1}=\by')\,\mathbb{P}(\bY^A_{n-1}=\by'\,|\,\bY^A_0=\by) \cr
&= \frac{1}{u_A} \sum_{n=1}^{u_A}\ \sum_{B\in\scrP(D\setminus A)}\ \sum_{\bx'\in\hat{\calB}^{A,B}\cap\mathbb{Z}^d}\ \sum_{j'\in S^{A\cup B}} \balpha^{A\cup B}_{j'}\,\mathbb{P}(\bY^A_{n-1}=(\bx',j')\,|\,\bY^A_0=\by) \cr
&= \frac{1}{u_A} \sum_{n=1}^{u_A}\ \sum_{B\in\scrP(D\setminus A)}\ \sum_{\bx'(D\setminus A)\in\calB^{A,B}\cap\mathbb{Z}^{d-|A|}}\ \sum_{j'\in S^{A\cup B}} \balpha^{A\cup B}_{j'} \cr 
&\qquad \cdot \sum_{\bx'(A)\in\mathbb{Z}^{|A|}} \mathbb{P}\big(\bX^A_{n-1}(A)=\bx'(A),(\bX^A_{n-1}(D\setminus A),J^A_{n-1})=(\bx'(D\setminus A),j')\,\big|\,\bY^A_0=\by\big) \cr
&= \sum_{B\in\scrP(D\setminus A)}\ \sum_{\bx'(D\setminus A)\in\calB^{A,B}\cap\mathbb{Z}^{d-|A|}}\ \sum_{j'\in S^{A\cup B}} \balpha^{A\cup B}_{j'} \cr 
&\qquad\qquad \cdot \frac{1}{u_A} \sum_{n=0}^{u_A-1} q_A^{(n)}\big((\bx(D\setminus A),j),(\bx'(D\setminus A),j')\big), 
\label{eq:gAy_induce}
\end{align}
where $\balpha^{A\cup B}_{j'}$ is the conditional mean increment vector of $\calL$ defined in Subsection \ref{sec:modeldescription} and we denote by $q_A^{(n)}((\bx(D\setminus A),j),(\bx'(D\setminus A),j'))$ the $n$-step transition probability of the induced Markov chain $\calL^A$, given by 
\begin{align*}
&q_A^{(n)}\big((\bx(D\setminus A),j),(\bx'(D\setminus A),j')\big) \cr
&\quad= 
\mathbb{P}\big((\bX^A_n(D\setminus A),J^A_n)=(\bx'(D\setminus A),j')\,\big|\,(\bX^A_0(D\setminus A),J^A_0)=(\bx(D\setminus A),j)\big).
\end{align*}
Hence, we obtain the following proposition. 
\begin{proposition} \label{pr:gA_induced}
For any nonempty $A\in\scrP(D)$ and any $\by\in\hat{\calS}^A$, $\bg^A_{\by}$ is given by the induced Markov chain $\calL^A$.
\end{proposition}

Next, for nonempty $A\in\scrP(D)$, we consider the embedded Markov chain of the expanded Markov chain $\hat{\calL}^A$, $\tilde{\calL}^A=\{\tilde{\bY}^A_n\}=\{(\tilde{\bX}^A_n,\tilde{J}^A_n)\}$, and the partition of the state space $\hat{\calS}^A$, $\{\hat{\calV}^A_B: B\in\scrP(D\setminus A)\}$. 
For $\by=(\bx,j)\in\hat{\calS}^A$, let $\tilde{\balpha}^A_{\by}$ be the conditional mean increment vector of $\tilde{\calL}^A$ given that the state is in $\by$, i.e., 
\[
\tilde{\balpha}^A_{\by} 
= \mathbb{E}(\tilde{\bX}^A_1-\tilde{\bX}^A_0\,|\,\tilde{\bY}^A_0=\by)
= \mathbb{E}(\bX^A_{u^A(\by)}-\bX^A_0\,|\,\bY^A_0=\by).
\]
For nonempty $B\in\scrP(D\setminus A)$, if $\by=(\bx,j)\in\hat{\calV}^A_B$, then we have $u^A(\by)=u_{A\cup B}$ and $\bx(B)\ge K_{A\cup B}>u_{A\cup B}$. Hence, by the same reason as that used in Remark \ref{re:emexMC}, the expanded Markov chain $\hat{\calL}^A$ behaves just like the expanded Markov chain $\hat{\calL}^{A\cup B}=\{\bY^{A\cup B}_n\}=\{(\bX^{A\cup B}_n,J^{A\cup B}_n)\}$ during time interval $[0,u_{A\cup B}]$, in stochastic sense, when $\bY^{A\cup B}_0=\bY^A_0=\by\in\hat{\calV}^A_B$. 
Thus, we have, for $\by\in\hat{\calV}^A_B$, 
\[
\mathbb{E}(\bX^A_{u_{A\cup B}}-\bX^A_0\,|\,\bY^A_0=\by) = \mathbb{E}(\bX^{A\cup B}_{u_{A\cup B}}-\bX^{A\cup B}_0\,|\,\bY^{A\cup B}_0=\by)
\]
and obtain the following proposition.
\begin{proposition} \label{pr:gA_expanded}
For nonempty $A\in\scrP(D)$ and nonempty $B\in\scrP(D\setminus A)$, we have, for $\by\in\hat{\calV}^A_B$, $\tilde{\balpha}^A_{\by}/u_{A\cup B}=\bg^{A\cup B}_{\by}$. 
\end{proposition}

From the definition of the partition of $\hat{\calS}^A$, it can be seen that if $\by=(\bx,j)\in\hat{\calV}^A_\emptyset$, then $\bx(A)\in\mathbb{Z}^{|A|}$ and $0\le \bx(D\setminus A)< \max_{|B|>|A|} K_B$. Further we have $\calV_A=\{(\bx,j)\in\hat{\calV}^A_\emptyset: \bx(A)\ge K_A\}$. 
Hence, we have $\{(\bx(D\setminus A),j): (\bx,j)\in\hat{\calV}^A_\emptyset\}=\{(\bx(D\setminus A),j): (\bx,j)\in\calV_A\}$ and both the sets are finite. Since $\{(\bx(D\setminus A),j): (\bx,j)\in\hat{\calV}^A_\emptyset\}$ as well as $\{(\bx(D\setminus A),j): (\bx,j)\in\calV_A\}$ is a set of states of the induced Markov chain $\calL^A$, we obtain, by Proposition \ref{pr:gA_induced}, the following proposition.
\begin{proposition} \label{pr:gA_finite}
For any nonempty $A\in\scrP(D)$, we have $\{\bg^A_{\by}: \by\in\hat{\calV}^A_\emptyset\}=\{\bg^A_{\by}: \by\in\calV_A\}$ and both the sets are finite. 
\end{proposition}

%%%%%%%%%%%%%%%%%%%%%%%%%%%%%%%%%%%%%%%%
%
We obtain the following approximation of $\bg^A_{\by}$ if the corresponding induced Markov chain $\calL^A$ is stable.
\begin{proposition} \label{pr:app_gA_stable}
Let $A$ be a nonempty element of $\scrP(D)$ and assume the induced Markov chain $\calL^A$ is stable in our sense. Furthermore, assume that the set of parameters $\{(K_B,u_B) : B\in\scrP(D),\,|B|>|A|\}$ is given. 
Then, for any $\varepsilon>0$, there exists a positive integer $u_A^*$ such that, for every $\by\in\hat{\calV}^A_\emptyset$, if $u_A\ge u_A^*$, then 
\begin{equation}
\big| g^A_{\by}(l) - a^A(l) \big| < \varepsilon\quad \mbox{for all $l\in D$}. 
\label{eq:app_gA_stable}
\end{equation}
\end{proposition}
\begin{proof}
First we note that, by Remark \ref{re:parameter_expand}, $\hat{\calV}^A_\emptyset$ is determined only using parameters in $\{(K_B,u_B): B\in\scrP(D),\, |B|>|A|\}$, which does not include $K_A$ and $u_A$. 
Since the induced Markov chain $\calL^A$ is stable, we obtain, by expression (\ref{eq:gAy_induce}) and the ergodic theorem of Markov chains, that for any $\by\in\hat{\calS}_A$, 
\[
\bg^A_{\by}\rightarrow\ba^A\quad\mbox{as}\quad u_A\rightarrow\infty,
\]
where we use expression (\ref{eq:aA}) of $\ba^A$. Thus, for any $\varepsilon>0$ and for any $\by\in\hat{\calV}^A_\emptyset\subset\hat{\calS}_A$, there exists a positive integer $u_A^*$ such that if $u_A\ge u_A^*$, then $|g^A_{\by}(l)-a^A(l)|<\varepsilon$ for all $l\in D$. Since, by Proposition \ref{pr:gA_finite}, $\{\bg^A_{\by}: \by\in\hat{\calV}^A_\emptyset\}$ is finite, we can commonly set this $u_A^*$ at the same value for every $\by\in\hat{\calV}^A_\emptyset$. 
\end{proof}

If the induced Markov chain $\calL^A$ is unstable, $\bg^A_{\by}$ is approximated as follows. 
\begin{proposition} \label{pr:app_gA_unstable}
Let $A$ be a nonempty element of $\scrP(D)$ and assume the induced Markov chain $\calL^A$ is unstable in our sense, which means that every state of $\calL^A$ is transient. Let the set of parameters $\{(K_B,u_B) : B\in\scrP(D),\,|B|>|A|\}$ be given. 
Then, for any $\varepsilon>0$, there exists a positive integer $u_A^*$ such that, for every $\by\in\hat{\calV}^A_\emptyset$, if $u_A\ge u_A^*$, then, for all $l\in D$,  
\begin{align}
\Big| g^A_{\by}(l)-\sum_{\genfrac{}{}{0pt}{}{B\in\scrP(D\setminus A)}{B\ne\emptyset}} \sum_{\by'\in\hat{\calV}^A_B} g^{A\cup B}_{\by'}(l) \cdot \frac{u_{A\cup B}}{u_A} \sum_{k=1}^{u_A} \mathbb{P}\big(\sigma^A_k\le u_A,\,\bY^A_{\sigma^A_{k-1}}=\by'\,\big|\,\bY^A_0=\by\big) \Big| < \varepsilon. 
\label{eq:app_gA_unstable}
\end{align}
\end{proposition}

\begin{proof}
Fix the value of $\varepsilon$ at an arbitrary positive number.
For $\by\in\hat{\calV}^A_\emptyset$, we obtain, by the definition of $\bg^A_{\by}$, 
\begin{align}
\bg^A_{\by}
&= \frac{1}{u_A} \mathbb{E}\Big(\sum_{n=0}^{u_A} 1(\sigma^A_n\le u_A<\sigma^A_{n+1}) (\bX^A_{u_A}-\bX^A_0)\,\Big|\,\bY^A_0=\by\Big) \cr
&= \frac{1}{u_A} \mathbb{E}\Big(\sum_{n=0}^{u_A} 1(\sigma^A_n\le u_A<\sigma^A_{n+1}) \Big( \sum_{k=1}^n (\bX^A_{\sigma^A_k}-\bX^A_{\sigma^A_{k-1}})+(\bX^A_{u_A}-\bX^A_{\sigma^A_n}) \Big)\,\Big|\,\bY^A_0=\by\Big) \cr
&= \frac{1}{u_A} \sum_{k=1}^{u_A} \mathbb{E}\big(1(\sigma^A_k\le u_A) (\bX^A_{\sigma^A_k}-\bX^A_{\sigma^A_{k-1}})\,\big|\,\bY^A_0=\by\big) \cr
&\qquad\qquad + \frac{1}{u_A} \sum_{n=0}^{u_A} \mathbb{E}\big(1(\sigma^A_n\le u_A<\sigma^A_{n+1})(\bX^A_{u_A}-\bX^A_{\sigma^A_n})\,\big|\,\bY^A_0=\by\big), 
\end{align}
where we use the fact that $\sigma^A_n>u_A$ for $n>u_A$; $1(\cdot)$ is an indicator function. By strong Markov property and Proposition \ref{pr:gA_expanded}, we obtain 
\begin{align}
&\mathbb{E}\big(1(\sigma^A_k\le u_A) (\bX^A_{\sigma^A_k}-\bX^A_{\sigma^A_{k-1}})\,\big|\,\bY^A_0=\by\big) \cr
&\quad= \sum_{\by'\in\hat{\calS}^A} \mathbb{E}\big(\bX^A_{\sigma^A_{k-1}+u^A(\by')}-\bX^A_{\sigma^A_{k-1}}\,\big|\,\sigma^A_{k-1}+u^A(\by')\le u_A,\,\bY^A_{\sigma^A_{k-1}}=\by',\,\bY^A_0=\by\big) \cr
&\qquad\qquad\qquad \cdot \mathbb{P}(\sigma^A_k\le u_A,\,\bY^A_{\sigma^A_{k-1}}=\by'\,|\,\bY^A_0=\by) \cr
&\quad= \sum_{\by'\in\hat{\calS}^A} \mathbb{E}\big(\bX^A_{u^A(\by')}-\bX^A_0\,\big|\,\bY^A_0=\by'\big)\,\mathbb{P}\big(\sigma^A_k\le u_A,\,\bY^A_{\sigma^A_{k-1}}=\by'\,\big|\,\bY^A_0=\by\big) \cr
%
%&\quad= \sum_{\by'\in\hat{\calS}^A} \tilde{\balpha}^A_{\by'}\,\mathbb{P}\big(\sigma^A_k\le u_A,\,\bY^A_{\sigma^A_{k-1}}=\by'\,\big|\,\bY^A_0=\by\big) \cr
%
&\quad= \sum_{\by'\in\hat{\calV}^A_\emptyset} \mathbb{E}\big(\bX^A_1-\bX^A_0\,\big|\,\bY^A_0=\by'\big)\,\mathbb{P}\big(\sigma^A_k\le u_A,\,\bY^A_{\sigma^A_{k-1}}=\by'\,\big|\,\bY^A_0=\by\big) \cr
&\qquad\qquad +\sum_{\genfrac{}{}{0pt}{}{B\in\scrP(D\setminus A)}{B\ne\emptyset}} \sum_{\by'\in\hat{\calV}^A_B} u_{A\cup B}\,\bg^{A\cup B}_{\by'}\,\mathbb{P}\big(\sigma^A_k\le u_A,\,\bY^A_{\sigma^A_{k-1}}=\by'\,\big|\,\bY^A_0=\by\big).
\end{align}
Thus, $\bg^A_{\by}$ is represented as 
\begin{equation}
\bg^A_{\by} = \bphi^A_{1,\by} + \bphi^A_{2,\by} + \bphi^A_{3,\by}, 
\label{eq:gA_phiA123}
\end{equation}
where the vectors $\bphi^A_{1,\by}$, $\bphi^A_{2,\by}$ and $\bphi^A_{3,\by}$ are defined as 
\begin{align*}
&\bphi^A_{1,\by} = \frac{1}{u_A} \sum_{k=1}^{u_A} \sum_{\by'\in\hat{\calV}^A_\emptyset} \mathbb{E}\big(\bX^A_1-\bX^A_0\,\big|\,\bY^A_0=\by'\big)\,\mathbb{P}\big(\sigma^A_k\le u_A,\,\bY^A_{\sigma^A_{k-1}}=\by'\,\big|\,\bY^A_0=\by\big), \cr
&\bphi^A_{2,\by} = \frac{1}{u_A} \sum_{k=1}^{u_A}\sum_{\genfrac{}{}{0pt}{}{B\in\scrP(D\setminus A)}{B\ne\emptyset}} \sum_{\by'\in\hat{\calV}^A_B} u_{A\cup B}\,\bg^{A\cup B}_{\by'}\,\mathbb{P}\big(\sigma^A_k\le u_A,\,\bY^A_{\sigma^A_{k-1}}=\by'\,\big|\,\bY^A_0=\by\big), \cr
&\bphi^A_{3,\by} = \frac{1}{u_A} \sum_{n=0}^{u_A} \mathbb{E}\big(1(\sigma^A_n\le u_A<\sigma^A_{n+1})(\bX^A_{u_A}-\bX^A_{\sigma^A_n})\,\big|\,\bY^A_0=\by\big). 
\end{align*}

First, we consider $\bphi^A_{1,\by}$. Temporally define $\calV_A^*$ as $\calV_A^*=\{(\bx'(D\setminus A),j'): (\bx',j')\in\hat{\calV}^A_\emptyset\}$; this $\calV_A^*$ is finite (see the derivation of Proposition \ref{pr:gA_finite}).  Since the process $\{\bX^A_n\}$ is skip free in all coordinates, we obtain, for $l\in D$, 
\begin{align}
|\phi^A_{1,\by}(l)| 
&\le \frac{1}{u_A} \sum_{k=1}^{u_A}\ \sum_{\by'\in\hat{\calV}^A_\emptyset} \mathbb{P}\big(\sigma^A_{k-1}+u^A(\by')\le u_A,\,\bY^A_{\sigma^A_{k-1}}=\by'\,\big|\,\bY^A_0=\by\big) \cr
&\le \frac{1}{u_A} \sum_{k=1}^{u_A} \mathbb{P}\big(\sigma^A_{k-1}\le u_A-1,\,\bY^A_{\sigma^A_{k-1}}\in\hat{\calV}^A_\emptyset\,\big|\,\bY^A_0=\by\big) \cr
&\le \frac{1}{u_A} \sum_{k=0}^{u_A-1} \mathbb{P}\big(\bY^A_k\in\hat{\calV}^A_\emptyset\,\big|\,\bY^A_0=\by\big) \cr
&= \frac{1}{u_A} \sum_{k=0}^{u_A-1} \mathbb{P}\big((\bX^A_k(D\setminus A),J^A_k)\in\calV_A^*\,\big|\,(\bX^A_0(D\setminus A),J^A_0)=(\bx(D\setminus A),j)\big), 
\label{eq:phiA1_g^A}
\end{align}
where $\{(\bX^A_n(D\setminus A),J^A_n)\}$ is the induced Markov chain $\calL^A$. Since, by the assumption, all the states of $\calL^A$ are transient and $\calV_A^*$ is finite, there exists a positive integer $u_{A,1}^*$ such that if $u_A\ge u_{A,1}^*$, then $|\phi^A_{1,\by}(l)|<\varepsilon/2$ for all $l\in D$; we can commonly give this $u_{A,1}^*$ for all $\by\in\hat{\calV}^A_\emptyset$ since if $\by=(\bx,j)\in\hat{\calV}^A_\emptyset$, then $(\bx(D\setminus A),j)\in\calV_A^*$ and $\calV_A^*$ is finite. 

Next, we consider $\bphi^A_{3,\by}$. Since $\sigma^A_{n+1}-\sigma^A_n\le\max_{|B|>|A|} u_B$ for all $n\ge 0$ and $\{\bX^A_n\}$ is skip free in all coordinates, we have, for all $\by=(\bx,j)\in\hat{\calV}^A_\emptyset$ and for all $l\in D$, 
\begin{align}
|\phi^A_{3,\by}(l)| 
&\le \frac{1}{u_A} \sum_{n=0}^{u_A} \mathbb{E}\big(1(\sigma^A_n\le u_A<\sigma^A_{n+1}) | X^A_{u_A}(l)-X^A_{\sigma^A_n}(l)| \,\big|\,\bY^A_0=\by\big) \cr
&\le \frac{1}{u_A} \sum_{n=0}^{u_A} \mathbb{E}\big(1(\sigma^A_n\le u_A<\sigma^A_{n+1}) (\sigma^A_{n+1}-\sigma^A_n) \,\big|\,\bY^A_0=\by\big) \cr
&\le \frac{\max_{|B|>|A|} u_B}{u_A} \mathbb{E}\Big(\sum_{n=0}^{u_A} 1(\sigma^A_n\le u_A<\sigma^A_{n+1}) \,\Big|\,\bY^A_0=\by\Big) \cr
&= \frac{\max_{|B|>|A|} u_B}{u_A}.
\label{eq:phiA2_g^A} 
\end{align}
Thus, there exists a positive integer $u_{A,2}^*$ such that, for all $\by\in\hat{\calV}^A_\emptyset$, if $u_A\ge u_{A,2}^*$, then $|\phi^A_{3,\by}(l)|<\varepsilon/2$ for all $l\in D$.
As a result, letting $u_A^*$ be set as $u^*_A=\max\{u_{A,1}^*,\,u_{A,2}^*\}$, we obtain from equation (\ref{eq:gA_phiA123}) that, for all $\by\in\hat{\calV}^A_\emptyset$, if $u_A\ge u_A^*$, then 
\begin{align}
|g^A_{\by}(l)-\phi^A_{2,\by}(l)| \le |\phi^A_{1,\by}(l)| + |\phi^A_{3,\by}(l)| < \varepsilon\quad\mbox{for all $l\in D$}, 
\end{align}
and this completes the proof.
\end{proof}

%%%%%%%%%%%%%%%%%%%%%%%%%%%%%%%
%
\begin{proof}[Proof of Proposition \ref{pr:app_gA_original}]
In order to prove the proposition, we use induction with respect to the cardinality of nonempty elements in $\scrP(D)$, where we use $k$ for the parameter of induction. 

Fix the value of $\varepsilon$ at an arbitrary  positive number and that of $\delta$ at a real number in $(0,1)$. 
First, we consider the case where $k=d$. For nonempty $A\in\scrP(D)$, if $|A|=d$, then we have $A=D$ and, by Assumption \ref{as:inducedMC}, induced Markov chain $\calL^D$ is stable. 
Thus, by Proposition \ref{pr:app_gA_stable}, there exists a positive integer $u_D$ such that, for every $\by\in\hat{\calV}^D_\emptyset$, 
\[
|g^D_{\by}(l)-p^{D,\emptyset}_{\by}\,a^D(l)| < \varepsilon/d \quad\mbox{for all $l\in D$}, 
\]
where we set $p^{D,\emptyset}_{\by}=1$. Since we have $\scrD^D_{stable}=\{\emptyset\}$, this implies 
\[
\delta^{1/d}< \sum_{B\in\scrD^D_{stable}} p^{D,B}_{\by} = p^{D,\emptyset}_{\by}\le 1.
\]
We set $u_D$ at such a value and set $K_D$ at a positive integer satisfying $u_D<K_D$.

Next, for $k$ such that $1<k\le d$, suppose we have the set of parameters $\{(K_A,u_A) : A\in\scrP(D),\,|A|\ge k\}$ satisfying the conditions mentioned in Subsection \ref{sec:embeddedMC} and the set of nonnegative vectors $\{(p^{A,B}_{\by},\,B\in\scrD^A_{stable}): A\in\scrP(D),\,|A|\ge k,\,\by\in\hat{\calV}^A_\emptyset\}$ such that, for every $A\in\scrP(D),\,|A|\ge k$, and for every $\by\in\hat{\calV}^A_\emptyset$, 
\begin{align*}
&\delta^{(d-|A|+1)/d}< \sum_{B\in\scrD^A_{stable}} p^{A,B}_{\by} \le 1, \\
&\bigg| g^A_{\by}(l) - \sum_{B\in\scrD^A_{stable}} p^{A,B}_{\by}\,a^{A\cup B}(l) \bigg| < \varepsilon(d-|A|+1)/d \quad \mbox{for all $l\in D$}. 
\end{align*}

Consider the case of $k-1$. Let $A$ be an element of $\scrP(D)$ satisfying $|A|=k-1$. 
If induced Markov chain $\calL^A$ is stable, then, by Proposition \ref{pr:app_gA_stable}, there exists a positive integer $u_A$  such that, for every $\by\in\hat{\calV}^A_\emptyset$, 
\[
\biggl| g^A_{\by}(l)-\sum_{B\in\scrD^A_{stable}} p^{A,B}_{\by}\,a^{A\cup B}(l) \biggr| < \varepsilon(d-k)/d \quad\mbox{for all $l\in D$}, 
\]
where we set $p^{A,\emptyset}_{\by}=1$ and $p^{A,B}_{\by}=0$ for nonempty $B\in\scrD^A_{stable}$; this implies 
\[
\delta^{(d-k)/d}< \sum_{B\in\scrD^A_{stable}} p^{A,B}_{\by}\le 1.
\]
We set $u_A$ at such a value and set $K_A$ at a positive integer satisfying $u_A<K_A$ and $K_A>K_{A\cup B}$ for every nonempty $B\in\scrP(D\setminus A)$.

When the induced Markov chain $\calL^A$ is unstable, we use Proposition \ref{pr:app_gA_unstable} for obtaining the desired approximation formula of $\bg^A_{\by}$. We give, for nonempty $B\in\scrD^A_{stable}$ and for $\by\in\hat{\calV}^A_\emptyset$, $p^{A,B}_{\by}$ as 
\begin{align}
p^{A,B}_{\by} &= \sum_{\genfrac{}{}{0pt}{}{C\in\scrP(D\setminus A)}{C\ne\emptyset}} \sum_{C'\in\scrD^{A\cup C}_{stable}} 1(C\cup C'=B) \sum_{\by'\in\hat{\calV}^A_C} p^{A\cup C,C'}_{\by'} q^{A,C}_{\by,\by'},
\end{align}
where $q^{A,C}_{\by,\by'}$ is defined as
\[
q^{A,C}_{\by,\by'} = \frac{u_{A\cup C}}{u_A} \sum_{k=1}^{u_A} \mathbb{P}\big(\sigma^A_k\le u_A,\,\bY^A_{\sigma^A_{k-1}}=\by'\,\big|\,\bY^A_0=\by\big). 
\]
Note that, for $C\in\scrP(D\setminus A)$, $C\ne\emptyset$ and for $C'\in\scrD^{A\cup C}_{stable}$, $p^{A\cup C,C'}_{\by}$ is given by the assumption of induction. 
%
%%%%%%%%%%%%%%%%%%%
Let us consider the sum of $q^{A,B}_{\by,\by'}$ with respect to $B$ and $\by'$. 
From the definition of function $u^A$, if $\bY^A_{\sigma^A_{k-1}}\in\hat{\calV}^A_\emptyset$, then we have $u^A(\bY^A_{\sigma^A_{k-1}})=1$; for nonempty $B\in\scrP(D\setminus A)$, if $\bY^A_{\sigma^A_{k-1}}\in\hat{\calV}^A_B$, then we have $u^A(\bY^A_{\sigma^A_{k-1}})=u_{A\cup B}$. Thus, we obtain 
\begin{align}
\sum_{\genfrac{}{}{0pt}{}{B\in\scrP(D\setminus A)}{B\ne\emptyset}} \sum_{\by'\in\hat{\calV}^A_B} q^{A,B}_{\by,\by'}
&= \sum_{\genfrac{}{}{0pt}{}{B\in\scrP(D\setminus A)}{B\ne\emptyset}} \sum_{\by'\in\hat{\calV}^A_B} \frac{u_{A\cup B}}{u_A} \sum_{k=1}^{u_A} \mathbb{P}\big(\sigma^A_k\le u_A,\,\bY^A_{\sigma^A_{k-1}}=\by'\,\big|\,\bY^A_0=\by\big) \cr
&= \frac{1}{u_A} \sum_{\genfrac{}{}{0pt}{}{B\in\scrP(D\setminus A)}{B\ne\emptyset}} \sum_{k=1}^{u_A} \mathbb{E}\big(u^A(\bY^A_{\sigma^A_{k-1}})\,1(\sigma^A_k\le u_A)\,1(\bY^A_{\sigma^A_{k-1}}\in\hat{\calV}^A_B)\,\big|\,\bY^A_0=\by\big) \cr
&= \psi^A_{1,\by} - \psi^A_{2,\by}, 
\label{eq:qAB_psi12}
\end{align}
where 
\begin{align}
&\psi^A_{1,\by} = \frac{1}{u_A} \mathbb{E}\Big( \sum_{k=1}^{u_A} u^A(\bY^A_{\sigma^A_{k-1}})\,1(\sigma^A_k\le u_A)\,\big|\,\bY^A_0=\by\Big), \\
&\psi^A_{2,\by} = \frac{1}{u_A} \sum_{k=1}^{u_A} \mathbb{P}\big(\sigma^A_k\le u_A,\,\bY^A_{\sigma^A_{k-1}}\in\hat{\calV}^A_\emptyset\,\big|\,\bY^A_0=\by\big). 
\end{align}
Since $\sigma^A_k=\sum_{l=1}^{k} u^A(\bY^A_{\sigma^A_{l-1}})$, we have
\[
\sum_{k=1}^{u_A} u^A(\bY^A_{\sigma^A_{k-1}})\,1(\sigma^A_k\le u_A) = \sum_{k=1}^{u_A} \sigma^A_k\,1(\sigma^A_k\le u_A<\sigma^A_{k+1}), 
\] 
and this leads us to
\begin{align}
&\frac{u_A-\max_{l\in D\setminus A} u_{A\cup\{l\}}}{u_A} < \psi^A_{1,\by} \le u_A/u_A = 1.
\end{align}
Thus, for any $\varepsilon_0>0$ satisfying $\delta^{1/d}+\varepsilon_0\in(0,1)$ and for any $\by\in\hat{\calV}^A_\emptyset$, there exists a positive integer $u^*_{A,1a}$ such that, if $u_A\ge u^*_{A,1a}$, then $\delta^{1/d}+\varepsilon_0< \psi^A_{1,\by} \le 1$. 
Let $\varepsilon_0$ be fixed at some value satisfying the condition above. 
Setting $\calV^*_A=\{(\bx'(D\setminus A),j'): (\bx',j')\in\hat{\calV}^A_\emptyset\}$, which is finite (see the proof of Proposition \ref{pr:app_gA_unstable}), we obtain
\begin{align}
\psi^A_{2,\by}
&=\frac{1}{u_A} \sum_{k=1}^{u_A} \mathbb{P}\big(\sigma^A_{k-1}+1\le u_A,\,\bY^A_{\sigma^A_{k-1}}\in\hat{\calV}^A_\emptyset\,\big|\,\bY^A_0=\by\big) \cr
&\le \frac{1}{u_A} \sum_{k=0}^{u_A-1} \mathbb{P}\big(\bY^A_k\in\hat{\calV}^A_\emptyset\,\big|\,\bY^A_0=\by\big) \cr
&= \frac{1}{u_A} \sum_{k=0}^{u_A-1} \mathbb{P}\big((\bX^A_k(D\setminus A),J^A_k)\in\calV_A^*\,\big|\,(\bX^A_0(D\setminus A),J^A_0)=(\bx(D\setminus A),j)\big), 
\end{align}
where $\{(\bX^A_n(D\setminus A),J^A_n)\}$ is the induced Markov chain $\calL^A$. 
Note that, by the assumption, all the states of $\calL^A$ are transient and $\calV_A^*$ is finite. Thus, for the fixed $\varepsilon_0>0$, there exists a positive integer $u_{A,1b}^*$ such that if $u_A\ge u_{A,1b}^*$, then $\psi^A_{2,\by}<\varepsilon_0$; we can commonly give this $u_{A,1b}^*$ for all $\by=(\bx,j)\in\hat{\calV}^A_\emptyset$ since if $\by=(\bx,j)\in\hat{\calV}^A_\emptyset$, then we have $(\bx(D\setminus A),j)\in\calV_A^*$ and $\calV_A^*$ is finite. 
As a result, by expression (\ref{eq:qAB_psi12}), if $u_A\ge\max\{u^*_{A,1a},u^*_{A,1b}\}$, then, for every $\by\in\hat{\calV}^A_\emptyset$, 
\begin{equation}
\delta^{1/d}< \sum_{\genfrac{}{}{0pt}{}{B\in\scrP(D\setminus A)}{B\ne\emptyset}} \sum_{\by'\in\hat{\calV}^A_B} q^{A,B}_{\by,\by'} \le 1;
\label{eq:app_sum_qAB}
\end{equation}
in other words, there exists a positive integer $u^*_{A,1}$ such that, for every $\by\in\hat{\calV}^A_\emptyset$, if $u_A\ge u^*_{A,1}$, then expression (\ref{eq:app_sum_qAB}) holds. 
Using this result, we evaluate the sum of $p^{A,B}_{\by}$ with respect to $B$. 
%
%%%%%%%%%%%%%%%%%%%%%%%%%%%%%
Since $\{A\cup B: B\in\scrD^A_{stable}\}=\{A\cup C\cup C': C\in\scrP(D\setminus A),\,C'\in\scrD^{A\cup C}_{stable}\}$, we have $\sum_{B\in\scrD^A_{stable}} 1(C\cup C'=B)=1$ for $C\in\scrP(D\setminus A)$ and $C'\in\scrD^{A\cup C}_{stable}$. Thus, we obtain, for $\by\in\hat{\calV}^A_\emptyset$, 
\begin{align}
&\sum_{B\in\scrD^A_{stable}} p^{A,B}_{\by}
= \sum_{\genfrac{}{}{0pt}{}{C\in\scrP(D\setminus A)}{C\ne\emptyset}} \sum_{\by'\in\hat{\calV}^A_C}\ \sum_{C'\in\scrD^{A\cup C}_{stable}} p^{A\cup C,C'}_{\by'} q^{A,C}_{\by,\by'}, 
\label{eq:sum_pAB}
\end{align}
and, by the assumption of induction and expression (\ref{eq:app_sum_qAB}), we obtain, for every $\by\in\hat{\calV}^A_\emptyset$, if $u_A\ge u^*_{A.1}$, then 
\begin{align}
&\sum_{B\in\scrD^A_{stable}} p^{A,B}_{\by}
\le \sum_{\genfrac{}{}{0pt}{}{C\in\scrP(D\setminus A)}{C\ne\emptyset}} \sum_{\by'\in\hat{\calV}^A_C} q^{A,C}_{\by,\by'}
\le 1, \\
&\sum_{B\in\scrD^A_{stable}} p^{A,B}_{\by}
> \delta^{(d-k+1)/d} \sum_{\genfrac{}{}{0pt}{}{C\in\scrP(D\setminus A)}{C\ne\emptyset}} \sum_{\by'\in\hat{\calV}^A_C} q^{A,C}_{\by,\by'} 
> \delta^{(d-k+1)/d}\cdot\delta^{1/d}
> \delta^{(d-|A|+1)/d}, 
\end{align}
where $|A|=k-1$ and we use the fact that $\hat{\calV}^A_C\subset \hat{\calV}^{A\cup C}_\emptyset$ and $|A\cup C|\ge k$. 

%%%%%%%%%%%%%%%%%%%%%%%%
In order to obtain an approximation formula of $\bg^A_{\by}$, we consider the following: 
\begin{align}
\bg^A_{\by} - \sum_{B\in\scrD^A_{stable}} p^{A,B}_{\by} \ba^{A\cup B}
= \bphi^A_{1,\by}+\bphi^A_{2,\by}, 
\end{align}
where 
\begin{align*}
&\bphi^A_{1,\by} = \bg^A_{\by} - \sum_{\genfrac{}{}{0pt}{}{B\in\scrP(D\setminus A)}{B\ne\emptyset}} \sum_{\by'\in\hat{\calV}^A_B} \bg^{A\cup B}_{\by'}\,q^{A,B}_{\by,\by'}, \\
&\bphi^A_{2,\by} = \sum_{\genfrac{}{}{0pt}{}{B\in\scrP(D\setminus A)}{B\ne\emptyset}} \sum_{\by'\in\hat{\calV}^A_B} \bg^{A\cup B}_{\by'}\,q^{A,B}_{\by,\by'} - \sum_{B\in\scrD^A_{stable}} p^{A,B}_{\by} \ba^{A\cup B}.
\end{align*}
By Proposition \ref{pr:app_gA_unstable}, there exists a positive integer $u^*_{A,2}$ such that, for every $\by\in\hat{\calV}^A_\emptyset$, if $u_A\ge u^*_{A,2}$, then, for all $l\in D$, 
\begin{align}
&|\phi^A_{1,\by}(l)| < \varepsilon/d.
\end{align}
By the same reason as that used for expression (\ref{eq:sum_pAB}), we obtain, for $\by\in\hat{\calV}^A_\emptyset$, 
\begin{align}
&\sum_{B\in\scrD^A_{stable}} p^{A,B}_{\by} \ba^{A\cup B} 
= \sum_{\genfrac{}{}{0pt}{}{C\in\scrP(D\setminus A)}{C\ne\emptyset}} \sum_{\by'\in\hat{\calV}^A_C}\ \sum_{C'\in\scrD^{A\cup C}_{stable}} p^{A\cup C,C'}_{\by'} q^{A,C}_{\by,\by'}\,\ba^{A\cup C\cup C'}, 
\label{eq:sum_pABaAB}
\end{align}
and, by the assumption of induction and expressions (\ref{eq:app_sum_qAB}) and (\ref{eq:sum_pABaAB}), we obtain, for all $l\in D$,  
\begin{align}
|\phi^A_{2,\by}(l)| 
&= \Bigl| \sum_{\genfrac{}{}{0pt}{}{B\in\scrP(D\setminus A)}{B\ne\emptyset}} \sum_{\by'\in\hat{\calV}^A_B} g^{A\cup B}_{\by'}(l)\,q^{A,B}_{\by,\by'} - \sum_{\genfrac{}{}{0pt}{}{B\in\scrP(D\setminus A)}{B\ne\emptyset}} \sum_{\by'\in\hat{\calV}^A_B}\ \sum_{C\in\scrD^{A\cup B}_{stable}} p^{A\cup B,C}_{\by'} a^{A\cup B\cup C}(l)\,q^{A,B}_{\by,\by'} \Bigr| \nonumber \\
&\le \sum_{\genfrac{}{}{0pt}{}{B\in\scrP(D\setminus A)}{B\ne\emptyset}} \sum_{\by'\in\hat{\calV}^A_B} \Bigl| g^{A\cup B}_{\by'}(l) - \sum_{C\in\scrD^{A\cup B}_{stable}} p^{A\cup B,C}_{\by'} a^{A\cup B,C}(l) \Bigr| q^{A,B}_{\by,\by'} \cr
&< \frac{\varepsilon (d-k+1)}{d} \sum_{\genfrac{}{}{0pt}{}{B\in\scrP(D\setminus A)}{B\ne\emptyset}} \sum_{\by'\in\hat{\calV}^A_B} q^{A,B}_{\by,\by'} \cr
&\le \frac{\varepsilon (d-k+1)}{d},
\end{align}
where we use the fact that $\hat{\calV}^A_B\subset \hat{\calV}^{A\cup B}_\emptyset$ and $|A\cup B|\ge k$. 
Thus, for every $\by\in\hat{\calV}^A_\emptyset$, if $u_A\ge u^*_{A,2}$, then, for $l\in D$, 
\begin{align}
&\Big| g^A_{\by}(l) - \sum_{B\in\scrD^A_{stable}} p^{A,B}_{\by} a^{A\cup B}(l)\,\Big| 
< \frac{\varepsilon}{d} + \frac{\varepsilon (d-k+1)}{d} 
= \frac{\varepsilon (d-|A|+1)}{d}, 
\end{align}
where $|A|=k-1$. 
We set $u_A$ so that it satisfies $u_A\ge\max\{u^*_{A,1},u^*_{A,2}\}$ and set $K_A$ at a positive integer satisfying $u_A<K_A$ and $K_A>K_{A\cup B}$ for every nonempty $B\in\scrP(D\setminus A)$. Then, the assumption of induction holds for the case where $|A|=k-1$, and this completes the proof.
\end{proof}

%%%%%%%%%%%%%%%%%%%%%%%
%
% Section 3, continuation
%
%%%%%%%%%%%%%%%%%%%%%%%
%
\subsection{Stability and instability conditions}

\subsubsection{Stability condition for semi-irreducible MMRRWs}

We consider a stability condition for a semi-irreducible MMRRW $\calL=\{\bY_n\}=\{(\bX_n,J_n)\}$ and give it as a condition that the mean drift vectors satisfy. For the purpose, we use a linear function on $\mathbb{R}^d$ as a test function and apply Proposition \ref{pr:Foster2}. Linear functions are simple but, in our experience, they work well, especially for low dimensional cases. The stability condition we get can easily be verified if the mean drift vectors are obtained. 
Hereafter, we denote by $\langle\bx_1,\bx_2\rangle$ the inner product of vectors $\bx_1,\bx_2\in\mathbb{R}^d$. The linear function we use is given by $\langle\bx,\bw\rangle$, where $\bw$ is a $d$-dimensional positive vector. The following theorem corresponds to Condition B and Theorem 2.1 of Malyshev and Menshikov \cite{Malyshev81} (also see Condition B and Theorem 4.3.4 of Fayolle et al.\ \cite{Fayolle95}). 
\begin{theorem}[Stable MMRRW] \label{th:stable}
If there exists a positive vector $\bw\in\mathbb{R}^d$ such that $\langle \ba^A,\bw \rangle <0$ for all $A\in\scrD_{stable}$, then the skip-free semi-irreducible MMRRW $\calL$ is stable in our sense.
\end{theorem}
\begin{proof}
We prove this theorem by Proposition \ref{pr:Foster2}. 
Let $\bw\in\mathbb{R}^d$ be a positive vector satisfying the condition of the theorem and consider the following linear function as a test function:
\begin{equation}
f(\by) = f((\bx,j)) = \langle \bx,\bw \rangle,\quad \by\in\calS, 
\end{equation}
where $\calS$ is the state space of the MMRRW $\calL$. 
We consider the embedded Markov chain of $\calL$, $\tilde{\calL}=\{\tilde{\bY}_n\}=\{(\tilde{\bX}_n,\tilde{J}_n)\}$, and evaluate the conditional mean increment of $f(\tilde{\bY}_n)$. Since $\langle\bx,\bw\rangle$ is linear with respect to $\bx$, we have, for $\by\in\calS$, 
\begin{align}
\mathbb{E}(f(\tilde{\bY}_{n+1})-f(\tilde{\bY}_n)\,|\,\tilde{\bY}_n=\by) 
&= \mathbb{E}(\langle\tilde{\bX}_{n+1},\bw\rangle-\langle\tilde{\bX}_n,\bw\rangle\,|\,\tilde{\bY}_n=\by) \cr
&= \langle \mathbb{E}(\tilde{\bX}_{n+1}-\tilde{\bX}_n\,|\,\tilde{\bY}_n=\by),\bw\rangle \cr
&= \langle\tilde{\balpha}_{\by},\bw\rangle, 
\end{align}
where $\tilde{\balpha}_{\by}$ is the conditional mean increment of $\tilde{\calL}$. 
Let $\varepsilon^*$ be defined as $\varepsilon^*=\min_{A\in\scrD_{stable}} |\langle\ba^A,\bw\rangle|>0$ and set $\varepsilon_0>0$ and $\delta_0\in(0,1)$ so that they satisfy $\delta_0\,\varepsilon^*-\varepsilon_0 \sum_{l=1}^d w(l) >0$; we denote by $\varepsilon$ the left hand side of the inequality, i.e., $\varepsilon=\delta_0\,\varepsilon^*-\varepsilon_0 \sum_{l=1}^d w(l)$.
By Proposition \ref{pr:app_talpha}, for the $\varepsilon_0$ and $\delta_0$, there exist $\{(K_A,u_A) : A\in\scrP(D),\,A\ne\emptyset\}$ and $\{(p^{A,B}_{\by},\,B\in\scrD^A_{stable}): A\in\scrP(D),\,A\ne\emptyset,\,\by\in\calV_A\}$ such that, for all nonempty $A\in\scrP(D)$ and for all $\by\in\calV_A$, 
\begin{align*}
&\delta_0< \sum_{B\in\scrD^A_{stable}} p^{A,B}_{\by} \le 1, \\
&\bigg| \tilde{\alpha}_{\by}(l)/u_A - \sum_{B\in\scrD^A_{stable}} p^{A,B}_{\by}\,a^{A\cup B}(l) \bigg| < \varepsilon_0\quad \mbox{for all $l\in D$}. 
\end{align*}
Thus, we obtain, for all nonempty $A\in\scrP(D)$ and for all $\by\in\calV_A$, 
\begin{align}
\langle\tilde{\balpha}_{\by},\bw\rangle
&= u_A \Bigl\langle \sum_{B\in\scrD^A_{stable}} p^{A,B}_{\by}\,\ba^{A\cup B},\bw \Bigr\rangle + u_A \Bigl\langle \tilde{\balpha}_{\by}/u_A-\sum_{B\in\scrD^A_{stable}} p^{A,B}_{\by}\,\ba^{A\cup B},\bw \Bigr\rangle \cr
&= u_A \sum_{B\in\scrD^A_{stable}} p^{A,B}_{\by} \langle\ba^{A\cup B},\bw\rangle + u_A \sum_{l=1}^d \Bigl(\tilde{\alpha}_{\by}(l)/u_A-\sum_{B\in\scrD^A_{stable}} p^{A,B}_{\by}\,a^{A\cup B}(l)\Bigr) w(l) \cr
&\le -u_A \Bigl(\delta_0\,\varepsilon^* - \varepsilon_0 \sum_{l=1}^d w(l) \Bigr) \cr
&= -u_A \varepsilon.
\end{align}
By the definition of the embedded Markov chain of the MMRRW $\calL$, if $\by\in\calV_\emptyset$, then $u(\by)=1$ and 
\begin{align}
\mathbb{E}\big(f(\tilde{\bY}_{n+1})-f(\tilde{\bY}_n)\,|\,\tilde{\bY}_n=\by\big) 
=\mathbb{E}\big(\langle\bX_{\sigma_n+1}-\bX_{\sigma_n},\bw\rangle\,|\,\bY_{\sigma_n}=\by\big)
\le \sum_{l=1}^d w(l)<\infty, 
\end{align}
where we use the fact that the MMRRW is skip free in all coordinates.
$\calV_\emptyset$ is finite and we have $\calS\setminus\calV_\emptyset=\sum_{\genfrac{}{}{0pt}{}{A\in\scrP(D) }{A\ne\emptyset}} \calV_A$. Thus, by Proposition \ref{pr:Foster2}, we see that the semi-irreducible MMRRW is stable in our sense. 
\end{proof}

\begin{example}[Two-dimensional QBD process; see Examples \ref{ex:2D_QBD1} and \ref{ex:2D_QBD2}] \label{ex:2D_QBD3} {\rm 
For the two-dimensional QBD process in Examples \ref{ex:2D_QBD1}, consider the case where $\ba^D<0$, $a^{\{1\}}(1)<0$ and $a^{\{2\}}(2)<0$. In this case, the induced Markov chains $\calL^{\{1\}}$ and $\calL^{\{2\}}$ as well as $\calL^D$ are stable, i.e., $\scrD_{stable}=\{\{1\},\{2\},D\}$, where $D=\{1,2\}$, and $a^{\{1\}}(2)$ and $a^{\{2\}}(1)$ are zero. 
Let a two-dimensional vector $\bw$ be set as $\bw=(1,1)>0$, then we have
\[
\langle\ba^{\{1\}},\bw\rangle = a^{\{1\}}(1)<0,\quad 
\langle\ba^{\{2\}},\bw\rangle = a^{\{2\}}(2)<0,\quad 
\langle\ba^D,\bw\rangle = a^D(1)+a^D(2)<0. 
\]
Hence, by Theorem \ref{th:stable}, the two-dimensional QBD process is stable in our sense. 
\hfill$\Box$
} \end{example}

%%%%%%%%%%%%%%%%%%%%%%%%%%%%%%%%%%%%%%%%%%%%%%%
%
\subsubsection{Some extensions of the results in Subsection \ref{sec:mean_increments}}

In order to obtain an instability condition for semi-irreducible MMRRWs, we here extend the results in Subsection \ref{sec:mean_increments}. 

First, we redefine the set of parameters $\{(K_A,u_A): A\in\scrP(D),\,A\ne\emptyset\}$, which is used for dividing the state space $\calS$ (resp.\ $\hat{\calS}^A$) into mutually exclusive subsets, $\calV_A,\,A\in\scrP(D)$ (resp.\ $\hat{\calV}^A_B,\,B\in\scrP(D\setminus A)$), where $\calS$ is the state space of the MMRRW $\calL=\{(\bX_n,J_n)\}$ and, for nonempty $A\in\scrP(D)$, $\hat{\calS}^A$ is that of the expanded Markov chain $\hat{\calL}^A=\{(\bX^A_n,J^A_n)\}$. 
For the parameter set, we so far assumed that $K_A>u_A$ for nonempty $A\in\scrP(D)$ and that, for nonempty $A,B\in\scrP(D)$, if $|A|<|B|$, then $K_A>K_B$; this assumption leads us to 
\[
K_A > \max\big\{ \max_{\genfrac{}{}{0pt}{}{B\in\scrP(D)}{|B|>|A|}} K_B,\ u_A\big\}\quad\mbox{for all nonempty $A\in\scrP(D)$}.
\]
Instead, we hereafter assume the following condition: 
\begin{equation}
K_A > \max_{\genfrac{}{}{0pt}{}{B\in\scrP(D)}{|B|>|A|}} K_B + u_A\quad\mbox{for all nonempty $A\in\scrP(D)$}.
\label{eq:Ku_condition2}
\end{equation}
Let $A$ be a nonempty element of $\scrP(D)$. Since the process $\{\bX^A_n\}$ is skip free in all coordinates, it can be seen from the new assumption that if $(\bX^A_0,J^A_0)\in\calV_A\subset\hat{\calV}^A_\emptyset$, then we have, for $n\in\{1,2,...,u_A\}$ and for $l\in A$, 
\[
X^A_n(l) \ge X^A_0(l)-u_A \ge K_A-u_A > \max_{\genfrac{}{}{0pt}{}{B\in\scrP(D)}{|B|>|A|}} K_B, 
\]
where we use the fact that if $(\bx,j)\in\calV_A$, then $\bx(A)\ge K_A$. Thus, for $n\in\{1,2,...,u_A\}$ and for nonempty $B\in\scrP(D\setminus A)$, if $(\bX^A_0,J^A_0)\in\calV_A$ and $(\bX^A_n,J^A_n)\in\hat{\calV}^A_B\subset\hat{\calV}^{A\cup B}_\emptyset$, then we have $\bX^A_n(A)>K_{A\cup B}$ and $(\bX^A_n,J^A_n)\in\calV_{A\cup B}\subset\hat{\calV}^A_B$, where we use the fact that if $(\bx,j)\in\hat{\calV}^A_B$, then $\bx(B)\ge K_{A\cup B}$.
As a result, we obtain the following proposition. 
\begin{proposition} \label{pr:XAn_in_VA}
Let $A$ be a nonempty element of $\scrP(D)$ and consider expanded Markov chain $\hat{\calL}^A=\{\bY^A_n\}=\{(\bX^A_n,J^A_n)\}$. Under condition (\ref{eq:Ku_condition2}), if $\bY^A_0\in\calV_A\subset\hat{\calV}^A_\emptyset$, then 
\[
\bY^A_n\in\hat{\calV}^A_\emptyset \cup \biggl(\bigcup_{\genfrac{}{}{0pt}{}{B\in\scrP(D\setminus A)}{B\ne\emptyset}} \calV_{A\cup B}\biggr)\quad\mbox{for all $n\in\{1,2,...,u_A\}$}. 
\]
\end{proposition}

%%%%%%%%%%%%%%%%%%%%
%
Next, we extend Propositions \ref{pr:app_gA_original} and \ref{pr:app_gA_unstable}. For the purpose, we consider a real function $f$ on $\mathbb{R}^d$ and assume that, for some constant $C>0$, $f$ satisfies
\begin{equation}
|f(\bx')-f(\bx)|\le C\, \|\bx'-\bx\|_1\quad\mbox{for $\bx,\bx'\in\mathbb{R}^d$}, 
\label{eq:f_condition}
\end{equation}
where $\|\bx\|_1=\sum_{l=1}^d |x(l)|$. Further, for nonempty $A\in\scrP(D)$ and for $\by=(\bx,j)\in\hat{\calS}^A$, we define $g^A_{f,\by}$ as 
\[
g^A_{f,\by} 
=\frac{1}{u_A} \sum_{k=1}^{u_A} \mathbb{E}(f(\bX^A_k)-f(\bX^A_{k-1})\,|\,\bY^A_0=\by) 
=\frac{1}{u_A} \mathbb{E}\big(f(\bX^A_{u_A})-f(\bX^A_0)\,\big|\,\bY^A_0=\by\big). 
\]
This $g^A_{f,\by}$ satisfies properties similar to those for $\bg^A_{\by}$, described in Subsection \ref{sec:mean_increments}. 
Consider Proposition \ref{pr:app_gA_unstable}. It can be seen from the proof of the proposition that the reason why the proposition holds comes from the fact that the process $\{\bX^A_n\}$ is skip free in all coordinates; in other words, $X^A_{n+1}(l)-X^A_n(l)$ is bounded for all $l\in D$; by the assumption for $f$, it also holds for $g^A_{f,\by}$. Thus, we obtain the following corollary of Proposition \ref{pr:app_gA_unstable}.
%
%%%%%%%%%%%%%%%%%%%
\begin{corollary} \label{co:app_gAf_unstable}
Let $f$ be a real function on $\mathbb{R}^d$ satisfying inequality (\ref{eq:f_condition}). Let $A$ be a nonempty element of $\scrP(D)$ and assume induced Markov chain $\calL^A$ is unstable in our sense. Furthermore, assume that the set of parameters $\{(K_B,u_B) : B\in\scrP(D),\,|B|>|A|\}$ is given and it satisfies condition (\ref{eq:Ku_condition2}). 
Then, for any $\varepsilon>0$, there exists a positive integer $u_{f,A}^*$ such that, for every $\by\in\calV_A\subset\hat{\calV}^A_\emptyset$, if $u_A\ge u_{f,A}^*$, then 
\begin{align}
\Big| g^A_{f,\by}-\sum_{\genfrac{}{}{0pt}{}{B\in\scrP(D\setminus A)}{B\ne\emptyset}} \sum_{\by'\in\calV_{A\cup B}} g^{A\cup B}_{f,\by'} q^{A,B}_{\by,\by'} \Big| < \varepsilon, 
\label{eq:app_gAf_unstable}
\end{align}
where $q^{A,B}_{\by,\by'}$ has already appeared in the proof of Proposition \ref{pr:app_gA_original} and it is given as 
\[
q^{A,B}_{\by,\by'}=\frac{u_{A\cup B}}{u_A} \sum_{k=1}^{u_A} \mathbb{P}\big(\sigma^A_k\le u_A,\,\bY^A_{\sigma^A_{k-1}}=\by'\,\big|\,\bY^A_0=\by\big).
\]
\end{corollary}

\begin{proof}
This corollary can be proved in the same manner as that used in the proof of Proposition \ref{pr:app_gA_unstable}; hence, we only outline the proof. 

Replacing $\bX^A_n$ with $f(\bX^A_n)$, we obtain from expression (\ref{eq:gA_phiA123}) the following: 
\begin{equation}
g^A_{f,\by} = \phi^A_{1,\by} + \phi^A_{2,\by} + \phi^A_{3,\by}, 
\label{eq:gAf_phiA123}
\end{equation}
where $\phi^A_{1,\by}$, $\phi^A_{2,\by}$ and $\phi^A_{3,\by}$ are defined as 
\begin{align*}
&\phi^A_{1,\by} = \frac{1}{u_A} \sum_{k=1}^{u_A} \sum_{\by'\in\hat{\calV}^A_\emptyset} \mathbb{E}\big(f(\bX^A_1)-f(\bX^A_0)\,\big|\,\bY^A_0=\by'\big)\,\mathbb{P}\big(\sigma^A_k\le u_A,\,\bY^A_{\sigma^A_{k-1}}=\by'\,\big|\,\bY^A_0=\by\big), \cr
&\phi^A_{2,\by} 
%= \frac{1}{u_A} \sum_{k=1}^{u_A}\sum_{\genfrac{}{}{0pt}{}{B\in\scrP(D\setminus A)}{B\ne\emptyset}} \sum_{\by'\in\hat{\calV}^A_B} u_{A\cup B}\,g^{A\cup B}_{f,\by'}\,\mathbb{P}\big(\sigma^A_k\le u_A,\,\bY^A_{\sigma^A_{k-1}}=\by'\,\big|\,\bY^A_0=\by\big), \cr
= \sum_{\genfrac{}{}{0pt}{}{B\in\scrP(D\setminus A)}{B\ne\emptyset}} \sum_{\by'\in\hat{\calV}^A_B} g^{A\cup B}_{f,\by'} q^{A,B}_{\by,\by'}, \cr
&\phi^A_{3,\by} = \frac{1}{u_A} \sum_{n=0}^{u_A} \mathbb{E}\big(1(\sigma^A_n\le u_A<\sigma^A_{n+1})(f(\bX^A_{u_A})-f(\bX^A_{\sigma^A_n}))\,\big|\,\bY^A_0=\by\big).
\end{align*}
From the assumption for $f$ and the fact that the process $\{\bX^A_n\}$ is skip free in all coordinates, we obtain 
\begin{align}
&|f(\bX^A_1)-f(\bX^A_0)| \le C\|\bX^A_1-\bX^A_0\|_1 \le C d,\\
&|f(\bX^A_{u_A})-f(\bX^A_{\sigma^A_n})| 
\le C\|\bX^A_{u_A}-\bX^A_{\sigma^A_n}\| 
\le C d (\sigma^A_{n+1}-\sigma^A_n)
\le C d \max_{\genfrac{}{}{0pt}{}{B\in\scrP(D)}{|B|>|A|}} u_B, 
\end{align}
where $\sigma^A_n$ is assumed to satisfy $\sigma^A_n\le u_A<\sigma^A_{n+1}$. 
Thus, by the same arguments as those used for deriving expressions (\ref{eq:phiA1_g^A}) and (\ref{eq:phiA2_g^A}), we see that there exists positive integer $u_{A}^*$ such that if $u_A\ge u_{A}^*$, then $|\phi^A_{1,\by}|<\varepsilon/2$ and $|\phi^A_{3,\by}|<\varepsilon/2$ and we can commonly give this $u_{A}^*$ for all $\by\in\hat{\calV}^A_\emptyset$. 
As a result, for all $\by\in\hat{\calV}^A_\emptyset$, if $u_A\ge u_{A}^*$, then 
\begin{align}
|g^A_{f,\by}-\phi^A_{2,\by}| \le |\phi^A_{1,\by}| + |\phi^A_{3,\by}| < \varepsilon. 
\end{align}
On the other hand, by Proposition \ref{pr:XAn_in_VA}, if $\by\in\calV_A\subset\hat{\calV}^A_\emptyset$, then, for every $\sigma^A_k\le u_A$,  
\[
\bY^A_{\sigma_k} \notin \bigcup_{\genfrac{}{}{0pt}{}{B\in\scrP(D\setminus A)}{B\ne\emptyset}} \bigl(\hat{\calV}^A_B\setminus \calV_{A\cup B}\bigr), 
\]
where $\calV_{A\cup B}\subset\hat{\calV}^A_B$. Thus, for nonempty $B\in\scrP(D\setminus A)$ and for $\by'\in\hat{\calV}^A_B\setminus \calV_{A\cup B}$, we have $q^{A,B}_{\by,\by'}=0$, and this leads us to expression (\ref{eq:app_gAf_unstable}).
\end{proof}

By Corollary \ref{co:app_gAf_unstable}, we obtain the following corollary of Proposition \ref{pr:app_gA_original}. 
%
%%%%%%%%%%%%%%%%%%
\begin{corollary} \label{co:app_gAf_original}
For nonempty $A\in\scrP(D)$, let a positive integers $u'_A$ be given. For any $\varepsilon>0$ and any $\delta\in(0,1)$, there exist the set of positive integers $\{(K_A,u_A) : A\in\scrP(D),\,A\ne\emptyset\}$ satisfying condition (\ref{eq:Ku_condition2}) as well as condition $u_A\ge u'_A$ for all nonempty $A\in\scrP(D)$ and the set of nonnegative vectors $\{(p^{A,B}_{\by,\by'},\,B\in\scrD^A_{stable},\by'\in\calV_{A\cup B}): A\in\scrP(D),\,A\ne\emptyset,\,\by\in\calV_A\}$ such that, for all nonempty $A\in\scrP(D)$ and for all $\by\in\calV_A$, 
\begin{align}
&\delta\le \delta^{(d-|A|+1)/d} < \sum_{B\in\scrD^A_{stable}}\,\sum_{\by'\in\calV_{A\cup B}} p^{A,B}_{\by,\by'} \le 1,
\label{eq:app_pABf_original} \\
&\bigg| g^A_{f,\by} - \sum_{B\in\scrD^A_{stable}}\,\sum_{\by'\in\calV_{A\cup B}} p^{A,B}_{\by,\by'}\,g^{A\cup B}_{f,\by'} \bigg| < \varepsilon(d-|A|+1)/d \le \varepsilon. 
\label{eq:app_gAf_original}
\end{align}
\end{corollary}

\begin{proof}
This corollary can be proved in a manner similar to that used in Proposition \ref{pr:app_gA_original}, as follows. 
Arbitrarily fix the value of $\varepsilon$ at a positive number and that of $\delta$ at a real number in $(0,1)$. Further, let $A$ be a nonempty element of $\scrP(D)$. We use induction with respect to the cardinality of nonempty elements in $\scrP(D)$, where $k$ is used for the parameter of induction. 
Before doing it, we note that if the induced Markov chain $\calL^A$ is stable, then we have $\emptyset\in\scrD^A_{stable}$; thus, by setting $p^{A,\emptyset}_{\by,\by}=1$ and other $p^{A,B}_{\by,\by'}=0$, we have, for $\by\in\calV_A$,  
\[
\sum_{B\in\scrD^A_{stable}}\,\sum_{\by'\in\calV_{A\cup B}} p^{A,B}_{\by,\by'}\,g^{A\cup B}_{f,\by'} = g^A_{f,\by}, 
\]
and inequalities (\ref{eq:app_pABf_original}) and (\ref{eq:app_gAf_original}) obviously hold irrespective of the value of $u_A$. Hence, we can set $u_A$ at any positive number satisfying $u_A\ge u'_A$. 

First, we consider the case where $k=d$. If $|A|=d$, then we have $A=D$ and, by Assumption \ref{as:inducedMC}, the induced Markov chain $\calL^D$ is stable. Thus, inequalities (\ref{eq:app_pABf_original}) and (\ref{eq:app_gAf_original}) hold and we set $u_D$ at a positive integer satisfying $u_D\ge u'_D$ and set $K_D$ at a positive integer satisfying $K_D>u_D$.
Next, for $k$ such that $1<k\le d$, suppose that we have the set of parameters $\{(K_A,u_A) : A\in\scrP(D),\,|A|\ge k\}$ and the set of nonnegative vectors $\{(p^{A,B}_{\by,\by'},\,B\in\scrD^A_{stable},\,\by'\in\calV_{A\cup B}): A\in\scrP(D),\,|A|\ge k,\,\by\in\calV_A\}$ satisfying the statements of the corollary. 

Consider the case of $k-1$ and let $A$ be a nonempty element of $\scrP(D)$ satisfying $|A|=k-1$. If the induced Markov chain $\calL^A$ is stable, then inequalities (\ref{eq:app_pABf_original}) and (\ref{eq:app_gAf_original}) hold; we set $u_A$ at a positive integer satisfying $u_A\ge u'_A$ and set $K_A$ at a positive integer satisfying $K_A>\max_{|B|>|A|}K_B+u_A$.
If the induced Markov chain $\calL^A$ is unstable, we use Corollary \ref{co:app_gAf_unstable} for obtaining the desired approximation formula of $g^A_{f,\by}$. 
For nonempty $B\in\scrD^A_{stable}$, $\by\in\calV_A\subset\hat{\calV}^A_\emptyset$ and $\by'\in\calV_{A\cup B}\subset\hat{\calV}^A_B$, we give $p^{A,B}_{\by,\by'}$ as 
\begin{align}
p^{A,B}_{\by,\by'} &= \sum_{\genfrac{}{}{0pt}{}{C\in\scrP(D\setminus A)}{C\ne\emptyset}} \sum_{C'\in\scrD^{A\cup C}_{stable}} 1(C\cup C'=B) \sum_{\by''\in\calV_{A\cup C}} p^{A\cup C,C'}_{\by'',\by'} q^{A,C}_{\by,\by''},
\end{align}
where $q^{A,C}_{\by,\by''}$ is given in Corollary \ref{co:app_gAf_unstable}. 
From the proof of Proposition \ref{pr:app_gA_original}, we see that there exists a positive integer $u^*_{A,1}$ such that if $u_A\ge u^*_{A,1}$, then, for every $\by\in\calV_A\subset\hat{\calV}^A_\emptyset$, 
\begin{equation}
\delta^{1/d}< \sum_{\genfrac{}{}{0pt}{}{C\in\scrP(D\setminus A)}{C\ne\emptyset}} \sum_{\by'\in\hat{\calV}^A_C} q^{A,C}_{\by,\by'} 
=\sum_{\genfrac{}{}{0pt}{}{C\in\scrP(D\setminus A)}{C\ne\emptyset}} \sum_{\by'\in\calV_{A\cup C}} q^{A,C}_{\by,\by'}\le 1, 
\label{eq:app_sum_qABf}
\end{equation}
where we use Proposition \ref{pr:XAn_in_VA}. Furthermore, by the same reason as that used for deriving expression (\ref{eq:sum_pAB}), we obtain, for $\by\in\calV_A\subset\hat{\calV}^A_\emptyset$, 
\begin{align}
&\sum_{B\in\scrD^A_{stable}}\,\sum_{\by'\in\calV_{A\cup B}} p^{A,B}_{\by,\by'}
= \sum_{\genfrac{}{}{0pt}{}{C\in\scrP(D\setminus A)}{C\ne\emptyset}} \sum_{\by''\in\calV_{A\cup C}}\ \sum_{C'\in\scrD^{A\cup C}_{stable}} \sum_{\by'\in\calV_{A\cup C\cup C'}} p^{A\cup C,C'}_{\by'',\by'} q^{A,C}_{\by,\by''}, 
\label{eq:sum_pABf}
\end{align}
and, by the assumption of induction and expression (\ref{eq:app_sum_qABf}), we see that, for every $\by\in\calV_A\subset\hat{\calV}^A_\emptyset$, if $u_A\ge u^*_{A,1}$, then 
\begin{align}
&\delta^{(d-k+2)/d}<
\sum_{B\in\scrD^A_{stable}}\,\sum_{\by'\in\calV_{A\cup B}} p^{A,B}_{\by,\by'}
 \le 1. 
\end{align}
With respect to $g^A_{f,\by}$, we consider the following expression: 
\begin{align}
g^A_{f,\by} - \sum_{B\in\scrD^A_{stable}}\,\sum_{\by'\in\calV_{A\cup B}} p^{A,B}_{\by,\by'}\,g^{A\cup B}_{f,\by'}
= \phi^A_{1,\by}+\phi^A_{2,\by}, 
\end{align}
where 
\begin{align*}
&\phi^A_{1,\by} = g^A_{f,\by} - \sum_{\genfrac{}{}{0pt}{}{B\in\scrP(D\setminus A)}{B\ne\emptyset}} \sum_{\by'\in\calV_{A\cup B}} g^{A\cup B}_{f,\by'}\,q^{A,B}_{\by,\by'}, \\
&\phi^A_{2,\by} = \sum_{\genfrac{}{}{0pt}{}{B\in\scrP(D\setminus A)}{B\ne\emptyset}} \sum_{\by'\in\calV_{A\cup B}} g^{A\cup B}_{f,\by'}\,q^{A,B}_{\by,\by'} - \sum_{B\in\scrD^A_{stable}}\,\sum_{\by'\in\calV_{A\cup B}} p^{A,B}_{\by,\by'}\,g^{A\cup B}_{f,\by'}.
\end{align*}
By Corollary \ref{co:app_gAf_unstable}, there exists a positive integer $u^*_{A,2}$ such that, for every $\by\in\calV_A\subset\hat{\calV}^A_\emptyset$, if $u_A\ge u^*_{A,2}$, then 
\begin{align*}
&|\phi^A_{1,f,\by}| < \varepsilon/d.
\end{align*}
On the other hand, by the same reason as that used for deriving expression (\ref{eq:sum_pAB}), we obtain, for $\by\in\calV_A$, 
\begin{align}
&\sum_{B\in\scrD^A_{stable}}\,\sum_{\by'\in\calV_{A\cup B}} p^{A,B}_{\by,\by'}\,g^{A\cup B}_{f,\by'} 
= \sum_{\genfrac{}{}{0pt}{}{C\in\scrP(D\setminus A)}{C\ne\emptyset}}\,\sum_{\by''\in\calV_{A\cup C}}\,\sum_{C'\in\scrD^{A\cup C}_{stable}}\,\sum_{\by'\in\calV_{A\cup C\cup C'}} p^{A\cup C,C'}_{\by'',\by'}\,g^{A\cup C\cup C'}_{f,\by'}\,q^{A,C}_{\by,\by''}, 
\label{eq:sum_pABaABf}
\end{align}
and, by the assumption of induction and expressions (\ref{eq:app_sum_qABf}) and (\ref{eq:sum_pABaABf}), we obtain 
\begin{align}
|\phi^A_{2,\by}| 
%&= \Bigl| \sum_{\genfrac{}{}{0pt}{}{B\in\scrP(D\setminus A)}{B\ne\emptyset}} \sum_{\by'\in\hat{\calV}^A_B} g^{A\cup B}_{\by'}(l)\,q^{A,B}_{\by,\by'} - \sum_{\genfrac{}{}{0pt}{}{B\in\scrP(D\setminus A)}{B\ne\emptyset}} \sum_{\by'\in\hat{\calV}^A_B}\ \sum_{C\in\scrD^{A\cup B}_{stable}} p^{A\cup B,C}_{\by'} a^{A\cup B\cup C}(l)\,q^{A,B}_{\by,\by'} \Bigr| \nonumber \\
%
&\le \sum_{\genfrac{}{}{0pt}{}{C\in\scrP(D\setminus A)}{C\ne\emptyset}} \sum_{\by''\in\calV_{A\cup C}} \Bigl| g^{A\cup C}_{f,\by''} - \sum_{C'\in\scrD^{A\cup C}_{stable}}\,\sum_{\by'\in\calV_{A\cup C\cup C'}} p^{A\cup C,C'}_{\by'',\by'}\,g^{A\cup C\cup C'}_{f,\by'} \Bigr| q^{A,C}_{\by,\by''} \cr
&< \frac{\varepsilon (d-k+1)}{d} \sum_{\genfrac{}{}{0pt}{}{C\in\scrP(D\setminus A)}{C\ne\emptyset}} \sum_{\by''\in\calV_{A\cup C}} q^{A,C}_{\by,\by''} \cr
&\le \frac{\varepsilon (d-k+1)}{d},
\end{align}
where we use the fact that $|A\cup C|\ge k$ for nonempty $C\in\scrP(D\setminus A)$. 
Thus, for every $\by\in\calV_A\subset\hat{\calV}^A_\emptyset$, if $u_A\ge u^*_{A,2}$, then 
\begin{align}
&\Bigl| g^A_{f,\by} - \sum_{B\in\scrD^A_{stable}}\,\sum_{\by'\in\calV_{A\cup B}} p^{A,B}_{\by,\by'}\,g^{A\cup B}_{f,\by'}\Bigr| 
< \frac{\varepsilon}{d} + \frac{\varepsilon (d-k+1)}{d} 
= \frac{\varepsilon (d-|A|+1)}{d}, 
\end{align}
where $|A|=k-1$. 
We set $u_A$ so that it satisfies $u_A\ge\max\{u^*_{A,1},u^*_{A,2}\}$ and set $K_A$ at a positive integer satisfying $K_A>\max_{|B|>|A|} K_B+u_A$. Then, the assumption of induction holds for the case where $|A|=k-1$, and this completes the proof.
\end{proof}

Finally, we consider a property corresponding to Proposition \ref{pr:talpha_gA}. By the same reason as that used for deriving the proposition, we obtain the following.
\begin{corollary} \label{co:gAf}
For $\by=(\bx,j)\in\calS$, if $\by\in\calV_A$ for some nonempty $A\in\scrP(D)$, then we have $\mathbb{E}(f(\tilde{\bX}_1)-f(\tilde{\bX}_0)\,|\,\tilde{\bY}_0=\by)/u_A=g^A_{f,\by}$.
\end{corollary}

%%%%%%%%%%%%%%%%%%%%%%%%%%%%%
%
\subsubsection{Instability condition for semi-irreducible MMRRWs}

One of the most crucial points in applying Proposition \ref{pr:Markov_unstable2} to semi-irreducible MMRRWs is selection of test function. We here propose to use a function composed of several linear functions. Let $m$ be a positive integer and, for $k\in\{1,2,...,m\}$, let $\bw_k$ be a vector in $\mathbb{R}^d$. For $k\in\{1,2,...,m\}$, we define a linear real function $f_k$ on $\mathbb{R}^d$ as $f_k(\bx)=\langle\bx,\bw_k\rangle$, $\bx\in\mathbb{R}^d$, and a real function $f$ on $\mathbb{R}^d$ as 
\begin{equation}
f(\bx) = \max_{1\le k\le m} f_k(\bx) = \max_{1\le k\le m} \langle\bx,\bw_k\rangle,\ \bx\in\mathbb{R}^d.
\label{eq:testfunction_f}
\end{equation}
We also use (abuse) notation $f$ for denoting the test function to analyze stability of semi-irreducible MMRRWs, i.e., for $\by=(\bx,j)\in\calS$, we have 
\[
f(\by)=f(\bx)=\max_{1\le k\le m} f_k(\bx).
\]
Note that the test function used for deriving the stability condition in Theorem \ref{th:stable} is this $f$ when $k=1$ and $\bw_1>0$. Define a constant $C$ and a function $k^*$ as 
\[
C = \max_{1\le k\le m} \| \bw_k \|_1,\quad 
k^*(\bx) = \arg\max_{1\le k\le m} f_k(\bx),\ \bx\in\mathbb{R}^d,  
\]
then we have, for $\bx,\bx'\in\mathbb{R}^d$,
\[
f(\bx)=f_{k^*(\bx)}(\bx)\ge f_{k^*(\bx')}(\bx),\quad 
|f(\bx)| \le \max_{1\le k\le m}  \| \bw_k \|_1\,\| \bx \|_1 \le C\,\| \bx \|_1. 
\]
This and linearity of $f_k$ lead us to that, for $\bx,\bx'\in\mathbb{R}^d$, 
\[
|f(\bx')-f(\bx)| \le \max\{f_{k^*(\bx')}(\bx')-f_{k^*(\bx')}(\bx),\,f_{k^*(\bx)}(\bx)-f_{k^*(\bx)}(\bx') \}
\le C\,\| \bx'-\bx \|_1, 
\]
and hence the function $f$ satisfies condition (\ref{eq:f_condition}).
Before stating an instability condition for MMRRWs, we define the following notation. For a real function $f$ given by expression (\ref{eq:testfunction_f}), let $\calA_f$ be a subset of $\calS$ defined as
\[
\calA_f = \{\by=(\bx,j)\in\calS: f(\by)=f(\bx)>c\}, 
\]
where $c$ is a sufficiently large constant. 
Consider the partition of $\calS$, $\{\calV_A: A\in\scrP(D)\}$. Let $\scrD_f$ be the index set of the elements of the partition that intersect $\calA_f$, i.e., 
\[
\scrD_f = \{A\in\scrP(D): \calV_A\cap\calA_f\ne\emptyset\}, 
\]
and $\scrD_{f,stable}$ the index set of stable induced Markov chains, defined by 
\[
\scrD_{f,stable} = \{A\cup B: A\in\scrD_f,\,B\in\scrD^A_{stable}\}.
\]
Furthermore, for $A\in\scrP(D)$, denote by $\scrI_{f,A}$ the index set of the linear functions that attain the maximum on some state in $\calV_A$, i.e.,
\[
\scrI_{f,A} = \{k\in\{1,2,...,m\}: f(\by)=f_k(\bx)\ \mbox{for some $\by=(\bx,j)\in\calV_A$} \}.  
\]
Note that, for $\by=(\bx,j)\in\calV_A$ and $\bx'\in\mathbb{R}^d$,  
\begin{equation}
f_{k^*(\bx)}(\bx') \ge \min_{k\in\scrI_{f,A}} f_k(\bx').
\label{eq:lb_fkstar}
\end{equation}
The following theorem gives an instability condition for semi-irreducible MMRRWs and it corresponds to Condition B$'$ and Theorem 2.1 of Malyshev and Menshikov \cite{Malyshev81} (also see Condition B$'$ and Theorem 4.3.4 of Fayolle et al.\ \cite{Fayolle95}). 
%
%%%%%%%%%%%%%%%%%%%%%%%%%%
\begin{theorem}[Unstable MMRRW] \label{th:unstable}
The skip-free semi-irreducible MMRRW $\calL=\{\bY_n\}=\{(\bX_n,J_n)\}$ is unstable in our sense if there exists a real function $f$ given by expression (\ref{eq:testfunction_f}) such that, for some $c>0$, $\calA_f\ne\emptyset$ and $\calA_f^C\cap\calS_0\ne\emptyset$, and furthermore, for every $A\in\scrD_{f,stable}$ and for every $k\in\scrI_{f,A}$, we have $f_k(\ba^A)>0$, where $\calS_0$ is the irreducible class of the MMRRW and $\ba^A$ the mean drift vector of the expanded Markov chain $\hat{\calL}^A$. 
\end{theorem}

\begin{proof}
Let $\varepsilon_0$ arbitrarily set at a small positive number and $\delta_0$ at a number in $(0,1)$ sufficiently close to $1$. 
First, we determine the parameter set $\{u'_A: A\in\scrP(D),\ A\ne\emptyset\}$, which is appeared in Corollary \ref{co:app_gAf_original}. For $A\in\scrD_{stable}$, set $u'_A$ so that, for every $\by\in\calV_A$ and every $l\in D$, inequality 
\[
\Bigl| \frac{1}{u'_A} \mathbb{E}\big(X^A_{u'_A}(l)-X^A_0(l)\,|\,\bY^A_0=\by\big) - a^A(l)\Bigr| = \bigl|g^A_{\by}(l)- a^A(l)\bigr|< \varepsilon_0/d
\]
holds, where we give $\bg^A_{\by}$ by using $u'_A$ instead of $u_A$; it is possible by Propositions \ref{pr:app_talpha} and \ref{pr:talpha_gA} since the induced Markov chain $\calL^A$ is stable. If $A\notin\scrD_{stable}$, then we set $u'_A=1$. 
Using this parameter set $\{u'_A\}$, we obtain, by Corollary \ref{co:app_gAf_original}, the set of positive integers $\{(K_A,u_A) : A\in\scrP(D),\,A\ne\emptyset\}$ and the set of nonnegative vectors $\{(p^{A,B}_{\by,\by'},\,B\in\scrD^A_{stable},\by'\in\calV_{A\cup B}): A\in\scrP(D),\,A\ne\emptyset,\,\by\in\calV_A\}$, for which expressions (\ref{eq:app_pABf_original}) and (\ref{eq:app_gAf_original}) hold. 

Next, we apply Proposition \ref{pr:Markov_unstable2} to the MMRRW. 
By the assumption of the theorem, we have $\calA_f\ne\emptyset$ and $\calA_f^C\cap\calS_0\ne\emptyset$, and condition (i) of Proposition \ref{pr:Markov_unstable2} holds. 
Let $\varepsilon^*$ be defined as 
\[
\varepsilon^* = \min_{A\in\scrD_{f,stable}} \min_{k\in\scrI_{f,A}} f_k(\ba^A)>0. 
\]
For $A\in\scrD_f$ and for $\by=(\bx,j)\in\calV_A$, we have, by Corollary \ref{co:gAf}, 
\begin{align}
&\mathbb{E}\big(f(\tilde{\bY}_{n+1})-f(\tilde{\bY}_n)\,|\,\tilde{\bY}_n=\by\big)
= \mathbb{E}\big(f(\tilde{\bX}_{n+1})-f(\tilde{\bX}_n)\,|\,\tilde{\bY}_n=\by\big)
= u_A\,g^A_{f,\by}. 
\end{align}
If the induced Markov chain $\calL^A$ is stable, then $A\in\scrD_{f,stable}$ and we have
\begin{align}
g^A_{f,\by}
&\ge \frac{1}{u_A} \mathbb{E}\big( f_{k^*(\bx)}(\bX^A_{u_A})-f_{k^*(\bx)}(\bX^A_0)\,|\,\bY_0=\by\big) \cr
&= f_{k^*(\bx)}\big( (1/u_A) \mathbb{E}(\bX^A_{u_A}-\bX^A_0\,|\,\bY_0=\by) \big) \cr
&= f_{k^*(\bx)}(\bg^A_{\by}) \cr
&\ge \min_{k\in\scrI_{f,A}} f_k(\ba^A) + f_{k^*(\bx)}(\bg^A_{\by}-\ba^A) \cr
&\ge \varepsilon^* - C\,\|\bg^A_{\by}-\ba^A\|_1 \cr
&\ge \varepsilon^* - C\,\varepsilon_0 > 0, 
\label{eq:lb_gAf_stable}
\end{align}
where we use expression (\ref{eq:lb_fkstar}) and the fact that $\varepsilon_0$ is arbitrarily small. 
If the induced Markov chain $\calL^A$ is unstable, then, by expression (\ref{eq:lb_gAf_stable}) and Corollary \ref{co:app_gAf_original}, we obtain 
\begin{align}
g^A_{f,\by} 
&= \sum_{B\in\scrD^A_{stable}}\,\sum_{\by'\in\calV_{A\cup B}} p^{A,B}_{\by,\by'}\,g^{A\cup B}_{f,\by'} + \biggl(g^A_{f,\by} - \sum_{B\in\scrD^A_{stable}}\,\sum_{\by'\in\calV_{A\cup B}} p^{A,B}_{\by,\by'}\,g^{A\cup B}_{f,\by'}\biggr) \cr
&\ge \sum_{B\in\scrD^A_{stable}}\,\sum_{\by'\in\calV_{A\cup B}} p^{A,B}_{\by,\by'} (\varepsilon^*-C\,\varepsilon_0) - \varepsilon_0 \cr
&> \delta_0 (\varepsilon^*-C\,\varepsilon_0) - \varepsilon_0 >0,
\end{align}
where we use the fact that $A\cup B\in\scrD_{f,stable}$ and that $\delta_0$ is sufficiently close to $1$ and $\varepsilon_0$ sufficiently small.
As a result, setting $\varepsilon=\delta_0 (\varepsilon^*-C\,\varepsilon_0) - \varepsilon_0>0$, we obtain, for $A\in\scrD_f$ and for $\by=(\bx,j)\in\calV_A$, 
\begin{equation}
\mathbb{E}\big(f(\tilde{\bY}_{n+1})-f(\tilde{\bY}_n)\,|\,\tilde{\bY}_n=\by\big)
> u_A\,\varepsilon > \varepsilon>0. 
\label{eq:lb_Ef}
\end{equation}
Since $\calA_f=\bigcup_{A\in\scrD_f}(\calA_f\cap\calV_A)$, inequality (\ref{eq:lb_Ef}) holds for all $\by\in\calA_f$, and condition (ii) of Proposition \ref{pr:Markov_unstable2} holds. 
Since the process $\{\bX_n\}$ is skip free in all coordinates, we have 
\[
|f(\bY_1)-f(\bY_0)| = |f(\bX_1)-f(\bX_0)| \le C\,\| \bX_1-\bX_0 \|_1 \le C d, 
\]
and this leads us to that, for $\by,\by'\in\calS$, if $|f(\by')-f(\by)|>C d$, then $\mathbb{P}(\bY_1=\by'\,|\,\bY_0=\by)=0$. Thus, condition (iii) of Proposition \ref{pr:Markov_unstable2} holds and this completes the proof. 
\end{proof}

%%%%%%%%%%%%%%%%%%%%%%%%%%%%%%
%
In applying Theorem \ref{th:unstable} to a semi-irreducible MMRRW, it seems difficult to select a proper test function $f$ satisfying the condition of the theorem. Therefore, we recommend using the following type of function as a test function. 
Let $\bw_0$ be a $d$-dimensional positive vector and $c_0$ a sufficiently large positive number. Let $m$ be a positive integer and $A_k,\,k=1,2,...,m,$ nonempty different elements of $\scrP(D)$ satisfying the following conditions: 
\begin{align}
&\mbox{$|A_k|=|A_{k'}|$ for all $k,k'\in\{1,2,...,m\}$}, \label{eq:Ak_condition1} \\
&\mbox{$|A_k\setminus A_{k'}|=|A_{k'}\setminus A_k|=1$ for all $k,k'\in\{1,2,...,m\}$ such that $k\ne k'$}. \label{eq:Ak_condition2}
\end{align}
For $k\in\{1,2,...,m\}$, let $\bw_k$ be defined as
\begin{equation}
\bw_k = \bw_0 - c_0 \sum_{l\in D\setminus A_k} \be_l, 
\label{eq:wk}
\end{equation}
where $\be_l$ is the $l$-th unit vector, and $f_k(\bx)$ as $f_k(\bx)=\langle \bx,\bw_k \rangle$; the test function $f$ is given by expression (\ref{eq:testfunction_f}). For $l\in A_k$, hyperplane $f_k(\bx)=c$ intersects $x_l$-axis at $x(l)=c/w_0(l)>0$ and, for $l\in D\setminus A_k$, it intersects $x_l$-axis at $x(l)=c/(w_0(l)-c_0)<0$, where $c_0$ is a sufficiently large positive number.
For $\bx\in\mathbb{R}_+^d$, if $f_k(\bx)>c$, then $\sum_{l\in A_k} x(l) w_0(l) >c$ and we have
\begin{equation}
\max_{l\in A_k} x(l) > \frac{c}{\sum_{l\in A_k} w_0(l)}. 
\label{eq:lb_xl}
\end{equation}
Note that, for nonempty $A\in\scrP(D)$, if $(\bx,j)\in\calV_A$, then $\bx(A)\ge K_A$ and $\bx(D\setminus A)<K_A$. From this, we see that, for $l\in A_k$ and for nonempty $A\in\scrP(D)$ such that $l\notin A$, if $(\bx,j)\in\calV_A$, then $x(l)<K_A$. Thus, for $(\bx,j)\in\calS$, if $x(l)\ge\max_{A\in\scrP(D)} K_A$, then $(\bx,j)\in\calV_{B\cup\{l\}}$ for some $B\in\scrP(D\setminus\{l\})$. This and inequality (\ref{eq:lb_xl}) lead us to that, for a sufficiently large $c$, 
\[
\{\by=(\bx,j)\in\calS: f_k(\bx)>c\} \subset \bigcup_{l\in A_k}\,\bigcup_{B\in\scrP(D\setminus\{l\})} \calV_{B\cup\{l\}}, 
\]
and we obtain 
\begin{equation}
\calA_f \subset \bigcup_{l\in\,\bigcup_{k=1}^m A_k}\ \bigcup_{B\in\scrP(D\setminus\{l\})} \calV_{B\cup\{l\}}. 
\label{eq:AfV}
\end{equation}
Thus, for a sufficiently large $c$, defining $\bar{\scrD}_f$ as 
\begin{equation}
\bar{\scrD}_f = \Bigl\{\{l\}\cup B: l\in\bigcup_{k=1}^m A_k,\ B\in\scrP(D\setminus\{l\})\Bigr\}, 
\label{eq:barDf}
\end{equation}
we have $\scrD_f\subset\bar{\scrD}_f$. Further, define $\bar{\scrD}_{f,stable}$ as 
\begin{equation}
\bar{\scrD}_{f,stable} 
= \{A\cup B: A\in\bar{\scrD}_f,\ B\in\scrD^A_{stable}\}
= \Bigl\{\{l\}\cup B: l\in\bigcup_{k=1}^m A_k,\ B\in\scrD^{\{l\}}_{stable}\Bigr\}, 
\label{eq:barDfstable}
\end{equation}
then we have $\scrD_{f,stable}\subset\bar{\scrD}_{f,stable}$.
Next, we consider $\scrI_{f,A}$ for nonempty $A\in\scrP(D)$. For $k,k'\in\{1,2,...,m\}$ such that $k\ne k'$, denote by $l_{k,k'}$ the unique element of $A_k\setminus A_{k'}$. Then, we have, for $k,k'\in\{1,2,...,m\}$ such that $k\ne k'$ and for $\bx\in\mathbb{R}_+^d$, 
\begin{equation}
f_k(\bx)-f_{k'}(\bx) 
= -c_0 \Big\langle \bx, \sum_{l\in D\setminus A_k} \be_l-\sum_{l\in D\setminus A_{k'}} \be_l \Big\rangle 
= c_0 (x(l_{k,k'})-x(l_{k',k})). 
\end{equation}
Thus, for $k\in\{1,2,...,m\}$ and for nonempty $A\in\scrP(D)$, if $l_{k,k'}\in D\setminus A$ and $l_{k',k}\in A$ for some $k'\in\{1,2,...,m\}$, then $f_k(\bx)-f_{k'}(\bx)<0$ for all $(\bx,j)\in\calV_A$ and $f_k$ does not assign $f$ in $\calV_A$. Therefore, defining $\bar{\scrI}_{f,A}$ as 
\begin{align}
\bar{\scrI}_{f,A} &= \{k\in\{1,2,...,m\}: \mbox{$A_k\setminus A_{k'}\subset A$ or ($A_k\setminus A_{k'}\subset D\setminus A$ and $A_{k'}\setminus A_k\subset D\setminus A$)} \cr
&\qquad\qquad\qquad\qquad \mbox{for all $k'\in\{1,2,...,m\}$ such that $k'\ne k$}\}, 
\label{eq:barIfA}
\end{align}
we have $\scrI_{f,A}\subset\bar{\scrI}_{f,A}$. 
Furthermore, by expression (\ref{eq:AfV}), we see that, for a sufficiently large $c$, $\calV_\emptyset\subset\calA_f^C$. 
As a result, we obtain from Theorem \ref{th:unstable} the following corollary.
%
%%%%%%%%%%%%%%%%%%
\begin{corollary}[Unstable MMRRW] \label{co:unstable2}
Assume $\calV_\emptyset\cap\calS_0\ne\emptyset$. Then, the skip-free semi-irreducible MMRRW $\calL=\{\bY_n\}=\{(\bX_n,J_n)\}$ is unstable in our sense if there exists a real function $f$ given by expressions (\ref{eq:testfunction_f}) and (\ref{eq:wk}) such that, for every $A\in\bar{\scrD}_{f,stable}$ and for every $k\in\bar{\scrI}_{f,A}$, we have $f_k(\ba^A)>0$. 
\end{corollary}

\begin{example}[Two-dimensional QBD process; see Examples \ref{ex:2D_QBD1}, \ref{ex:2D_QBD2} and \ref{ex:2D_QBD3}] \label{ex:2D_QBD4} {\rm 
For the two-dimensional QBD process in Examples \ref{ex:2D_QBD1}, consider the case where $\ba^D<0$, $a^{\{1\}}(1)>0$ and $a^{\{2\}}(2)<0$. In this case, we have $\scrD_{stable}=\{\{1\},\{2\},D\}$, where $D=\{1,2\}$, and $a^{\{1\}}(2)$ and $a^{\{2\}}(1)$ are zero. Furthermore, we assume $\calV_0\cap\calS_0\ne\emptyset$. 

Set $m=1$ and $A_1=\{1\}$; then, the test function $f$ is given by
\[
f(\by)=f_1(\bx) = \langle\bx,\bw_1\rangle = w_0(1)\,x(1)+(w_0(2)-c_0) x(2),\ \by=(\bx,j)\in\calS, 
\]
where $\bw_1=( w_0(1),w_0(2))-c_0(0,1)=( w_0(1),w_0(2)-c_0)$. For a sufficiently large number $c$, which is used for defining $\calA_f$, we have 
\[
\bar{\scrD}_{f,stable}=\{\{1\},D\},\quad
\bar{\scrI}_{f,\{1\}}=\{1\},\quad \bar{\scrI}_{f,D}=\{1\}, 
\]
where we set $c$ so large that $\calV_0\cap\calA_f=\emptyset$. 
Set $w_0(1)=w_0(2)=1$ and $c_0=2+a^D(1)/a^D(2)>0$, then we have
\begin{align*}
&f_1(\ba^{\{1\}})=\langle\ba^{\{1\}},\bw_1\rangle = a^{\{1\}}(1)>0,\cr
&f_1(\ba^D)=\langle\ba^D,\bw_1\rangle = a^D(1)+(w_0(2)-c_0) a^D(2) = -a^D(2)>0. 
\end{align*}
Hence, by Corollary \ref{co:unstable2}, the two-dimensional QBD process is unstable in our sense. 
\hfill$\Box$
} \end{example}

%%%%%%%%%%%%%%%%%%%%%%%
%
% Section 3, continuation
%
%%%%%%%%%%%%%%%%%%%%%%%
%
\subsection{Low dimensional cases} \label{sec:low_dimensional_cases}

Applying Theorem \ref{th:stable} and Corollary \ref{co:unstable2}, we obtain classification of low-dimensional MMRRWs with respect to stability, where the mean drift vectors are used for classifying them. The dimensions we consider are one through three, i,e, $d=1,2,3$.

%%%%%%%%%%%%%%%%%%%%%%%%%%%%
%
\subsubsection{One dimensional case}

Here we deal with a one-dimensional semi-irreducible MMRRW (1D-MMRRW for short) $\calL=\{\bY_n\}$ on the state space $\calS$ and assume $(\{0\}\times S^\emptyset)\cap\calS_0\ne\emptyset$ for simplicity, where $\calS_0$ is the unique irreducible class of $\calL$. 
As mentioned in Example \ref{ex:QBD1}, the 1D-MMRRW is a QBD process whose transition probability matrix $P$ is given by expression (\ref{eq:QBD_P}). Our model is semi-irreducible, but the classification is similar to that of irreducible QBD processes (see, for example, Neuts \cite{Neuts94} and Latouche and Rawaswami \cite{Latouche99}). 

The 1D-MMRRW has only one induced Markov chain $\calL^D$, where $D=\{1\}$, and by Assumption \ref{as:inducedMC}, it is semi-irreducible and stable in our sense. The transition probability matrix of $\calL^D$ is given by (for the notation, see Example \ref{ex:QBD1})
\[
P^D = P_{-1}^{\{1\},\{1\}}+P_0^{\{1\},\{1\}}+P_1^{\{1\},\{1\}},
\]
and the mean drift $a^D$ is given as 
\[
a^D = \bpi^D \bigl(-P_{-1}^{\{1\},\{1\}} + \bpi^D P_1^{\{1\},\{1\}}\bigr) \bone, 
\]
where $\bpi^D$ is the stationary probability vector of $P^D$ and $\bone$ is a column vector of $1$'s whose dimension is determined in context. 
The classification of 1D-MMRRWs with respect to stability is given as follows.
%
%%%%%%%%%%%%%%%%%%%%%%%%%%%
\begin{theorem} \label{th:classification_1DMMRRW}
The 1D-MMRRW is stable in our sense if $a^D<0$, and it is unstable in our sense if $a^D>0$. 
\end{theorem}
This theorem can be proved by Theorem \ref{th:stable} and Corollary \ref{co:unstable2}, but the proof is straightforward and we omit it.

%%%%%%%%%%%%%%%%%%%%%%%%%%%%%%%
%
\subsubsection{Two dimensional case}

Here we deal with a two-dimensional semi-irreducible MMRRW (2D-MMRRW for short) $\calL=\{\bY_n\}$ on the state space $\calS$ and assume $(\{0\}\times\{0\}\times S^\emptyset)\cap\calS_0\ne\emptyset$ for simplicity, where $\calS_0$ is the unique irreducible class of $\calL$. 
As mentioned in Example \ref{ex:2D_QBD1}, the 2D-MMRRW is a two-dimensional QBD process, and as mentioned in Example \ref{ex:2D_QBD2}, it has three induced Markov chains $\calL^{\{1\}}$, $\calL^{\{2\}}$ and $\calL^D$, where $D=\{1,2\}$, whose transition probability matrices $P^{\{1\}}$, $P^{\{2\}}$ and $P^D$ are given by expressions (\ref{eq:2DQBD_P1}), (\ref{eq:2DQBD_P2}) and (\ref{eq:2DQBD_PD}), respectively. 
By Assumption \ref{as:inducedMC}, $\calL^D$ is semi-irreducible and stable in our sense, and the mean drift vector $\ba^D$ is given as (for the notation, see Examples \ref{ex:2D_QBD1})
\[
a^D(1) = \bpi^D \bigl(-P_{(-1,*)}^{D,D}+P_{(1,*)}^{D,D}\bigr) \bone,\quad 
a^D(2) = \bpi^D \bigl(-P_{(*,-1)}^{D,D}+P_{(*,1)}^{D,D}\bigr) \bone, 
\]
where $\bpi^D$ is the stationary probability vector of $P^D$. 
The induced Markov chains $\calL^{\{1\}}$ and $\calL^{\{2\}}$ are one-dimensional MMRRWs. Hence, by Theorem \ref{th:classification_1DMMRRW}, $\calL^{\{1\}}$ is stable if $a^D(2)<0$ and it is unstable if $a^D(2)>0$; $\calL^{\{2\}}$ is stable if $a^D(1)<0$ and it is unstable if $a^D(1)>0$.  
For $l\in\{1,2\}$, by Theorem 6.2.1 of Latouche and Rawaswami \cite{Latouche99}, which also holds in our models, when $\calL^{\{l\}}$ is stable, its stationary distribution $\bpi^{\{l\}}$ is given as, in matrix geometric form, 
\[
\bpi^{\{l\}}=(\bpi_k^{\{l\}},k\in\mathbb{Z}_+),\quad \bpi_k^{\{l\}}=\bpi_k^{\{l\}} (R^{\{l\}})^{k-1},\ k\ge 2, 
\]
where $R^{\{l\}}$ is the rate matrix of $P^{\{l\}}$. Hence, the mean increment vectors $\ba^{\{1\}}$ and $\ba^{\{2\}}$ are given as 
\begin{align*}
a^{\{1\}}(1) &= \bpi_0^{\{1\}}\big(-P_{(-1,0)}^{\{1\},\{1\}}-P_{(-1,1)}^{\{1\},D}+P_{(1,0)}^{\{1\},\{1\}}+P_{(1,1)}^{\{1\},D} \big)\,\bone \cr
&\qquad + \bpi_1^{\{1\}}\big(-P_{(-1,-1)}^{D,\{1\}}-P_{(-1,0)}^{D,D}-P_{(-1,1)}^{D,D}+P_{(1,-1)}^{D,\{1\}}+P_{(1,0)}^{D,D}+P_{(1,1)}^{D,D} \big)\,\bone \cr
&\qquad + \bpi_2^{\{1\}}\big(I-R^{\{1\}}\big)^{-1} \big(-P_{(-1,*)}^{D,D}+P_{(1,*)}^{D,D} \big)\,\bone,\\
a^{\{1\}}(2) &= 0,\quad 
a^{\{2\}}(1) = 0, \\
a^{\{2\}}(2) &= \bpi_0^{\{2\}}\big(-P_{(0,-1)}^{\{2\},\{2\}}-P_{(1,-1)}^{\{2\},D}+P_{(0,1)}^{\{2\},\{2\}}+P_{(1,1)}^{\{2\},D} \big)\,\bone \cr
&\qquad + \bpi_1^{\{2\}}\big(-P_{(-1,-1)}^{D,\{2\}}-P_{(0,-1)}^{D,D}-P_{(1,-1)}^{D,D}+P_{(-1,1)}^{D,\{2\}}+P_{(0,1)}^{D,D}+P_{(1,1)}^{D,D} \big)\,\bone \cr
&\qquad + \bpi_2^{\{2\}}\big(I-R^{\{2\}}\big)^{-1} \big(-P_{(*,-1)}^{D,D}+P_{(*,1)}^{D,D} \big)\,\bone.
\end{align*}
The classification of 2D-MMRRWs with respect to stability is given by the following theorem; a similar result was obtained in Fayolle \cite{Fayolle89} for ordinary two-dimensional reflecting random walks (2D-RRWs for short). 
%
%%%%%%%%%%%%%%%%%%
\begin{theorem} \label{th:classification_2DMMRRW}
Assume $(\{0\}\times\{0\}\times S^\emptyset)\cap\calS_0\ne\emptyset$. 
\begin{itemize}
\item[(C1)] When $a^D(1)<0$ and $a^D(2)<0$, the 2D-MMRRW $\calL$ is
\begin{itemize}
\item[(a)] stable in our sense if $a^{\{1\}}(1)<0$ and $a^{\{2\}}(2)<0$; 
\item[(b)] unstable in our sense if either  $a^{\{1\}}(1)>0$ or $a^{\{2\}}(2)>0$.
\end{itemize}
\item[(C2)] When $a^D(1)>0$ and $a^D(2)<0$, the 2D-MMRRW $\calL$ is
\begin{itemize}
\item[(a)] stable in our sense if $a^{\{1\}}(1)<0$; 
\item[(b)] unstable in our sense if $a^{\{1\}}(1)>0$.
\end{itemize}
\item[(C3)] When $a^D(1)<0$ and $a^D(2)>0$, the 2D-MMRRW $\calL$ is
\begin{itemize}
\item[(a)] stable in our sense if $a^{\{2\}}(2)<0$; 
\item[(b)] unstable in our sense if $a^{\{2\}}(2)>0$.
\end{itemize}
\item[(C4)] When $a^D(1)>0$ and $a^D(2)>0$, the 2D-MMRRW $\calL$ is unstable in our sense.
\end{itemize}
\end{theorem}

%%%%%%%%%
\begin{proof}
Note that $(\{0\}\times\{0\}\times S^\emptyset)\subset\calV_\emptyset$ and hence, by the assumption, we have $\calV_\emptyset\cap\calS_0\ne\emptyset$. 
(C1)-(a) was already proved in Example \ref{ex:2D_QBD3}, and (C1)-(b) when $a^{\{2\}}(1)>0$ was already proved in Example \ref{ex:2D_QBD4}; (C1)-(b) when $a^{\{2\}}(2)>0$ can be proved in a manner similar to that used in Example \ref{ex:2D_QBD4}.
In (C2), we have $\scrD_{stable}=\{ \{1\},D\}$. 
First, consider (C2)-(a). Set $\bw=(1,1-a^D(1)/a^D(2))>0$, then we have
\[
\langle\ba^{\{1\}},\bw\rangle = a^{\{1\}}(1) < 0,\quad 
\langle\ba^D,\bw\rangle = a^D(2) < 0. 
\]
Hence, by Theorem \ref{th:stable}, the 2D-MMRRW is stable in our sense. 
Next, consider (C2)-(b). Set $m=1$ and $A_1=\{1\}$; then, the test function $f$ is given as
\[
f(\by)=f_1(\bx) = \langle\bx,\bw_1\rangle = w_0(1) x(1)+(w_0(2)-c_0) x(2),\ \by=(\bx,j)\in\calS, 
\]
where $\bw_1=(w_0(1),w_0(2)-c_0)$. For a sufficiently large number $c$, which is used for defining $\calA_f$, we have $\bar{\scrD}_{f,stable}=\{\{1\},D\}$, $\bar{\scrI}_{f,\{1\}}=\{1\}$ and $\bar{\scrI}_{f,D}=\{1\}$, where we set $c$ so large that $\calV_0\cap\calA_f=\emptyset$. 
Set $w_0(1)=w_0(2)=1$ and $c_0=2$, then we have
\begin{align*}
&f_1(\ba^{\{1\}})=a^{\{1\}}(1)>0,\quad 
f_1(\ba^D)=a^D(1)-a^D(2)>0. 
\end{align*}
Thus, by Corollary \ref{co:unstable2}, the 2D-MMRRW is unstable in our sense. 
(C3) is symmetric to (C2), thus the proof is analogous. 
In (C4),  we have $\scrD_{stable}=\{D\}$. Set $m=1$ and $A_1=\emptyset$; then, the test function $f$ is given by
\[
f(\by)=f_1(\bx) = \langle\bx,\bw_1\rangle = w_0(1) x(1)+w_0(2) x(2),\ \by=(\bx,j)\in\calS, 
\]
where $\bw_1=(w_0(1),w_0(2))$. For a sufficiently large number $c$, we have $\bar{\scrD}_{f,stable}=\{D\}$ and $\bar{\scrI}_{f,D}=\{1\}$, where we set $c$ so large that $\calV_0\cap\calA_f=\emptyset$. 
Set $w_0(1)=w_0(2)=1$, then we have
\begin{align*}
&f_1(\ba^D)=\langle\ba^D,\bw_1\rangle = a^D(1)+a^D(2)>0. 
\end{align*}
Thus, by Corollary \ref{co:unstable2}, the 2D-MMRRW is unstable in our sense. 
\end{proof}

\begin{remark}
The case where $a^D(1)$ or $a^D(2)$ is zero cannot be classified by our results. The case where at least one of $a^{\{1\}}(1)$ and $a^{\{2\}}(2)$ is zero cannot also be classified by our results.
\end{remark}

%%%%%%%%%%%%%%%%%%%%%%%%%%%%%%%%%%%%%%%%%
%
\subsubsection{Three dimensional case}

Here we deal with a three-dimensional semi-irreducible MMRRW (3D-MMRRW for short) $\calL=\{\bY_n\}$ on the state space $\calS$ and assume $(\{0\}\times\{0\}\times\{0\}\times S^\emptyset)\cap\calS_0\ne\emptyset$ for simplicity, where $\calS_0$ is the unique irreducible class of $\calL$. 
In this case, $D=\{1,2,3\}$ and the 3D-MMRRW has seven induced Markov chains: $\calL^D$, $\calL^{\{1,2\}}$, $\calL^{\{2,3\}}$, $\calL^{\{3,1\}}$, $\calL^{\{1\}}$, $\calL^{\{2\}}$ and $\calL^{\{3\}}$. 
$\calL^D$ is a semi-irreducible finite Markov chain and its stationary distribution can be calculated by a standard method. $\calL^{\{1,2\}}$, $\calL^{\{2,3\}}$ and $\calL^{\{3,1\}}$ are semi-irreducible QBD processes, and if they are stable in our sense, their stationary distributions are given in the matrix geometric form. Therefore, we can obtain the mean drift vectors $\ba^D$, $\ba^{\{1,2\}}$, $\ba^{\{2,3\}}$ and $\ba^{\{3,1\}}$ if they exist. 
Furthermore, $\calL^{\{1\}}$, $\calL^{\{2\}}$ and $\calL^{\{3\}}$ are 2D-MMRRWs, and we can classify them by Theorem \ref{th:classification_2DMMRRW} with respect to stability. 
However, as far as we know, there are no general methods to obtain the stationary distributions of 2D-MMRRWs and it is generally difficult to get the mean drift vectors $\ba^{\{1\}}$, $\ba^{\{2\}}$ and $\ba^{\{3\}}$ even if we know that the corresponding induced Markov chains are stable. Hence, hereafter, supposing the signs of $\ba^{\{1\}}$, $\ba^{\{2\}}$ and $\ba^{\{3\}}$ are known, we classify 3D-MMRRWs with respect to stability. 

The classification can be done in a manner similar to that used for classifying 2D-MMRRWs; but in this case, there are many combinations of the signs of the mean drift vectors. Hence, we summarize only typical cases on Tables \ref{tab:table1} through \ref{tab:table4}, where, for example, ``$\ba^{\{1,2\}}=(+\,-\,0)$" means that $a^{\{1,2\}}(1)>0$, $a^{\{1,2\}}(2)<0$ and $a^{\{1,2\}}(3)=0$; ``NA'' means that the corresponding induced Markov chain is unstable in our sense, blanks mean anything, and ``stable" (resp.\ ``unstable") means that the 3D-MMRRW $\calL$ is stable (resp.\ unstable) in our sense. 
Of course, the tables do not include all the cases, but every case is shown on some table or symmetric to someone on some table; for example, the case where $\ba^D=(-\,-\,-)$, $\ba^{\{1,2\}}=(-\,-\,0)$, $\ba^{\{2,3\}}=(0\,+\,-)$, $\ba^{\{3,1\}}=(-\,0\,-)$, $\ba^{\{1\}}=(-\,0\,0)$, $\ba^{\{2\}}=(0\,-\,0)$ and $\ba^{\{3\}}=\mbox{NA}$ is not shown on any table but it is symmetric to C1-2-1 on Table \ref{tab:table1}. 
We prove only C1-1-1 and C1-1-2; other cases except for C1-7-1 are analogously proved. C1-7-1 will be explained in Theorem \ref{th:3D_spiral}.

%%%%%%%%%%%%%%%%
\begin{proof}[Proof of C1-1-1 and C1-1-2]
Note that $\{0\}\times\{0\}\times\{0\}\times S^\emptyset\subset\calV_\emptyset$; hence, by the assumption, we have $\calV_\emptyset\cap\calS_0\ne\emptyset$. 
First, we consider C1-1-1. In this case, we have 
\[
\scrD_{stable}=\{ \{1\},\{2\},\{3\},\{1,2\},\{2,3\},\{3,1\},D\}.
\]
Set $\bw=(1,1,1)$, then we have, for all $A\in\scrD_{stable}$, $f(\ba^A)=\langle\ba^A,\bw\rangle<0$. Therefore, by Theorem \ref{th:stable}, the 3D-MMRRW is stable in our sense. 
Next, we consider C1-1-2. In this case, we also have 
\[
\scrD_{stable}=\{ \{1\},\{2\},\{3\},\{1,2\},\{2,3\},\{3,1\},D\}.
\]
Set $m=1$ and $A_1=\{1\}$; then, we obtain, from expressions (\ref{eq:barDfstable}) and (\ref{eq:barIfA}), 
\[
\bar{\scrD}_{f,stable}=\{\{1\},\{1,2\},\{1,3\},D\},\quad
\bar{\scrI}_{f,A}=\{1\}\ \mbox{for $A\in\bar{\scrD}_{f,stable}$}. 
\]
Set $\bw_0=(1,1,1)$ and $c_0$ as
\[
c_0=\frac{\|\ba^{\{1,2\}}\|_1+\|\ba^{\{1,3\}}\|_1+\|\ba^D\|_1}{\min\{|a^{\{1,2\}}(2)|,|a^{\{1,3\}}(3)|,|a^D(2)|,|a^D(3)|\}}.
\]
Since $\bw_1=\bw_0-c_0(\be_2+\be_3)=(1,1-c_0,1-c_0)$ and $f_1(\bx)=\langle\bx,\bw_1\rangle$, we have
\begin{align*}
&f_1(\ba^{\{1\}}) = a^{\{1\}}(1)>0,\cr
&f_1(\ba^{\{1,2\}}) \ge a^{\{1,2\}}(1)+a^{\{1,2\}}(2)+\|\ba^{\{1,2\}}\|_1+\|\ba^{\{1,3\}}\|_1+\|\ba^D\|_1>0,\cr
&f_1(\ba^{\{1,3\}}) \ge a^{\{1,3\}}(1)+a^{\{1,3\}}(3)+\|\ba^{\{1,2\}}\|_1+\|\ba^{\{1,3\}}\|_1+\|\ba^D\|_1>0,\cr
&f_1(\ba^D) \ge a^D(1)+a^D(2)+a^D(3)+2(\|\ba^{\{1,2\}}\|_1+\|\ba^{\{1,3\}}\|_1+\|\ba^D\|_1)>0.
\end{align*}
Hence, by Corollary \ref{co:unstable2}, the 3D-MMRRW is unstable in our sense. 
\end{proof}

C1-7-1 on Table \ref{tab:table1} is a unique case, and in relation to semimartingale reflecting Brownian motions (SRBMs), it corresponds to the case where SRBMs have spiral fluid paths (see, for example, Kharroubi et al.\ \cite{Kharroubi02} and Bramson et al.\ \cite{Bramson10}). 
The following theorem provides the classification in C1-7-1; for ordinary three-dimensional reflecting random walks, a similar result was obtained in Malyshev and Menshikov \cite{Malyshev81} (also, see Fayolle et al.\ \cite{Fayolle95}). 
%
%%%%%%%%%%%%%%%%%%%%%%%%%%%%%
\begin{theorem} \label{th:3D_spiral}
Assume $(\{0\}\times\{0\}\times\{0\}\times S^\emptyset)\cap\calS_0\ne\emptyset$. Furthermore, assume the 3D-MMRRW satisfies $\ba^D<0$ and define $r$ as  
\begin{equation}
r=\left|\frac{a^{\{1,2\}}(2)}{a^{\{1,2\}}(1)}\right|\,\left|\frac{a^{\{2,3\}}(3)}{a^{\{2,3\}}(2)}\right|\,\left|\frac{a^{\{3,1\}}(1)}{a^{\{3,1\}}(3)}\right|. 
\label{eq:spiral_r}
\end{equation}
\begin{itemize} 
\item[(a)] Assume the 3D-MMRRW satisfies 
\begin{align*}
&a^{\{1,2\}}(1)<0,\ a^{\{1,2\}}(2)>0,\ a^{\{2,3\}}(2)<0,\ a^{\{2,3\}}(3)>0, 
a^{\{3,1\}}(3)<0,\ a^{\{3,1\}}(1)>0, 
\end{align*}
then the 3D-MMRRW is stable in our sense if $r<1$ and it is unstable in our sense if $r>1$.
\item[(b)] C1-7-1: Assume the 3D-MMRRW satisfies 
\begin{align*}
&a^{\{1,2\}}(1)>0,\ a^{\{1,2\}}(2)<0,\ a^{\{2,3\}}(2)>0,\ a^{\{2,3\}}(3)<0,
a^{\{3,1\}}(3)>0,\ a^{\{3,1\}}(1)<0, 
\end{align*}
then the 3D-MMRRW is stable in our sense if $r>1$ and it is unstable in our sense if $r<1$.
\end{itemize}
\end{theorem}

%%%%%%%%%%%%%%%%%%%%%%%
\begin{proof}
We prove only part (a); part (b) is analogously proved. 

Note that $(\{0\}\times\{0\}\times\{0\}\times S^\emptyset)\subset\calV_\emptyset$; hence, by the assumption, we have $\calV_\emptyset\cap\calS_0\ne\emptyset$. 
Consider part (a), where we have 
\[
\scrD_{stable}=\{\{1,2\},\{2,3\},\{3,1\},D\}.
\] 
Assume $r<1$ and let $\delta$ be a positive number satisfying $\sqrt{r}<\delta<1$. Set $\bw$ as 
\begin{equation}
\bw = \biggl( 1,\,\frac{-a^{\{1,2\}}(1)}{a^{\{1,2\}}(2)}\,\delta,\ \frac{a^{\{2,3\}}(2)}{a^{\{2,3\}}(3)}\,\frac{a^{\{1,2\}}(1)}{a^{\{1,2\}}(2)}\,\delta^2 \biggr)>0, 
\label{eq:3DMMRRW_w}
\end{equation}
then we have $\langle\ba^D,\bw\rangle<0$ and 
\begin{align*}
&\langle\ba^{\{1,2\}},\bw\rangle = a^{\{1,2\}}(1) (1-\delta) < 0,\cr
&\langle\ba^{\{2,3\}},\bw\rangle 
= -a^{\{2,3\}}(2)\frac{a^{\{1,2\}}(1)}{a^{\{1,2\}}(2)} \delta(1-\delta) < 0,\cr
&\langle\ba^{\{3,1\}},\bw\rangle 
= a^{\{3,1\}}(1) (1-\delta^2/r) < 0.
\end{align*}
Therefore, by Theorem \ref{th:stable}, the 3D-MMRRW is stable in our sense.

Next, we assume $r>1$ and let $\delta$ be a positive number satisfying $1<\delta<\sqrt{r}$. Set $m=3$ and set $A_1=\{1,2\}$, $A_2=\{2,3\}$ and $A_3=\{3,1\}$; this satisfies conditions (\ref{eq:Ak_condition1}) and (\ref{eq:Ak_condition2}). By expressions (\ref{eq:barDfstable}) and (\ref{eq:barIfA}), we obtain $\bar{\scrD}_{f,stable}=\scrD_{stable}$ and 
\[
\bar{\scrI}_{f,\{1,2\}} = \{1\},\quad 
\bar{\scrI}_{f,\{2,3\}} = \{2\},\quad 
\bar{\scrI}_{f,\{3,1\}} = \{3\},\quad 
\bar{\scrI}_{f,D} = \{1,2,3\}.
\]
Set $\bw_0$ as $\bw_0=\bw$, where $\bw$ is given by expression (\ref{eq:3DMMRRW_w}), and set $c_0$ as
\[
c_0 = \frac{|\langle\ba^D,\bw_0\rangle|}{\min_{k\in\{1,2,3\}}|a^D(k)|}\delta>0. 
\]
For $k\in\{1,2,3\}$, $\bw_k$ is given as $\bw_k=\bw_0-c_0\sum_{k'\in D\setminus A_k} \be_{k'}$ and $f_k$ as $f_k(\bx)=\langle\bx,\bw_k\rangle$ for $\bx\in\mathbb{R}^3$.  
Hence, we have
\begin{align*}
&f_1(\ba^{\{1,2\}}) = \langle\ba^{\{1,2\}},\bw_1\rangle = \langle\ba^{\{1,2\}},\bw\rangle > 0,\cr
&f_2(\ba^{\{2,3\}}) = \langle\ba^{\{2,3\}},\bw_2\rangle = \langle\ba^{\{2,3\}},\bw\rangle  > 0,\cr
&f_3(\ba^{\{3,1\}}) = \langle\ba^{\{3,1\}},\bw_3\rangle = \langle\ba^{\{3,1\}},\bw\rangle > 0.
\end{align*}
Furthermore, for $k\in\{1,2,3\}$, we have
\[
f_k(\ba^D) = \langle\ba^D,\bw_k\rangle 
= |\langle\ba^D,\bw_0\rangle| \biggl(-1+\frac{\sum_{k'\in D\setminus A_k} |a^D(k')|}{\min_{k'\in\{1,2,3\}}|a^D(k')|}\delta\biggr)>0.
\]
Therefore, by Corollary \ref{co:unstable2}, the 3D-MMRRW is unstable in our sense. 
\end{proof}

\begin{table}[htdp]
\caption{Classification of 3D-MMRRWs: $\ba^D=(-\,-\,-)$.}
\begin{center}
\begin{tabular}{cccccccc} \hline
Case & $\ba^{\{1,2\}}$ & $\ba^{\{2,3\}}$ & $\ba^{\{3,1\}}$ & $\ba^{\{1\}}$ & $\ba^{\{2\}}$ & $\ba^{\{3\}}$ & $\calL$ \cr \hline
C1-1-1 & $(-\,-\,0)$ & $(0\,-\,-)$ & $(-\,0\,-)$ & $(-\,0\,0)$ & $(0\,-\,0)$ & $(0\,0\,-)$ & stable \cr
C1-1-2 & $(-\,-\,0)$ & $(0\,-\,-)$ & $(-\,0\,-)$ & $(+\,0\,0)$ &  &  & unstable \cr
C1-1-3 & $(-\,-\,0)$ & $(0\,-\,-)$ & $(-\,0\,-)$ &  & $(0\,+\,0)$ &  & unstable \cr
C1-1-4 & $(-\,-\,0)$ & $(0\,-\,-)$ & $(-\,0\,-)$ &  &  & $(0\,0\,+)$ & unstable \cr \hline
C1-2-1 & $(+\,-\,0)$ & $(0\,-\,-)$ & $(-\,0\,-)$ & $(-\,0\,0)$ & NA & $(0\,0\,-)$ & stable \cr
C1-2-2 & $(+\,-\,0)$ & $(0\,-\,-)$ & $(-\,0\,-)$ & $(+\,0\,0)$ & NA &  & unstable \cr
C1-2-3 & $(+\,-\,0)$ & $(0\,-\,-)$ & $(-\,0\,-)$ &  & NA & $(0\,0\,+)$ & unstable \cr \hline
C1-3-1 & $(+\,+\,0)$ &  &  & NA & NA &  & unstable \cr \hline
C1-4-1 & $(+\,-\,0)$ & $(0\,+\,-)$ & $(-\,0\,-)$ & $(-\,0\,0)$ & NA & NA & stable \cr 
C1-4-2 & $(+\,-\,0)$ & $(0\,+\,-)$ & $(-\,0\,-)$ & $(+\,0\,0)$ & NA & NA & unstable \cr \hline
C1-5-1 & $(+\,-\,0)$ & $(0\,-\,+)$ & $(-\,0\,-)$ & $(-\,0\,0)$ & NA & $(0\,0\,-)$ & stable \cr
C1-5-2 & $(+\,-\,0)$ & $(0\,-\,+)$ & $(-\,0\,-)$ & $(+\,0\,0)$ & NA &  & unstable \cr
C1-5-3 & $(+\,-\,0)$ & $(0\,-\,+)$ & $(-\,0\,-)$ &  & NA & $(0\,0\,+)$ & unstable \cr \hline
C1-6-1 & $(+\,-\,0)$ & $(0\,+\,-)$ & $(+\,0\,-)$ & $(-\,0\,0)$ & NA & NA & stable \cr
C1-6-1 & $(+\,-\,0)$ & $(0\,+\,-)$ & $(+\,0\,-)$ & $(+\,0\,0)$ & NA & NA & unstable \cr \hline
C1-7-1 & $(+\,-\,0)$ & $(0\,+\,-)$ & $(-\,0\,+)$ & NA & NA & NA & (The.\ \ref{th:3D_spiral}) \cr \hline
\end{tabular}
\end{center}
\label{tab:table1}
\end{table}%

\begin{table}[htdp]
\caption{Classification of 3D-MMRRWs: $\ba^D=(+\,-\,-)$.}
\begin{center}
\begin{tabular}{cccccccc} \hline
Case & $\ba^{\{1,2\}}$ & $\ba^{\{2,3\}}$ & $\ba^{\{3,1\}}$ & $\ba^{\{1\}}$ & $\ba^{\{2\}}$ & $\ba^{\{3\}}$ & $\calL$ \cr \hline
C2-1-1 & $(-\,-\,0)$ & NA & $(-\,0\,-)$ & $(-\,0\,0)$ & $(0\,-\,0)$ & $(0\,0\,-)$ & stable \cr
C2-1-2 & $(-\,-\,0)$ & NA & $(-\,0\,-)$ & $(+\,0\,0)$ &  &  & unstable \cr
C2-1-3 & $(-\,-\,0)$ & NA & $(-\,0\,-)$ &  & $(0\,+\,0)$ &  & unstable \cr
C2-1-4 & $(-\,-\,0)$ & NA & $(-\,0\,-)$ &  &  & $(0\,0\,+)$ & unstable \cr \hline
C2-2-1 & $(+\,-\,0)$ & NA & $(-\,0\,-)$ & $(-\,0\,0)$ & NA & $(0\,0\,-)$ & stable \cr
C2-2-2 & $(+\,-\,0)$ & NA & $(-\,0\,-)$ & $(+\,0\,0)$ & NA &  & unstable \cr
C2-2-3 & $(+\,-\,0)$ & NA & $(-\,0\,-)$ &  & NA & $(0\,0\,+)$ & unstable \cr \hline
C2-3-1 & $(+\,+\,0)$ & NA &  & NA & NA &  & unstable \cr \hline
C2-4-1 & $(+\,-\,0)$ & NA & $(+\,0\,-)$ & $(-\,0\,0)$ & NA & NA & stable \cr
C2-4-2 & $(+\,-\,0)$ & NA & $(+\,0\,-)$ & $(+\,0\,0)$ & NA & NA & unstable \cr \hline
C2-5-1 & $(+\,-\,0)$ & NA & $(-\,0\,+)$ & NA & NA & $(0\,0\,-)$ & stable \cr
C2-5-2 & $(+\,-\,0)$ & NA & $(-\,0\,+)$ & NA & NA & $(0\,0\,+)$ & unstable \cr \hline
\end{tabular}
\end{center}
\label{tab:table2}
\end{table}%

\begin{table}[htdp]
\caption{Classification of 3D-MMRRWs: $\ba^D=(+\,+\,-)$.}
\begin{center}
\begin{tabular}{cccccccc} \hline
Case & $\ba^{\{1,2\}}$ & $\ba^{\{2,3\}}$ & $\ba^{\{3,1\}}$ & $\ba^{\{1\}}$ & $\ba^{\{2\}}$ & $\ba^{\{3\}}$ & $\calL$ \cr \hline
C3-1-1 & $(-\,-\,0)$ & NA & NA & $(-\,0\,0)$ & $(0\,-\,0)$ & NA & stable \cr
C3-1-2 & $(-\,-\,0)$ & NA & NA & $(+\,0\,0)$ &  & NA & unstable \cr
C3-1-3 & $(-\,-\,0)$ & NA & NA &  & $(0\,+\,0)$ & NA & unstable \cr \hline
C3-2-1 & $(+\,-\,0)$ & NA & NA & $(-\,0\,0)$ & NA & NA & stable \cr
C3-2-2 & $(+\,-\,0)$ & NA & NA & $(+\,0\,0)$ & NA & NA & unstable \cr \hline
C3-3-1 & $(+\,+\,0)$ & NA & NA & NA & NA & NA & unstable \cr \hline
\end{tabular}
\end{center}
\label{tab:table3}
\end{table}%

\begin{table}[htdp]
\caption{Classification of 3D-MMRRWs: $\ba^D=(+\,+\,+)$.}
\begin{center}
\begin{tabular}{cccccccc} \hline
Case & $\ba^{\{1,2\}}$ & $\ba^{\{2,3\}}$ & $\ba^{\{3,1\}}$ & $\ba^{\{1\}}$ & $\ba^{\{2\}}$ & $\ba^{\{3\}}$ & $\calL$ \cr \hline
C4-1-1 & NA & NA & NA & NA & NA & NA & unstable \cr \hline
\end{tabular}
\end{center}
\label{tab:table4}
\end{table}%

%%%%%%%%%%%%%%%%%%%%%%%
%
% Section 4
%
%%%%%%%%%%%%%%%%%%%%%%%
%
\section{Application to queueing networks} \label{sec:application_QNW}

%%%%%%%%%%%%%%%%%%%%%%%%%%%%
%
\subsection{Relation between queueing networks and MMRRWs}

\subsubsection{Multiclass queueing networks}

In this section, we consider a multiclass open queueing network whose stochastic behavior is represented as a CTMC. There are $d$ customer classes and each customer class has its own queue. There are a number of stations and each queue (customer class) belongs to one of the stations. Each station has several servers and they serve customers in the station according to some service discipline. Let $D$ be defined as $D=\{1,2,...,d\}$, which is the index set of customer classes and also that of queues. 
Let $l$ and $l'$ be elements in $D$. A class-$l$ customer arriving from the outside of the network joins queue $l$ (Q$_l$) and receives service there; after completing the service, the customer leaves the network with probability $r_{l,0}\ge 0$ or changes the service class from $l$ to $l'$ with probability $r_{l,l'}\ge 0$, joins queue $l'$ (Q$_{l'}$) and receives service there again. Routing probabilities satisfy $\sum_{l'\in D\cup\{0\}}r_{l,l'}=1$ for all $l\in D$ and $r_{l,0}\ne 0$ for some $l\in D$. We denote by $R$ the routing probability matrix, i.e., $R=(r_{l,l'},i,i'\in D)$. We do not consider batch arrivals and batch services. 

For $l\in D$, let $\bar{X}_t(l)$ be the number of customers in Q$_l$ at time $t$ and $\bar{\bX}_t$ the vector of $\bar{X}_t(l)$'s, i.e., $\bar{\bX}_t=(\bar{X}_t(l),l\in D)$. Let $\bar{J}_t$ be the supplementary (background) state at time $t$ and assume the process $\{\bar{\bY}_t\}=\{(\bar{\bX}_t,\bar{J}_t)\}$ is a semi-irreducible CTMC on the state space $\calS$. 
For $l,l'\in D$, let $\bar{E}_t(l)$ be the cumulative number of class-$l$ customers arriving from the outside until time $t$, $\bar{O}_t(l)$ the cumulative number of customers departing from Q$_l$ until time $t$ and $\Phi_{l,l'}(n)$ the number of first $n$ customers departing from Q$_l$ that are next routed to Q$_{l'}$. We assume that $\bar{E}_0(l)=0$, $\bar{O}_0(l)=0$ and $\Phi_{l,l'}(0)=0$ for $l,l'\in D$. Then, we have the following queueing network equation (see, for example, Bramson \cite{Bramson08}): 
\begin{equation}
\bar{X}_t(l)=\bar{X}_0(l)+\bar{E}_t(l)+\sum_{l'\in D} \Phi_{l',l}(\bar{O}_t(l'))-\bar{O}_t(l),\ l\in D.
\label{eq:QNeq}
\end{equation}
We denote by $\bar{\bE}_t$ the vector of $\bar{E}_t(l)$'s and by $\bar{\bO}_t$ that of $\bar{O}_t(l)$'s.
Denote by $\bar{\ba}$ the mean drift vector of the process $\{\bX_t\}$ defined as 
\[
\bar{\ba} = \lim_{t\to\infty} \frac{1}{t} (\bar{\bX}_t-\bar{\bX}_0)  
\]
if the limit exists with probability one. Furthermore, we assume the following limits exist: 
\begin{align*}
&\bar{\blambda} = \lim_{t\to\infty} \bar{\bE}_t/t,\quad a.s., \qquad
\bar{\bmu} = \lim_{t\to\infty} \bar{\bO}_t/t,\quad a.s., 
\end{align*}
where, for $l\in D$, $\bar{\lambda}(l)$, the $l$-th element of $\bar{\blambda}$, is the external arrival rate of Q$_l$ and $\bar{\mu}(l)$ is the departure rate of Q$_l$. By expression (\ref{eq:QNeq}), we obtain 
\begin{equation*}
\bar{\ba} = \bar{\blambda} + (R^\top-I) \bar{\bmu}, 
\end{equation*}
where we use the fact that $\lim_{t\to\infty}\Phi_{l',l}(\bar{O}_t(l'))/t=r_{l',l}\,\bar{\mu}(l')$ for $l,l'\in D$. This equation asserts that the mean drift vector is given by the external arrival rates and the departure rates of the queues in the network. 
Note that, if the CTMC $\{\bar{\bY}_t\}$ is stable in our sense, the mean drift vector is a zero vector, i.e., $\bar{\ba}=0$, and we obtain the traffic equation $\bar{\bmu} = (I-R^\top)^{-1} \bar{\blambda}$.

Next, assuming the CTMC $\{\bar{\bY}_t\}$ is a continuous-time version of a $d$-dimensional MMRRW, we obtain a MMRRW corresponding to $\{\bar{\bY}_t\}$ in a manner similar to that used in Example \ref{ex:reentrant1}. For $\by,\by'\in\calS$ such that $\by\ne\by'$, we denote by $q(\by,\by')$ the transition rate that the state of the CTMC changes from $\by$ to $\by'$; for $\by\in\calS$, we define $q(\by,\by)$ as $q(\by,\by)=-\sum_{\by'\ne\by}q(\by,\by')$. 
We assume that $q(\by,\by),\,\by\in\calS,$ are bounded and set $\nu$ at a positive number satisfying $\nu\ge\sup_{\by\in\calS}|q(\by,\by)|$. Then, by uniformization, we can obtain a discrete-time Markov chain $\calL=\{\bY_n\}=\{(\bX_n,J_n)\}$ on the state space $\calS$ that has the same stationary distribution as the original CTMC, if it exists. 
The transition probabilities of $\{\bY_n\}$ are given by $\mathbb{P}(\bY_1=\by'|\bY_0=\by)=\delta_{\by,\by'}+q(\by,\by')/\nu$ for $\by,\by'\in\calS$, where $\delta_{\cdot,\cdot}$ is the delta function. 
Since the CTMC $\{\bar{\bY}_t\}$ is a continuous-time version of a MMRRW, the Markov chain $\{\bY_n\}$ becomes a $d$-dimensional MMRRW on the state space $\calS$; it is also the discrete-time queueing network corresponding to the original queueing network.  The mean drift vector of the MMRRW $\{\bY_n\}$ is given by 
\[
\ba = \lim_{n\to\infty} \frac{1}{n}(\bX_n-\bX_0) = \bar{\ba}/\nu,\quad a.s.
\]

%%%%%%%%%%%%%%%%%%%%%%%%%%%%
%
\subsubsection{Expanded and induced Markov chains}

First, we note that expanded and induced Markov chains can be defined for the CTMC $\{\bar{\bY}_t\}=\{(\bar{\bX}_t,\bar{J}_t)\}$ in a manner similar to that used for discrete-time MMRRWs. Let $A$ be a nonempty element of $\scrP(D)$ and $\{\bar{\bY}^A_t\}=\{(\bar{\bX}^A_t,\bar{J}^A_t)\}$ the expanded Markov chain of $\{\bar{\bY}_t\}$ with index $A$. 
We consider $\{\bar{\bY}^A_t\}$ represents a queueing network with extraordinary queues; that is, for $l\in A$, queue Q$_l$ is saturated with customers and $\bar{X}^A_t(l)$ is the number of arrivals at Q$_l$ minus that of departures from Q$_l$ until time $t$ (a similar interpretation has been utilized in Dumas \cite{Dumas97} and Tezcan \cite{Tezcan13}). For $l\in A$, $\bar{X}^A_t(l)$ may, therefore, take negative integers. For $l\in D\setminus A$, Q$_l$ is an ordinary queue. $\{\bar{\bY}^A_t\}$ satisfies the following queueing network equation:
\begin{equation}
\bar{X}^A_t(l)=\bar{X}^A_0(l)+\bar{E}_t(l)+\sum_{l'\in D} \Phi_{l',l}(\bar{O}^A_t(l'))-\bar{O}^A_t(l),\ l\in D, 
\label{eq:QNeq_expanded}
\end{equation}
where, for $l\in A$, $\bar{O}^A_t(l)$ is the cumulative number of departures from the saturated queue with index $l$ until time $t$. Note that, since we consider the multiclass queueing network, the output process of each saturated queue may depend on the state of other queues in the same station. 
Let $\bar{\bmu}^A$ be the departure rate vector of $\{\bar{\bY}^A_t\}$. Then, the mean drift vector of $\{\bar{\bY}^A_t\}$, $\bar{\ba}^A$, is given by 
\[
\bar{\ba}^A = \lim_{t\to\infty} \frac{1}{t}(\bar{\bX}^A_t-\bar{\bX}^A_0) = \bar{\blambda} + (R^\top-I) \bar{\bmu}^A,\quad a.s. 
\]
Induced Markov chain $\{(\bar{\bX}^A_t(D\setminus A),\bar{J}^A_t)\}$ is a queueing network obtained by replacing each extraordinary queue with several external arrival processes. The queueing network equation of the induced Markov chain $\{(\bar{\bX}^A_t(D\setminus A),\bar{J}^A_t)\}$ is, therefore, given by 
\begin{equation}
\bar{X}^A_t(l)=\bar{X}^A_0(l)+\bar{E}_t(l)+\sum_{l'\in A} \Phi_{l',l}(\bar{O}^A_t(l'))+\sum_{l'\in D\setminus A} \Phi_{l',l}(\bar{O}^A_t(l'))-\bar{O}^A_t(l),\ l\in D\setminus A, 
\label{eq:QNeq_induced}
\end{equation}
where $\bar{E}_t(l)+\sum_{l'\in A} \Phi_{l',l}(\bar{O}^A_t(l'))$ is considered as the number of external arrivals at Q$_l$ until time $t$ in the queueing network corresponding to the induced Markov chain. From this queueing network equation, we obtain 
\[
\bar{\ba}^A(D\setminus A) = \bar{\blambda}(D\setminus A) + R_0^\top\bar{\bmu}^A(A) + (R_1^\top-I) \bar{\bmu}^A(D\setminus A), 
\]
where $R_0=\big(r_{l,l'},l\in A,l'\in D\setminus A\big)$ and $R_1=\big(r_{l,l'},l,l'\in D\setminus A\big)$. 
If the induced Markov chain $\{(\bar{\bX}^A_t(D\setminus A),\bar{J}^A_t)\}$ is semi-irreducible and stable in our sense, then $\bar{\ba}^A(D\setminus A)=0$ and we obtain the following traffic equation:
\begin{equation}
\bar{\bmu}^A(D\setminus A) = (I-R_1^\top)^{-1} \bigl( \bar{\blambda}(D\setminus A) + R_0^\top\bar{\bmu}^A(A) \bigr).
\end{equation}
Thus, we obtain 
\begin{align}
\bar{\ba}^A(A) &= \bar{\blambda}(A) + R_2^\top\bar{\bmu}^A(D\setminus A) + (R_3^\top-I) \bar{\bmu}^A(A) \cr
&= \bar{\blambda}(A) + R_2^\top(I-R_1^\top)^{-1}\bar{\blambda}(D\setminus A) + \big( R_2^\top(I-R_1^\top)^{-1}R_0^\top +R_3^\top-I \big) \bar{\bmu}^A(A), 
\label{eq:aAA_nw}
\end{align}
where $R_2=\big(r_{l,l'},l\in D\setminus A,l'\in A\big)$ and $R_3=\big(r_{l,l'},l,l'\in A\big)$. 

For nonempty $A\in\scrP(D)$, let $\hat{\calL}^A=\{\bY^A_n\}=\{(\bX^A_n,J^A_n)\}$ be the discrete-time Markov chain obtained from $\{(\bar{\bX}^A_t,\bar{J}^A_t)\}$ by uniformization with parameter $\nu$. Then, $\hat{\calL}^A$ is the expanded Markov chain of $\calL=\{\bY_n\}=\{(\bX_n,J_n)\}$ with index $A$, where $\calL$ is the discrete-time Markov chain obtained from $\{(\bar{\bX}_t,\bar{J}_t)\}$ by uniformization with the parameter $\nu$. $\calL^A=\{(\bX^A_n(D\setminus A),J^A_t)\}$  is the discrete-time Markov chain corresponding to $\{(\bar{\bX}^A_t(D\setminus A),\bar{J}^A_t)\}$ and it is an induced Markov chain of $\calL$. 
Therefore, letting $\ba^A$ be the mean drift vector of $\hat{\calL}^A$, we have 
\begin{align}
\ba^A = \lim_{n\to\infty} \frac{1}{n} (\bX^A_n-\bX^A_0) = \bar{\ba}^A/\nu,\quad a.s., 
\label{eq:aA_nw}
\end{align}
and if the induced Markov chain $\calL^A$ is semi-irreducible and stable in our sense, we have $\ba^A(D\setminus A)=0$ and $\ba^A(A)=\bar{\ba}^A(A)/\nu$, where $\bar{\ba}^A(A)$ is given by expression (\ref{eq:aAA_nw}). This implies that, {\it in stability analysis of queueing networks, we do not need the stationary distributions of the stable induced Markov chains to get the mean drift vectors if the departure rates of the saturated queues can be obtained by other methods.} In such a case, it is easy to apply our results to queueing networks.

%%%%%%%%%%%%%%%%%%%%%%%%%%%%
%
\subsection{A two-station network}
\subsubsection{Model description}

\begin{figure}[bht]
\begin{center}
\setlength{\unitlength}{0.8mm}
\begin{picture}(100,42)(0,0)
%
% Station 1
\thicklines
\put(2,15){\makebox(0,0){\normalsize $\bar{\lambda}_1$}}
\put(5,15){\vector(1,0){7}}

\put(14,6){\makebox(0,0){\normalsize Q$_1$}}
\put(27,11){\makebox(0,0){\normalsize $\bar{h}_1$}}
\put(12,20){\line(1,0){6}}
\put(12,10){\line(1,0){6}}
\put(18,10){\line(0,1){10}}
\multiput(22,15)(2,0){6}{\line(1,0){1}}
\put(36,15){\vector(1,0){15}}

\put(42,40){\makebox(0,0){\normalsize Q$_4$}}
%\put(22,33){\makebox(0,0){\small High priority}}
\put(27,26){\makebox(0,0){\normalsize $\bar{h}_4$}}
\put(37,25){\line(1,0){6}}
\put(37,35){\line(1,0){6}}
\put(37,25){\line(0,1){10}}
\multiput(22,30)(2,0){6}{\line(1,0){1}}
\put(19,30){\vector(-1,0){15}}

\put(28,0){\makebox(0,0){\normalsize Station 1}}
\put(20,5){\line(0,1){35}}
\put(20,5){\line(1,0){15}}
\put(20,40){\line(1,0){15}}
\put(35,5){\line(0,1){35}}
%
% Station 2
\thicklines
\put(54,6){\makebox(0,0){\normalsize Q$_2$}}
%\put(75,12){\makebox(0,0){\small High priority}}
\put(67,18){\makebox(0,0){\normalsize $\bar{h}_2$}}
\put(52,20){\line(1,0){6}}
\put(52,10){\line(1,0){6}}
\put(58,10){\line(0,1){10}}
\multiput(62,15)(2,0){6}{\line(1,0){1}}
\put(88,20){\makebox(0,0){\normalsize $p$}}
\put(99,11){\makebox(0,0){\normalsize $1-p$}}
\put(76,15){\vector(1,0){28}}
\put(91,15){\line(0,1){13}}
\put(91,28){\vector(-1,0){8}}

\put(99,32){\makebox(0,0){\normalsize $\bar{\lambda}_3$}}
\put(95,32){\vector(-1,0){12}}

\put(82,40){\makebox(0,0){\normalsize Q$_3$}}
\put(67,33){\makebox(0,0){\normalsize $\bar{h}_3$}}
\put(77,25){\line(1,0){6}}
\put(77,35){\line(1,0){6}}
\put(77,25){\line(0,1){10}}
\multiput(62,30)(2,0){6}{\line(1,0){1}}
\put(59,30){\vector(-1,0){15}}

\put(68,0){\makebox(0,0){\normalsize Station 2}}
\put(60,5){\line(0,1){35}}
\put(60,5){\line(1,0){15}}
\put(60,40){\line(1,0){15}}
\put(75,5){\line(0,1){35}}

\end{picture}
\caption{A two-station network}
\label{fig:twostation}
\end{center}
\end{figure}
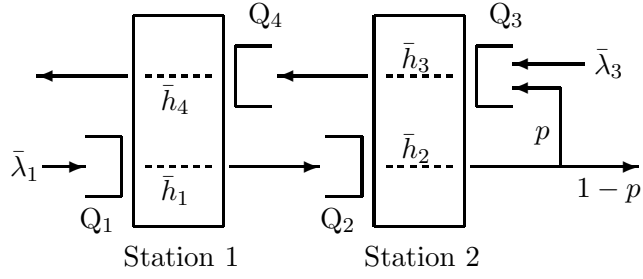

As an example of multiclass queueing network, we consider the following two-station network depicted in Fig.~\ref{fig:twostation}, which includes as special cases Lu-Kumar network \cite{Lu91}, Kumar-Seidman network \cite{Kumar90} and Rybko-Stolyar network \cite{Rybko92}  (also see Bramson \cite{Bramson08}). 

%\medskip
{\it Network configuration.}\quad 
Exogenous customers arrive at queue 1 (Q$_1$) as class-1 customers and at queue 3 (Q$_3$) as class-3 customers. After completing service at station~1, customers in Q$_1$ move to station~2 and join queue 2 (Q$_2$) as class-2 customers; customers in Q$_2$ next reenter station~2 and join Q$_3$ as class-3 customers with probability $p$ or depart from the network with probability $1-p$. After completing service at station~2, customers in Q$_3$ move to station~1 and join queue 4 (Q$_4$) as class-4 customers; customers in Q$_4$ next depart from the network. 
Let $D$ be defined as $D=\{1,2,3,4\}$ and, for $l\in D$, denote by $\bar{X}_t(l)$ the number of customers in Q$_l$ including one being served at time $t$; $\bar{\bX}_t$ is the vector of $\bar{X}_t(l)$'s, i.e., $\bar{\bX}_t=(\bar{X}_t(1),\bar{X}_t(2),\bar{X}_t(3),\bar{X}_t(4))$.

%%%%%%%%%%%%%%%%%%%%%%%%%%%%%%%%%%%%%%%%%%%%%%%%%%%%%%%%%%%
%
% fig_twostation.tex
%
%%%%%%%%%%%%%%%%%%%%%%%%%%%%%%%%%%%%%%%%%%%%%%%%%%%%%%%%%%%
%%%%%%%%%%%%%%%%%%%%%%%%%%%%%%%%%%%%%%%%%%%%%%%%%%%%%%%%%%%

%\medskip
{\it Arrival processes.}\quad 
The exogenous arrival process of class-1 customers and that of class-3 customers are subject to independent Markovian arrival processes (MAPs for short); MAPs are tractable arrival processes and they can represent correlated interarrival times. 
For $i=1,3$, we denote by $(\bar{C}_i,\bar{D}_i)$ the representation of the MAP for class-$i$ customers, where $\bar{C}_i$ is a matrix of phase transition rates without arrivals and $\bar{D}_i$ is that of phase transition rates with arrivals (see, for example, Latouche and Ramaswami \cite{Latouche99}).  
We assume that the phase set of the MAP for class-$i$ customers is given by $S^a_i=\{1,2,...,s^a_i\}$, where $s^a_i$ is some positive integer.
$\bar{C}_i+\bar{D}_i$ is the infinitesimal generator of the Markov chain that represents the phase process of the MAP for class-$i$ customers. We assume $\bar{C}_i+\bar{D}_i$ to be irreducible. Then, $\bar{C}_i+\bar{D}_i$ is positive recurrent and we denote by $\bpi_i^*=(\pi_{ij}^*,j\in S^a_i)$ the stationary phase distribution. The mean arrival rate of the MAP for class-$i$ customers, denoted by $\bar{\lambda}_i$, is given as $\bar{\lambda}_i=\bpi_i^* \bar{D}_i\bone$, where $\bone$ is a column vector of $1$'s whose dimension is determined in context. 
We denote by $\bar{J}_t(3)$ the phase state of the MAP for class-$1$ customers at time $t$ and by $\bar{J}_t(4)$ that of the MAP for class-$3$ customers. In order to exclude trivial cases, we assume $\bar{\lambda}_1\ne 0$. 

\begin{figure}[bht]
\begin{center}
\setlength{\unitlength}{0.8mm}
\begin{picture}(130,80)(0,0)
%
%%%%%
\thicklines
\put(10,40){\makebox(0,0){\normalsize $(0,0)$}}
\put(8,34){\makebox(0,0){\small $\bar{T}_i^{00}$}}

%%%%%
\put(40,50){\makebox(0,0){\normalsize $(1,0)$}}
\put(38,44){\makebox(0,0){\small $\bar{T}_i^{+0}$}}
\put(33,50){\vector(-2,-1){15}}
\put(24,49){\makebox(0,0){\small $\bar{T}_i^{1^*0}$}}

\put(40,30){\makebox(0,0){\normalsize $(0,1)$}}
\put(38,24){\makebox(0,0){\small $\bar{T}_i^{0+}$}}
\put(33,30){\vector(-2,1){15}}
\put(30,36){\makebox(0,0){\small $\bar{T}_i^{01^*}$}}

%%%%%
\put(70,60){\makebox(0,0){\normalsize $(2,0)$}}
\put(68,54){\makebox(0,0){\small $\bar{T}_i^{+0}$}}
\put(63,60){\vector(-2,-1){15}}
\put(54,59){\makebox(0,0){\small $\bar{T}_i^{2^*0}$}}

\put(70,40){\makebox(0,0){\normalsize $(1,1)$}}
\put(68,34){\makebox(0,0){\small $\bar{T}_i^{++}$}}
\put(63,40){\vector(-2,1){15}}
\put(60,46){\makebox(0,0){\small $\bar{T}_i^{+1^*}$}}
\put(63,40){\vector(-2,-1){15}}
\put(54,39){\makebox(0,0){\small $\bar{T}_i^{1^*+}$}}

\put(70,20){\makebox(0,0){\normalsize $(0,2)$}}
\put(68,14){\makebox(0,0){\small $\bar{T}_i^{0+}$}}
\put(63,20){\vector(-2,1){15}}
\put(60,26){\makebox(0,0){\small $\bar{T}_i^{02^*}$}}

%%%%%
\put(100,70){\makebox(0,0){\normalsize $(3,0)$}}
\put(98,64){\makebox(0,0){\small $\bar{T}_i^{++}$}}
\put(93,70){\vector(-2,-1){15}}
\put(84,69){\makebox(0,0){\small $\bar{T}_i^{2^*0}$}}

\put(100,50){\makebox(0,0){\normalsize $(2,1)$}}
\put(98,44){\makebox(0,0){\small $\bar{T}_i^{++}$}}
\put(93,50){\vector(-2,1){15}}
\put(90,56){\makebox(0,0){\small $\bar{T}_i^{+1^*}$}}
\put(93,50){\vector(-2,-1){15}}
\put(84,49){\makebox(0,0){\small $\bar{T}_i^{2^*+}$}}

\put(100,30){\makebox(0,0){\normalsize $(1,2)$}}
\put(98,24){\makebox(0,0){\small $\bar{T}_i^{++}$}}
\put(93,30){\vector(-2,1){15}}
\put(90,36){\makebox(0,0){\small $\bar{T}_i^{+2^*}$}}
\put(93,30){\vector(-2,-1){15}}
\put(84,29){\makebox(0,0){\small $\bar{T}_i^{1^*+}$}}

\put(100,10){\makebox(0,0){\normalsize $(0,3)$}}
\put(98,4){\makebox(0,0){\small $\bar{T}_i^{0+}$}}
\put(93,10){\vector(-2,1){15}}
\put(90,16){\makebox(0,0){\small $\bar{T}_i^{02^*}$}}

%%%%%
\put(115,70){\makebox(0,0){\normalsize $\cdots$}}
\put(115,50){\makebox(0,0){\normalsize $\cdots$}}
\put(115,30){\makebox(0,0){\normalsize $\cdots$}}
\put(115,10){\makebox(0,0){\normalsize $\cdots$}}
\end{picture}
\caption{Phase transition rates of the two-class MSP at station $i$.}
\label{fig:MSP1}
\end{center}
\end{figure}
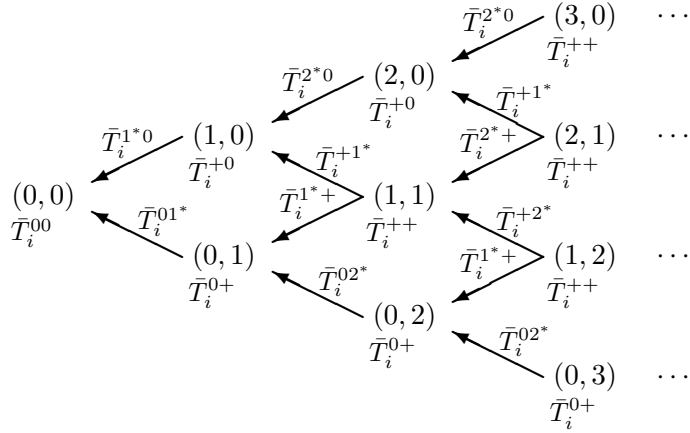

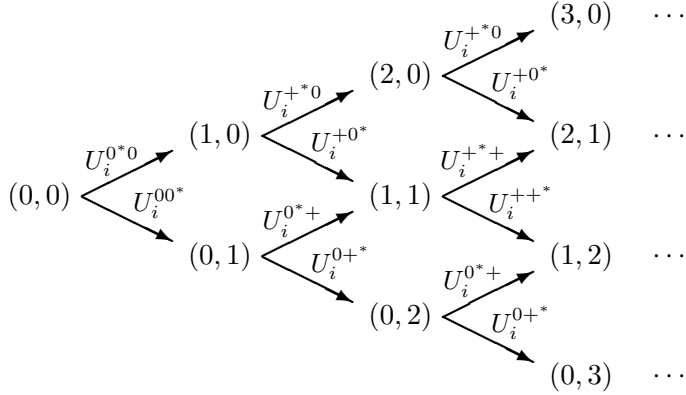
\begin{figure}[bht]
\begin{center}
\setlength{\unitlength}{0.8mm}
\begin{picture}(130,70)(0,0)
%
%%%%%
\thicklines
\put(10,35){\makebox(0,0){\normalsize $(0,0)$}}
\put(17,35){\vector(2,1){15}}
\put(22,41){\makebox(0,0){\small $U_i^{0^*0}$}}
\put(17,35){\vector(2,-1){15}}
\put(30,34){\makebox(0,0){\small $U_i^{00^*}$}}

%%%%%
\put(40,45){\makebox(0,0){\normalsize $(1,0)$}}
\put(47,45){\vector(2,1){15}}
\put(52,51){\makebox(0,0){\small $U_i^{+^*0}$}}
\put(47,45){\vector(2,-1){15}}
\put(60,44){\makebox(0,0){\small $U_i^{+0^*}$}}

\put(40,25){\makebox(0,0){\normalsize $(0,1)$}}
\put(47,25){\vector(2,1){15}}
\put(52,31){\makebox(0,0){\small $U_i^{0^*+}$}}
\put(47,25){\vector(2,-1){15}}
\put(60,24){\makebox(0,0){\small $U_i^{0+^*}$}}

%%%%%
\put(70,55){\makebox(0,0){\normalsize $(2,0)$}}
\put(77,55){\vector(2,1){15}}
\put(82,61){\makebox(0,0){\small $U_i^{+^*0}$}}
\put(77,55){\vector(2,-1){15}}
\put(90,54){\makebox(0,0){\small $U_i^{+0^*}$}}

\put(70,35){\makebox(0,0){\normalsize $(1,1)$}}
\put(77,35){\vector(2,1){15}}
\put(82,41){\makebox(0,0){\small $U_i^{+^*+}$}}
\put(77,35){\vector(2,-1){15}}
\put(90,34){\makebox(0,0){\small $U_i^{++^*}$}}

\put(70,15){\makebox(0,0){\normalsize $(0,2)$}}
\put(77,15){\vector(2,1){15}}
\put(82,21){\makebox(0,0){\small $U_i^{0^*+}$}}
\put(77,15){\vector(2,-1){15}}
\put(90,14){\makebox(0,0){\small $U_i^{0+^*}$}}

%%%%%
\put(100,65){\makebox(0,0){\normalsize $(3,0)$}}
\put(100,45){\makebox(0,0){\normalsize $(2,1)$}}
\put(100,25){\makebox(0,0){\normalsize $(1,2)$}}
\put(100,5){\makebox(0,0){\normalsize $(0,3)$}}

%%%%%
\put(115,65){\makebox(0,0){\normalsize $\cdots$}}
\put(115,45){\makebox(0,0){\normalsize $\cdots$}}
\put(115,25){\makebox(0,0){\normalsize $\cdots$}}
\put(115,5){\makebox(0,0){\normalsize $\cdots$}}
\end{picture}
\caption{Phase transition probabilities of the two-class MSP at station $i$.}
\label{fig:MSP2}
\end{center}
\end{figure}

%\medskip
{\it Service processes.}\quad 
We assume that there is a single server in each station and the service process there is represented as a two-class Markovian service process (MSP for short; see, for example, Ozawa \cite{Ozawa04}). 
While there are no standard representations for two-class MSPs, we use the following one, by which a preemptive-resume priority service, non-preemptive priority service and $(1,K)$-limited service can be described. 
For $i=1,2$, we denote by $S^s_i$ the phase set of the service process in station~$i$ and by $\bar{J}_t(i)$ the phase of the service process at time $t$; we assume $S^s_i$ is finite. 
Furthermore, we assume that $S^s_1$ is composed of mutually disjoint subsets $S^s_{10}$, $S^s_{11}$ and $S^s_{14}$ and that $S^s_2$ is also composed of mutually disjoint subsets $S^s_{20}$, $S^s_{22}$ and $S^s_{23}$, i.e., 
\[
S^s_1 = S^s_{10} \cup S^s_{11} \cup S^s_{14},\quad 
S^s_2 = S^s_{20} \cup S^s_{22} \cup S^s_{23}.
\]
When $\bar{J}_t(1)\in S^s_{10}$ ($\bar{J}_t(2)\in S^s_{20}$), the server in station~1 (resp.\ station~2) is idle or engaging in work other than service for customers; when $\bar{J}_t(1)\in S^s_{11}$ ($\bar{J}_t(2)\in S^s_{22}$), it is engaging in service for a class-1 (resp.\ class-2) customer; when $\bar{J}_t(1)\in S^s_{14}$ ($\bar{J}_t(2)\in S^s_{23}$), it is engaging in service for a class-4 (resp.\ class-3) customer. 
For $i\in\{1,2\}$, we assume the phase transition rates of the service process $\{\bar{J}_t(i)\}$ are given as Fig.~\ref{fig:MSP1}; in the figure, if $i=1$, then ``$(k,l)$" indicates the set of the phase states of $\{\bar{J}_t(1)\}$ in which $\bar{X}_t(1)=k$ and $\bar{X}_t(4)=l$, and if $i=2$, then it indicates that of the phase states of $\{\bar{J}_t(2)\}$ in which $\bar{X}_t(2)=k$ and $\bar{X}_t(3)=l$. 
Here, we explain only the case of $i=1$. 
$\bar{T}_1^{00}$ is the $|S^s_1|$-dimensional square matrix of transition rates for $\{\bar{J}_t(1)\}$ without service completions when $\bar{X}_t(1)=0$ and $\bar{X}_t(4)=0$; $\bar{T}_1^{+0}$ is that when $\bar{X}_t(1)\ge 1$ and $\bar{X}_t(4)=0$, and $\bar{T}_1^{0+}$ is that when $\bar{X}_t(1)=0$ and $\bar{X}_t(4)\ge 1$; $\bar{T}_1^{++}$ is that when $\bar{X}_t(1)\ge 1$ and $\bar{X}_t(4)\ge 1$. 
$\bar{T}_1^{1^*0}$ is the $|S^s_1|$-dimensional square matrix of transition rates for $\{\bar{J}_t(1)\}$ with a service completion of class-1 customer when $\bar{X}_t(1)=1$ and $\bar{X}_t(4)=0$, and $\bar{T}_1^{2^*0}$ is that when $\bar{X}_t(1)\ge 2$ and $\bar{X}_t(4)=0$; $\bar{T}_1^{1^*+}$ is that when $\bar{X}_t(1)=1$ and $\bar{X}_t(4)\ge 1$, and $\bar{T}_1^{2^*+}$ is that when $\bar{X}_t(1)\ge 2$ and $\bar{X}_t(4)\ge 1$. 
$\bar{T}_1^{01^*}$ is the $|S^s_1|$-dimensional square matrix of transition rates for $\{\bar{J}_t(1)\}$ with a service completion of class-4 customer when $\bar{X}_t(1)=0$ and $\bar{X}_t(4)=1$, and $\bar{T}_1^{02^*}$ is that when $\bar{X}_t(1)=0$ and $\bar{X}_t(4)\ge 2$; $\bar{T}_1^{+1^*}$ is that when $\bar{X}_t(1)\ge 1$ and $\bar{X}_t(4)=1$, and $\bar{T}_1^{+2^*}$ is that when $\bar{X}_t(1)\ge 1$ and $\bar{X}_t(4)\ge 2$. 
The phase of the service process in each station may change when a customer arrives at the station; for $i\in\{1,2\}$, we assume the phase transition probabilities of the service process $\{\bar{J}_t(i)\}$ at such arrival epochs are given as Fig.~\ref{fig:MSP2}. Here, we explain only the case of $i=1$. 
$U_1^{0^*0}$ is the transition probability matrix for $\{\bar{J}_t(1)\}$ when $\bar{X}_t(1)=0$, $\bar{X}_t(4)=0$ and a class-1 customer arrived at time $t$, and $U_1^{00^*}$ is that when $\bar{X}_t(1)=0$, $\bar{X}_t(4)=0$ and a class-4 customer arrived; $U_1^{+^*0}$ is that when $\bar{X}_t(1)\ge 1$, $\bar{X}_t(4)=0$ and a class-1 customer arrived at time $t$, and $U_1^{+0^*}$ is that when $\bar{X}_t(1)\ge 1$, $\bar{X}_t(4)=0$ and a class-4 customer arrived. 
$U_1^{0^*+}$, $U_1^{0+^*}$, $U_1^{+^*+}$ and $U_1^{++^*}$ are analogously given. 
For $i\in D$, we denote by $\bar{h}_i$ the mean service time of class-$i$ customers and by $\rho_i$ the offered load of Q$_i$; by the traffic equation, the offered loads are given by
\begin{align}
\rho_1=\bar{\lambda}_1 \bar{h}_1,\quad 
\rho_2=\bar{\lambda}_1 \bar{h}_2,\quad 
\rho_3=(p\bar{\lambda}_1+\bar{\lambda}_3) \bar{h}_3,\quad 
\rho_4=(p\bar{\lambda}_1+\bar{\lambda}_3) \bar{h}_4. 
\label{eq:offered_load}
\end{align}

%%%%%%%%%%%%%%%%%%%%%%%%
% figure: two-class MSP
%%%%%%%%%%%%%%%%%%%%%%%%

%%%%%%%%%%%%%%%%%%%%%%%%%%%%%%%%%%%%%%

%\medskip
{\it CTMC.}\quad 
Define a vector $\bar{\bJ}_t$ as $\bar{\bJ}_t=(\bar{J}_t(1),\bar{J}_t(2),\bar{J}_t(3),\bar{J}_t(4))$, then the stochastic process $\{\bar{\bY}_t\}=\{(\bar{\bX}_t,\bar{\bJ}_t)\}$ representing the behavior of the two-station network becomes a CTMC on the state space $\calS = \mathbb{Z}_+^4 \times S_J$, where $S_J = S^s_1 \times S^s_2 \times S^a_1 \times S^a_3$. 
Denote by $\bar{Q}=(\bar{Q}(\bx,\bx'),\,\bx,\bx'\in\mathbb{Z}_+^4)$ the infinitesimal generator of the Markov chain $\{\bar{\bY}_t\}$, where, for each $\bx,\bx'\in\mathbb{Z}_+^4$, $\bar{Q}(\bx,\bx')$ is a square block with the dimension of $s^a_1\,s^a_3\,|S^s_1|\,|S^s_2|$. Since the expression of $\bar{Q}$ is very lengthy and we do not use it in the paper, we omit to describe it.
We assume the Markov chain $\{\bar{\bY}_t\}$ is semi-irreducible and denote by $\calS_0$ its unique irreducible class. 
The state space $\calS$ is represented as 
\[
\calS = \bigcup_{A\in\scrP(D)} (\calB^A\cap\mathbb{Z}^d)\times S_J 
\]
and the transition rates of $\{\bar{\bY}_t\}$ are space-homogeneous in each subset $(\calB^A\cap\mathbb{Z}^d)\times S_J$; hence, $\{\bar{\bY}_t\}$ is a continuous-time version of a 4-dimensional skip-free semi-irreducible MMRRW. 
Let $A$ be a nonempty element of $\scrP(D)$ and $\{\bar{\bY}^A_t\}=\{(\bar{\bX}^A_t,\bar{\bJ}^A_t)\}$ the expanded Markov chain of $\{\bar{\bY}_t\}$ with index $A$; $\{(\bar{\bX}^A_t(D\setminus A),\bar{\bJ}^A_t)\}$ is the induced Markov chain of $\{\bar{\bY}_t\}$ with index $A$. We assume each induced Markov chain satisfies the condition of Assumption \ref{as:inducedMC}.
If the induced Markov chain $\{(\bar{\bX}^A_t(D\setminus A),\bar{\bJ}^A_t)\}$ is stable in our sense, we denote by $\bar{\bmu}^A$ the departure rate vector of $\{\bar{\bY}^A_t\}$.

%%%%%%%%%%%%%%%%%%%%%%%%%%%%
%
\subsubsection{Stability and instability conditions}

Consider the CTMC $\{\bar{\bY}_t\}=\{(\bar{\bX}_t,\bar{\bJ}_t)\}$. Since the process $\{\bar{\bX}_t\}$ is skip free in all coordinates, the diagonal elements of $\bar{Q}$ are bounded, i.e., for some positive number $\nu$, we have
\[
\sup_{\bx\in\mathbb{Z}_+^4} \max_{j\in S_0} \left|[\bar{Q}(\bx,\bx)]_{j,j} \right| \le \nu <\infty, 
\]
where, for a matrix $A$, we denote by $[A]_{i,j}$ the $(i,j)$-element of $A$.
Define a transition probability matrix $P=(P(\bx,\bx'), \bx,\bx'\in\mathbb{Z}_+^4)$ as $P = I + \bar{Q}/\nu$, 
where we have, for $\bx,\bx'\in\mathbb{Z}_+^4$, 
\[
P(\bx,\bx') = \left\{ \begin{array}{ll} 
I+\bar{Q}(\bx,\bx')/\nu & \mbox{if $\bx=\bx'$}, \cr 
\bar{Q}(\bx,\bx')/\nu & \mbox{if $\bx\ne\bx'$}. \end{array} \right.
\]
%
%By setting $\nu$ at an appropriate value, $P$ becomes aperiodic; hence, hereafter, we assume $P$ is aperiodic. 
%
Let $\calL=\{\bY_n\}=\{(\bX_n,\bJ_n)\}$ be a discrete-time Markov chain governed by the transition probability matrix $P$, where $\bX_n=(X_n(1),X_n(2),X_n(3),X_n(4))$ and $\bJ_n=(J_n(1),J_n(2),J_n(3),J_n(4))$. Since $\{\bY_t\}$ is a continuous-time version of a semi-irreducible MMRRW, the Markov chain $\{\bY_n\}$ is a semi-irreducible MMRRW with the same irreducible class $\calS_0$. 
Since $\{\bar{\bY}_t\}$ and $\{\bY_n\}$ has the same stationary distribution if it exists, we see that $\{\bar{\bY}_t\}$ is stable in our sense if $\{\bY_n\}$ is; we also see that $\{\bar{\bY}_t\}$ is unstable in our sense if $\{\bY_n\}$ is. Therefore, we analyze stability of $\{\bY_n\}$ instead of that of $\{\bar{\bY}_t\}$. 
The Markov chain $\{\bY_n\}$ represents a discrete-time two-station network corresponding to the continuous-time two-station network we consider. The discrete-time two-station network has the same network configuration as that of the original two-station network. 
For $i=1,3$, the arrival process of Q$_i$ is a discrete-time MAP with representation $(C_i,D_i)$, where $C_i=I+\bar{C}_i/\nu$ and $D_i=\bar{D}_i/\nu$, and, for $i=1,2$, the service process of station $i$ is a discrete-time MSP whose representation is given by $U_i^{k^*l},k,l=0,+$, \ $U_i^{kl^*},k,l=0,+$, and 
\begin{align*}
&T_i^{kl}=I+\bar{T}_i^{kl}/\nu,\ k,l=0,+,\quad
T_i^{k^*l}=\bar{T}_i^{k^*l}/\nu,\ k=1,2,\,l=0,+,\cr
&T_i^{kl^*}=\bar{T}_i^{kl^*}/\nu,\ k=0,+,\,l=1,2.
\end{align*}
For $i=1,3$, the mean arrival rate of class-$i$ customers is given by $\lambda_i=\bar{\lambda}_i/\nu$ and, for $i\in D$, the mean service time of class-$i$ customers is given by $h_i=\nu \bar{h}_i$. Hence, for $i\in D$, the offered load of Q$_i$ is given by $\rho_i$ of expression (\ref{eq:offered_load}). 

Since the offered load of station 1 is given by $\rho_1+\rho_4$ and that of station 2 by $\rho_2+\rho_3$, the nominal condition for the two-station network is given by 
\begin{equation}
\rho_1+\rho_4<1,\quad 
\rho_2+\rho_3<1. 
\label{eq:nominal}
\end{equation}
Since we are interested in the case where the nominal condition is not sufficient for the two-station network to be stable, we focus on the case where the MMRRW $\{\bY_n\}$ satisfies the following conditions.
\begin{assumption} \label{as:two_queue_spiral}
The index set of the stable induced Markov chains is given by 
\begin{equation}
\scrD_{stable} = \{D,\{1,2,3\},\{1,3,4\},\{1,4\},\{2,3\}\}, 
\end{equation}
and the mean drift vectors, $\ba^A,\,A\in\scrD_{stable}$, satisfy 
\begin{align}
& a^D(1)>0,\quad a^D(2)<0,\quad a^D(3)>0,\quad a^D(4)<0, \label{eq:aD} \\
& a^{\{1,2,3\}}(1)<0,\quad a^{\{1,2,3\}}(2)>0,\quad a^{\{1,2,3\}}(3)>0,\quad a^{\{1,2,3\}}(4)=0, \label{eq:a123} \\
& a^{\{1,3,4\}}(1)>0,\quad a^{\{1,3,4\}}(2)=0,\quad a^{\{1,3,4\}}(3)<0,\quad a^{\{1,3,4\}}(4)>0, \label{eq:a134} \\
& a^{\{1,4\}}(1)>0,\quad a^{\{1,4\}}(2)=0,\quad a^{\{1,4\}}(3)=0,\quad a^{\{1,4\}}(4)<0, \label{eq:a14} \\
& a^{\{2,3\}}(1)=0,\quad a^{\{2,3\}}(2)<0,\quad a^{\{2,3\}}(3)>0,\quad a^{\{2,3\}}(4)=0.  \label{eq:a23}
\end{align}
\end{assumption}

\bigskip
{\it Consistency of the conditions in Assumption \ref{as:two_queue_spiral}.} 
We verify consistency of the conditions in the assumption above. 
Since induced Markov chain $\calL^D$ is a semi-irreducible finite Markov chain, it is stable in our sense and the mean drift vector $\ba^D$ exists; hence, we assume condition (\ref{eq:aD}).

Next, consider induced Markov chain $\calL^{\{1,2,3\}}$, which is a one-dimensional MMRRW with just one induced Markov chain $\calL^{\{1,2,3\},\{4\}}$. By Remark \ref{re:inducedMC}, $\calL^{\{1,2,3\},\{4\}}$ is equivalent to the induced Markov chain $\calL^D$ and the mean drift of the corresponding expanded Markov chain $\hat{\calL}^{\{1,2,3\},\{4\}}$ is given by $a^D(4)$. Since $a^D(4)<0$, we see, by Theorem \ref{th:classification_1DMMRRW}, the induced Markov chain $\calL^{\{1,2,3\}}$ is stable in our sense. 
Analogously, we see induced Markov chain $\calL^{\{1,3,4\}}$ is stable in our sense. Thus, the mean drift vectors $\ba^{\{1,2,3\}}$ and $\ba^{\{1,3,4\}}$ exist and we assume conditions (\ref{eq:a123}) and (\ref{eq:a134}).
In a similar manner, it can also be seen that induced Markov chains $\calL^{\{1,2,4\}}$ and $\calL^{\{2,3,4\}}$ are unstable in our sense. 

Consider induced Markov chain $\calL^{\{1,4\}}$, which is a two-dimensional MMRRW. $\calL^{\{1,4\}}$ has three induced Markov chains $\calL^{\{1,4\},\{2\}}$, $\calL^{\{1,4\},\{3\}}$ and $\calL^{\{1,4\},\{2,3\}}$, which are equivalent to $\calL^{\{1,2,4\}}$, $\calL^{\{1,3,4\}}$ and $\calL^D$, respectively. 
Thus, $\calL^{\{1,4\},\{2\}}$ is unstable; $\calL^{\{1,4\},\{3\}}$ and $\calL^{\{1,4\},\{2,3\}}$ are stable and the mean drift vectors of the corresponding expanded Markov chains $\hat{\calL}^{\{1,4\},\{3\}}$ and $\hat{\calL}^{\{1,4\},\{2,3\}}$ are given by $\ba^{\{1,3,4\}}(\{2,3\})$ and $\ba^D(\{2,3\})$, respectively. Since $a^{\{1,3,4\}}(2)=0$, $a^{\{1,3,4\}}(3)<0$, $a^D(2)<0$ and $a^D(3)>0$, we see, by Theorem \ref{th:classification_2DMMRRW}, the induced Markov chain $\calL^{\{1,4\}}$ is stable in our sense. 
Analogously, we see induced Markov chain $\calL^{\{2,3\}}$ is stable in our sense. Thus, the mean drift vectors $\ba^{\{1,4\}}$ and $\ba^{\{2,3\}}$ exist and we assume conditions (\ref{eq:a14}) and (\ref{eq:a23}). 
In a similar manner, it can also be seen that induced Markov chains $\calL^{\{1,2\}}$, $\calL^{\{1,3\}}$, $\calL^{\{2,4\}}$ and $\calL^{\{3,4\}}$ are unstable in our sense. 

Finally, consider induced Markov chain $\calL^{\{1\}}$, which is a three-dimensional MMRRW. $\calL^{\{1\}}$ has four stable induced Markov chains $\calL^{\{1\},\{4\}}$, $\calL^{\{1\},\{2,3\}}$, $\calL^{\{1\},\{3,4\}}$ and $\calL^{\{1\},\{2,3,4\}}$ and the mean drift vectors of the corresponding expanded Markov chains $\hat{\calL}^{\{1\},\{4\}}$, $\hat{\calL}^{\{1\},\{2,3\}}$, $\hat{\calL}^{\{1\},\{3,4\}}$ and $\hat{\calL}^{\{1\},\{2,3,4\}}$ are given by $\ba^{\{1,4\}}(\{2,3,4\})$, $\ba^{\{1,2,3\}}(\{2,3,4\})$, $\ba^{\{1,3,4\}}(\{2,3,4\})$ and $\ba^D(\{2,3,4\})$, respectively. 
By conditions (\ref{eq:aD}) through (\ref{eq:a14}), we see that this case corresponds to C2-3-1 in Table \ref{tab:table2} and hence the induced Markov chain $\calL^{\{1\}}$ is unstable in our sense. 
In a similar manner, we see that induced Markov chains $\calL^{\{2\}}$, $\calL^{\{3\}}$ and $\calL^{\{4\}}$ are also unstable in our sense. 
This completes the verification of the consistency.
\hfill$\Box$

\bigskip
By expressions (\ref{eq:aAA_nw}) and (\ref{eq:aA_nw}), the nonzero elements of the mean drift vectors $\ba^D$, $\ba^{\{1,2,3\}}$, $\ba^{\{1,3,4\}}$,$\ba^{\{1,4\}}$ and $\ba^{\{2,3\}}$ are given as follows: 
\begin{align}
&a^D(1) = \big(\bar{\lambda}_1-\bar{\mu}^D(1)\big)/\nu,\quad 
a^D(2) = \big(\bar{\mu}^D(1)-\bar{\mu}^D(2)\big)/\nu,\cr
&a^D(3) = \big(\bar{\lambda}_3+p \bar{\mu}^D(2)-\bar{\mu}^D(3)\big)/\nu,\quad
a^D(4) = \big(\bar{\mu}^D(3)-\bar{\mu}^D(4)\big)/\nu, 
\label{eq:aD_nw} \\
&a^{\{1,2,3\}}(1) = \big(\bar{\lambda}_1-\bar{\mu}^{\{1,2,3\}}(1)\big)/\nu,\quad 
a^{\{1,2,3\}}(2) = \big(\bar{\mu}^{\{1,2,3\}}(1)-\bar{\mu}^{\{1,2,3\}}(2)\big)/\nu,\cr 
&a^{\{1,2,3\}}(3) = \big(\bar{\lambda}_3+p \bar{\mu}^{\{1,2,3\}}(2)-\bar{\mu}^{\{1,2,3\}}(3)\big)/\nu,
\label{eq:a123_nw} \\
&a^{\{1,3,4\}}(1) = \big(\bar{\lambda}_1-\bar{\mu}^{\{1,3,4\}}(1)\big)/\nu,\quad 
a^{\{1,3,4\}}(3) = \big(\bar{\lambda}_3+p \bar{\mu}^{\{1,3,4\}}(1)-\bar{\mu}^{\{1,3,4\}}(3)\big)/\nu,\cr 
&a^{\{1,3,4\}}(4) = \big(\bar{\mu}^{\{1,3,4\}}(3)-\bar{\mu}^{\{1,3,4\}}(4)\big)/\nu,
\label{eq:a134_nw} \\
&a^{\{1,4\}}(1) = \big(\bar{\lambda}_1-\bar{\mu}^{\{1,4\}}(1)\big)/\nu,\quad 
a^{\{1,4\}}(4) = \big(\bar{\lambda}_3+p \bar{\mu}^{\{1,4\}}(1)-\bar{\mu}^{\{1,4\}}(4)\big)/\nu,
\label{eq:a14_nw} \\
&a^{\{2,3\}}(2) = \big(\bar{\lambda}_1-\bar{\mu}^{\{2,3\}}(2)\big)/\nu,\quad 
a^{\{2,3\}}(3) = \big(\bar{\lambda}_3+p \bar{\mu}^{\{2,3\}}(2)-\bar{\mu}^{\{2,3\}}(3)\big)/\nu, 
\label{eq:a23_nw}
\end{align}
where, for $A\in\scrD_{stable}$ and for $i\in A$, $\bar{\mu}^A(i)$ is the departure rate of Q$_i$ in the original two-station network when Q$_j$, $j\in A$, are saturated with customers.
For $A\in\scrD_{stable}$ and for $i,j\in A$, let $r^A_{i,j}$ be defined as $r^A_{i,j} = \left| a^A(i)/a^A(j) \right|$ and let $r_1$ and $r_2$ be defined as 
\begin{equation}
r_1 = r^{\{1,2,3\}}_{2,1}\,r^{\{2,3\}}_{3,2}+r^{\{1,2,3\}}_{3,1},\quad 
r_2 = r^{\{1,3,4\}}_{4,3}\,r^{\{1,4\}}_{1,4}+r^{\{1,3,4\}}_{1,3}. 
\end{equation}
Note that $r_1$ and $r_2$ as well as $r^A_{i,j}$ does not depend on the uniformization parameter $\nu$. Under Assumption \ref{as:two_queue_spiral}, we have the following results. 

\begin{theorem} \label{th:two_station_stability} 
Assume at least one of the states in which every queue in the two-station network is empty belongs to the irreducible class $\calS_0$. Assume the conditions in Assumption \ref{as:two_queue_spiral} and further assume that $r^{\{1,4\}}_{1,4}\ge r^D_{1,4}$ and $r^{\{2,3\}}_{3,2} \ge r^D_{3,2}$. 
Then, the semi-irreducible MMRRW $\{\bY_n\}$ is stable in our sense if $r_1 r_2<1$ and it is unstable in our sense if $r_1 r_2>1$. 
Hence, the continuous-time two-station network we consider is stable in our sense if $r_1 r_2<1$ and it is unstable in our sense if $r_1 r_2>1$. 
\end{theorem}

\begin{proof}
First, assuming $r_1r_2<1$, we prove, by Theorem \ref{th:stable}, that $\{\bY_n\}$ is stable in our sense. 
Let $\delta$ be a positive number satisfying inequalities $\delta>1$ and 
\[
\big(r^{\{1,2,3\}}_{2,1} r^{\{2,3\}}_{3,2} \delta+r^{\{1,2,3\}}_{3,1}\big) \big(r^{\{1,3,4\}}_{4,3} r^{\{1,4\}}_{1,4} \delta+r^{\{1,3,4\}}_{1,3}\big) <1. 
\]
It is possible since we assume $r_1r_2<1$. 
Set positive vector $\bw$ as 
\begin{equation}
\bw = \Big( \big(r^{\{1,2,3\}}_{2,1} r^{\{2,3\}}_{3,2} \delta+r^{\{1,2,3\}}_{3,1}\big) \delta,\ r^{\{2,3\}}_{3,2} \delta,\ 1,\ r^{\{1,4\}}_{1,4} \big(r^{\{1,2,3\}}_{2,1} r^{\{2,3\}}_{3,2} \delta+r^{\{1,2,3\}}_{3,1}\big) \delta^2 \biggr), 
\label{eq:two_station_w}
\end{equation}
then we obtain, by Assumption \ref{as:two_queue_spiral}, 
\begin{align*}
&\langle \ba^{\{1,4\}},\bw \rangle = a^{\{1,4\}}(1) \big(r^{\{1,2,3\}}_{2,1} r^{\{2,3\}}_{3,2} \delta+r^{\{1,2,3\}}_{3,1}\big) \delta(1-\delta) < 0, \\
&\langle \ba^{\{2,3\}},\bw \rangle = -a^{\{2,3\}}(2)\,r^{\{2,3\}}_{3,2} (1-\delta) < 0,\\
&\langle \ba^{\{1,2,3\}},\bw \rangle = -a^{\{1,2,3\}}(1) \big( r^{\{1,2,3\}}_{2,1} r^{\{2,3\}}_{3,2} \delta+r^{\{1,2,3\}}_{3,1} \big) (1-\delta) < 0,\\
&\langle \ba^{\{1,3,4\}},\bw \rangle = -a^{\{1,3,4\}}(3) \Big( \big( r^{\{1,2,3\}}_{2,1} r^{\{2,3\}}_{3,2} \delta+r^{\{1,2,3\}}_{3,1} \big) \big(r^{\{1,3,4\}}_{4,3} r^{\{1,4\}}_{1,4} \delta+r^{\{1,3,4\}}_{1,3}\big) -1 \Big) < 0.
\end{align*}
Furthermore, we obtain, by the assumption of the theorem, 
\begin{align*}
&\langle \ba^D,\bw \rangle = -a^D(4)\big(r^{\{1,2,3\}}_{2,1} r^{\{2,3\}}_{3,2} \delta+r^{\{1,2,3\}}_{3,1}\big) \delta\big(r^D_{1,4}-r^{\{1,4\}}_{1,4}\delta\big) -a^D(2) \big(r^D_{3,2}-r^{\{2,3\}}_{3,2}\delta\big) < 0.
\end{align*}
Therefore, by Theorem \ref{th:stable}, the MMRRW $\{\bY_n\}$ is stable in our sense. 

Next, assuming $r_1r_2>1$, we prove, by Corollary \ref{co:unstable2}, that $\{\bY_n\}$ is unstable in our sense. 
By the assumption of the theorem, ($\{0\}\times\{0\}\times\{0\}\times\{0\}\times S_0)\cap\calS_0\ne\emptyset$ and hence we have $\calV_\emptyset\cap\calS_0\ne\emptyset$. 
Let $\delta$ be a positive number satisfying inequalities $\delta<1$ and 
\[
\big(r^{\{1,2,3\}}_{2,1} r^{\{2,3\}}_{3,2} \delta+r^{\{1,2,3\}}_{3,1}\big) \big(r^{\{1,3,4\}}_{4,3} r^{\{1,4\}}_{1,4} \delta+r^{\{1,3,4\}}_{1,3}\big) >1. 
\]
It is possible since we assume $r_1r_2>1$. 
Set $m=2$ and set $A_1=\{1,2,3\}$ and $A_2=\{1,3,4\}$; this satisfies conditions (\ref{eq:Ak_condition1}) and (\ref{eq:Ak_condition2}). By expressions (\ref{eq:barDfstable}) and (\ref{eq:barIfA}), we obtain $\bar{\scrD}_{f,stable}=\scrD_{stable}$ and 
\[
\bar{\scrI}_{f,\{1,4\}} = \{2\},\quad 
\bar{\scrI}_{f,\{2,3\}} = \{1\},\quad 
\bar{\scrI}_{f,\{1,2,3\}} = \{1\},\quad 
\bar{\scrI}_{f,\{1,3,4\}} = \{2\},\quad 
\bar{\scrI}_{f,D} = \{1,2\}.
\]
Set $\bw_0$ as $\bw_0=\bw>0$, where $\bw$ is given by expression (\ref{eq:two_station_w}) and set $c_0$ as
\[
c_0 = \frac{|\langle\ba^D,\bw_0\rangle|}{\min\{|a^D(2)|,\,|a^D(4)|\}}\delta^{-1}>0. 
\]
For $k\in\{1,2\}$, $\bw_k$ is given as $\bw_k=\bw_0-c_0\sum_{k'\in D\setminus A_k} \be_{k'}$ and $f_k$ as $f_k(\bx)=\langle\bx,\bw_k\rangle$ for $\bx\in\mathbb{R}^4$.  
Considering $\delta<1$, we obtain, by Assumption \ref{as:two_queue_spiral},
\begin{align*}
&f_2(\ba^{\{1,4\}}) = \langle\ba^{\{1,4\}},\bw_2\rangle = \langle\ba^{\{1,4\}},\bw\rangle > 0,\cr
&f_1(\ba^{\{2,3\}}) = \langle\ba^{\{2,3\}},\bw_1\rangle = \langle\ba^{\{2,3\}},\bw\rangle > 0,\cr
&f_1(\ba^{\{1,2,3\}}) = \langle\ba^{\{1,2,3\}},\bw_1\rangle = \langle\ba^{\{1,2,3\}},\bw\rangle > 0,\cr
&f_2(\ba^{\{1,3,4\}}) = \langle\ba^{\{1,3,4\}},\bw_2\rangle = \langle\ba^{\{1,3,4\}},\bw\rangle > 0. 
\end{align*}
Furthermore, for $k\in\bar{\scrI}_{f,D}=\{1,2\}$, we have
\[
f_k(\ba^D) = \langle\ba^D,\bw_k\rangle 
\ge |\langle\ba^D,\bw_0\rangle| \biggl(-1+\frac{\sum_{k'\in D\setminus A_k} |a^D(k')|}{\min\{|a^D(2)|,\,|a^D(4)|\}}\delta^{-1}\biggr)>0.
\]
Therefore, by Corollary \ref{co:unstable2}, the MMRRW $\{\bY_n\}$ is unstable in our sense. 
\end{proof}

%%%%%%%%%%%%%%%%%%%%%%%%%%%%%%%%%%%%%%%%%%%%%%%%%%%%%%%%%%%
%
% fig2.tex
%
%%%%%%%%%%%%%%%%%%%%%%%%%%%%%%%%%%%%%%%%%%%%%%%%%%%%%%%%%%%
\begin{figure}[bht]
\begin{center}
\setlength{\unitlength}{1.2mm}
\begin{picture}(70,35)(0,0)
\thinlines
\put(32,35){\makebox(0,0){\normalsize $x_1$}}
\put(40,3){\line(-1,4){7.5}}
\put(13.5,20){\makebox(0,0){\normalsize $x_2$}}
\put(40,3){\line(-3,2){24}}
\put(60,12){\makebox(0,0){\normalsize $x_3$}}
\put(40,3){\line(2,1){17}}
\put(58,28.5){\makebox(0,0){\normalsize $x_4$}}
\put(40,3){\line(2,3){16}}
\put(40,0){\makebox(0,0){\normalsize $O$}}
\thicklines
\put(36.5,29){\makebox(0,0){\normalsize $P_1$}}
\put(33.8,28){\vector(-1,-4){3.5}}
\multiput(20,16)(2,-0.5){5}{\line(1,0){1}}
\put(31,11.5){\makebox(0,0){\normalsize $P_2$}}
\put(30.5,13.8){\vector(4,-1){20.5}}
\put(52,6){\makebox(0,0){\normalsize $P_3$}}
\put(51,8.5){\vector(-1,3){4.5}}
\multiput(51,21)(-2,0.5){3}{\line(1,0){1}}
\put(47,25){\makebox(0,0){\normalsize $P_4$}}
\put(47,22){\vector(-4,1){12.5}}
\put(38.5,21.5){\makebox(0,0){\normalsize $P_5$}}
\end{picture}
\caption{A spiral path on the second vector field.}
\label{fig:spiral}
\end{center}
\end{figure}
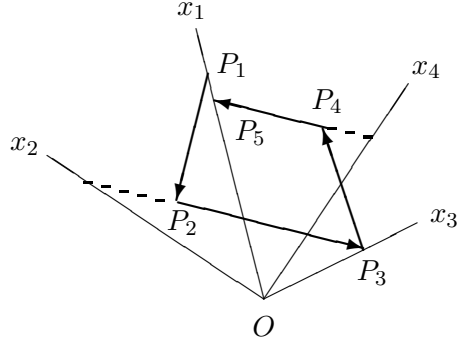
%%%%%%%%%%%%%%%%%%%%%%%%%%%%%%%%%%%%%%%%%%%%%%%%%%%%%%%%%%%

\begin{remark} \label{re:vectorfield}
Inequality $r_1 r_2<1$ corresponds to the condition that spiral paths on the second vector field introduced by Malyshev and Menshikov \cite{Malyshev81} (also see  Fayolle et al.\ \cite{Fayolle95}) reach the origin. Here we briefly explain this point. 
Recall that, since the dimension $d$ is four, sub-boundary $\calB^A$ for $A\in\scrP(D)$ is given as 
\[
\calB^A = \{(x(1),x(2),x(3),x(4))\in\mathbb{R}_+^4 : \mbox{$x(l)>0$ for $l\in A$; $x(l)=0$ for $l\in D\setminus A$} \}, 
\]
where $\calB^\emptyset$ includes only the origin and $\calB^D$ is the interior of $\mathbb{R}_+^4$. 
The second vector field generated from $\{\bY_n\}$ is a vector field constructed by assigning the mean drift vector $\ba^A$ to each point $\bx\in\calB^A$ if the induced Markov chain $\calL^A$ is stable, where $A$ is a nonempty element of $\scrP(D)$; if $\calL^A$ is unstable, some vector is assigned to each $\bx\in\calB^A$ but we here omit to describe it; see Malyshev and Menshikov \cite{Malyshev81} for details. 
Let us consider a spiral path on the second vector field, depicted in Fig.~\ref{fig:spiral}. We assume the path starts from point $P_1$ on the $x_1$-axis toward sub-boundary $\calB^{\{1,2,3\}}$.
Since vector $\ba^{\{1,2,3\}}$ is assigned to each point on $\calB^{\{1,2,3\}}$, the path continues further through $\calB^{\{1,2,3\}}$ toward sub-boundary $\calB^{\{2,3\}}$, following the vector $\ba^{\{1,2,3\}}$, and reaches point $P_2$ on $\calB^{\{2,3\}}$. 
Then, following vector $\ba^{\{2,3\}}$, the path continues through $\calB^{\{2,3\}}$ toward sub-boundary $\calB^{\{3\}}$ and reaches point $P_3$ on $\calB^{\{3\}}$, on which we suppose the path continues toward sub-boundary $\calB^{\{1,3,4\}}$.
The path continues further through $\calB^{\{1,3,4\}}$ toward sub-boundary $\calB^{\{1,4\}}$, following vector $\ba^{\{1,3,4\}}$, and reaches point $P_4$ on $\calB^{\{1,4\}}$. 
Then, following vector $\ba^{\{1,4\}}$, the path continues through $\calB^{\{1,4\}}$ toward $\calB^{\{1\}}$ and reaches point $P_5$ on $\calB^{\{1\}}$. 
After that, the path continues in the same way. 
We, therefore, see that if $P_5$ is closer to the origin than $P_1$, the spiral path eventually reaches the origin. Supposing $P_1=(1,0,0,0)$, we obtain, through some calculation, that $P_3=(0,0,x(3),0)$, where $x(3)$ is given by
\[
x(3) = \frac{a^{\{1,2,3\}}(2)}{a^{\{1,2,3\}}(1)} \frac{a^{\{2,3\}}(3)}{a^{\{2,3\}}(2)}-\frac{a^{\{1,2,3\}}(3)}{a^{\{1,2,3\}}(1)} = r_1. 
\]
Supporting $P_3=(0,0,1,0)$, we also obtain $P_5=(x(1),0,0,0)$, where $x(1)$ is given by
\[
x(1) = \frac{a^{\{1,3,4\}}(4)}{a^{\{1,3,4\}}(3)} \frac{a^{\{1,4\}}(1)}{a^{\{1,4\}}(4)}-\frac{a^{\{1,3,4\}}(1)}{a^{\{1,3,4\}}(3)} = r_2. 
\]
Hence, we see that if $r_1 r_2 <1$, the spiral path eventually reaches the origin. 
This is an intuitive explanation for the condition $r_1 r_2 <1$. 
\end{remark}

%%%%%%%%%%%%%%%%%%%%%%%%%%%%%%%%%%%%%
%
\subsubsection{Examples of service discipline}

We consider two models in this section: one is a two-station network with a non-preemptive priority service and the other that with a $(1,K)$-limited service. 

%%%%%%%%%%%%%%%%%%%%%%%%%%%%
{\it Non-preemptive priority service model.}\quad
Consider a two-station network in which class-4 customers have non-preemptive priority over class-1 customers and class-2 customers have non-preemptive priority over class-3 customers. 
For $i\in D$, we assume the service times of class-$i$ customers are subject to a phase-type (PH) distribution with representation $(\bbeta_i,\bar{H}_i)$, where $\bbeta_i$ and $\bar{H}_i$ are respectively the initial distribution and infinitesimal generator (restricted to transient states) of an absorbing Markov chain, by which the PH distribution is defined (see, for example, Latouche and Ramaswami \cite{Latouche99}). 
For $i\in D$, we denote by $\bar{h}_i$ the mean service time of class-$i$ customers; $\bar{h}_i$ is given by $\bar{h}_i=\bbeta_i (-\bar{H}_i)^{-1} \bone$. 
PH distributions are known to be dense in the space of distributions on $[0,\infty)$.  
Let the phase set of the PH-distribution for class-1 customers be given by $S^s_{11}=\{2,3,...,s^s_1+1\}$ and that for class-4 customers by $S^s_{14}=\{s^s_1+2, s^s_1+3,...,s^s_1+s^s_4+1\}$, where $s^s_1$ and $s^s_4$ are some positive integers. 
Let the phase set of the PH-distribution for class-3 customers be given by $S^s_{23}=\{2,3,...,s^s_3+1\}$ and that for class-2 customers by $S^s_{22}=\{s^s_3+2, s^s_3+3,...,s^s_3+s^s_2+1\}$, where $s^s_2$ and $s^s_3$ are some positive integers.
Then, letting $S^s_{10}=S^s_{20}=\{1\}$, we have $S^s_1=\{1,2,...,s^s_1+s^s_4+1\}$ and $S^s_2=\{1,2,...,s^s_3+s^s_2+1\}$. 
The representation of the MSP in station~1 is given in block form, as follows. 
\[
\bar{T}_1^{00} = \begin{pmatrix} 
0 & \bzero^\top & \bzero^\top \cr 
\bone & -I & O 
\cr \bone & O & -I \end{pmatrix},\quad 
\bar{T}_1^{+0} = \begin{pmatrix} 
-1 & \bbeta_1 & \bzero^\top \cr 
\bzero & \bar{H}_1 & O \cr 
\bzero & \bone \bbeta_1 & -I \end{pmatrix},\quad 
\bar{T}_1^{1^*0} = \begin{pmatrix} 
0 & \bzero^\top & \bzero^\top \cr 
-\bar{H}_1 \bone & O & O \cr 
\bzero & O & O \end{pmatrix},
\]
\[ 
\bar{T}_1^{2^*0} = \begin{pmatrix} 
0 & \bzero^\top & \bzero^\top \cr 
\bzero & -\bar{H}_1 \bone \bbeta_1 & O \cr 
\bzero & O & O \end{pmatrix},\quad 
\bar{T}_1^{0+} = \begin{pmatrix} 
-1 & \bzero^\top & \bbeta_4 \cr 
\bzero & -I & \bone \bbeta_4 \cr 
\bzero & O & \bar{H}_4 \end{pmatrix},\quad 
\bar{T}_1^{01^*} = \begin{pmatrix} 
0 & \bzero^\top & \bzero^\top \cr 
\bzero & O & O \cr 
-\bar{H}_4 \bone & O & O \end{pmatrix},
\]
\[
\bar{T}_1^{02^*} = \begin{pmatrix} 
0 & \bzero^\top & \bzero^\top \cr 
\bzero & O & O \cr 
\bzero & O & -\bar{H}_4 \bone \bbeta_4 \end{pmatrix},\quad 
\bar{T}_1^{++} = \begin{pmatrix} 
-1 & \bzero^\top & \bbeta_4 \cr 
\bzero & \bar{H}_1 & O \cr 
\bzero & O & \bar{H}_4 \end{pmatrix},
\]
\[
\bar{T}_1^{1^*+} = \bar{T}_1^{2^*+} = \begin{pmatrix} 
0 & \bzero^\top & \bzero^\top \cr 
\bzero & O & -\bar{H}_1 \bone \bbeta_4 \cr 
\bzero & O & O \end{pmatrix},\quad 
\bar{T}_1^{+1^*} = \begin{pmatrix} 
0 & \bzero^\top & \bzero^\top \cr 
\bzero & O & O \cr 
\bzero & -\bar{H}_4 \bone \bbeta_1 & O \end{pmatrix},
\]
\[
\bar{T}_1^{+2^*} = \begin{pmatrix} 
0 & \bzero^\top & \bzero^\top \cr 
\bzero & O & O \cr 
\bzero & O & -\bar{H}_4 \bone \bbeta_4 \end{pmatrix},\quad
U_1^{0^*0} = \begin{pmatrix} 
0 & \bbeta_1 & \bzero^\top \cr 
\bzero & \bone \bbeta_1 & O \cr 
\bzero & \bone \bbeta_1 & O \end{pmatrix},\quad 
U_1^{00^*} = \begin{pmatrix} 
0 & \bzero^\top & \bbeta_4 \cr 
\bzero & O & \bone \bbeta_4 \cr 
\bzero & O & \bone \bbeta_4 \end{pmatrix},
\]
where $O$ a matrix of $0$'s whose dimension is determined in context and $\bzero$ a column vector of $0$'s whose dimension is also determined in context; $U_1^{+^*0}$, $U_1^{+0^*}$, $U_1^{0^*+}$, $U_1^{0+^*}$, $U_1^{+^*+}$ and $U_1^{++^*}$ are identity matrices. 
Note that the matrices above have several elements corresponding to dummy transitions; for example, $\bone$'s and $-I$'s of $\bar{T}_1^{00}$ are blocks of such elements.
The representation of the MSP in station~2 are analogously given. 
In the CTMC $\{\bar{\bY}(t)\}$ arising from the two-station network with a non-preemptive priority service, the state $((0,0,0,0),(1,1,1,1))$ is accessible from all other states; hence, by Remark \ref{re:semiirreducible}, the CTMC $\{\bar{\bY}(t)\}$ is semi-irreducible and the state $((0,0,0,0),(1,1,1,1))$ belongs to the irreducible class $\calS_0$. 

Assuming nominal condition (\ref{eq:nominal}), we consider the case where $\bar{h}_1<\bar{h}_2$ and $\bar{h}_3<\bar{h}_4$; in this case, the nominal condition is not sufficient for the model to be stable. 
Supposing $\scrD_{stable}=\{D,\{1,2,3\},\{1,3,4\},\{1,4\},\{2,3\}\}$, we evaluate the mean drift vectors $\ba^A,\,A\in\scrD_{stable}$. 
\begin{itemize}
\item $\ba^D$. 
Consider the case where all the queues are saturated with customers. Since class-4 customers have non-preemptive priority over class-1 customers, the server in station 1 is engaged in service only for class-4 customers after some finite period of time. Hence, the mean departure rates of Q$_1$ and Q$_4$ are given by $\bar{\mu}^D(1)=0$ and $\bar{\mu}^D(4)=1/\bar{h}_4$, respectively. Analogously, we have $\bar{\mu}^D(2)=1/\bar{h}_2$ and $\bar{\mu}^D(3)=0$. 
By expressions (\ref{eq:aD_nw}), we obtain 
\begin{align*}
&a^D(1) = \bar{\lambda}_1/\nu >0,\quad 
a^D(2) = (-1/\bar{h}_2)/\nu <0,\quad 
a^D(3) = (\bar{\lambda}_3+p/\bar{h}_2)/\nu >0,\cr 
&a^D(4) = (-1/\bar{h}_4)/\nu <0, 
\end{align*}
where $\nu$ is the uniformization parameter. 
\item $\ba^{\{1,2,3\}}$.  
Consider the case where queues Q$_1$, Q$_2$ and Q$_3$ are saturated with customers. In this case, the server in station 2 is engaged in service only for class-2 customers after some finite period of time, and we have $\bar{\mu}^{\{1,2,3\}}(2)=1/\bar{h}_2$ and $\bar{\mu}^{\{1,2,3\}}(3)=0$. This implies that any customers do not arrive at Q$_4$ after some finite period of time and Q$_4$ will eventually become empty since customers in Q$_4$ have priority in service over customers in Q$_1$. 
After Q$_4$ becomes empty, the server in station 1 is engaged in service for customers in Q$_1$, which is saturated with customers, and we have $\bar{\mu}^{\{1,2,3\}}(1)=1/\bar{h}_1$.  Hence, by expressions (\ref{eq:a123_nw}), we obtain 
\begin{align*}
&a^{\{1,2,3\}}(1) = (\bar{\lambda}_1-1/\bar{h}_1)/\nu <0,\quad 
a^{\{1,2,3\}}(2) = (1/\bar{h}_1-1/\bar{h}_2)/\nu >0,\cr 
&a^{\{1,2,3\}}(3) = (\bar{\lambda}_3+p/\bar{h}_2)\nu >0,
\end{align*}
where we use the nominal condition and the assumption of  $\bar{h}_1<\bar{h}_2$.
\item $\ba^{\{1,3,4\}}$.  
Consider the case where queues Q$_1$, Q$_3$ and Q$_4$ are saturated with customers. This case is symmetric to that for $\ba^{\{1,2,3\}}$ and we have $\bar{\mu}^{\{1,3,4\}}(1)=0$, $\bar{\mu}^{\{1,3,4\}}(3)=1/\bar{h}_3$ and $\bar{\mu}^{\{1,3,4\}}(4)=1/\bar{h}_4$. Hence, by expression (\ref{eq:a134_nw}), we obtain 
\begin{align*}
&a^{\{1,3,4\}}(1) = \bar{\lambda}_1/\nu>0,\quad 
a^{\{1,3,4\}}(3) = (\bar{\lambda}_3-1/\bar{h}_3)/\nu<0,\cr
&a^{\{1,3,4\}}(4) = (1/\bar{h}_3-1/\bar{h}_4)/\nu >0, 
\end{align*}
where we use the nominal condition and the assumption of $\bar{h}_3<\bar{h}_4$.
\item $\ba^{\{1,4\}}$.  
Consider the case where queues Q$_1$ and Q$_4$ are saturated with customers. In this case, the server in station 1 is engaged in service only for class-4 customers after some finite period of time, and we have $\bar{\mu}^{\{1,4\}}(1)=0$ and $\bar{\mu}^{\{1,4\}}(4)=1/\bar{h}_4$. 
Hence, by expression (\ref{eq:a14_nw}), we obtain 
\[
a^{\{1,4\}}(1) = \bar{\lambda}_1/\nu>0,\quad 
a^{\{1,4\}}(4) = (\bar{\lambda}_3-1/\bar{h}_4)/\nu<0, 
\]
where we use the nominal condition. 
\item $\ba^{\{2,3\}}$.  
Consider the case where queues Q$_2$ and Q$_3$ are saturated with customers. This case is symmetric to that for $\ba^{\{1,4\}}$ and we have $\bar{\mu}^{\{2,3\}}(2)=1/\bar{h}_2$ and $\bar{\mu}^{\{2,3\}}(3)=0$. Hence, by expression (\ref{eq:a23_nw}), we have 
\[
a^{\{2,3\}}(2) = (\bar{\lambda}_1-1/\bar{h}_2)/\nu <0,\quad 
a^{\{2,3\}}(3) = (\bar{\lambda}_3+p/\bar{h}_2)/\nu >0, 
\]
where we use the nominal condition. 
\end{itemize}

By the arguments above, we see that the model satisfies the conditions in Assumption \ref{as:two_queue_spiral}. Furthermore, the following conditions of Theorem \ref{th:two_station_stability} are also satisfied: 
\begin{align*}
&\left|\frac{a^D(1)}{a^D(4)}\right|=\frac{\bar{\lambda}_1}{1/\bar{h}_4} \le \frac{\bar{\lambda}_1}{1/\bar{h}_4-\bar{\lambda}_3} = \left|\frac{a^{\{1,4\}}(1)}{a^{\{1,4\}}(4)}\right|, \cr
&\left|\frac{a^D(3)}{a^D(2)}\right|=\frac{\bar{\lambda}_3+p/\bar{h}_2}{1/\bar{h}_2} \le \frac{\bar{\lambda}_3+p/\bar{h}_2}{1/\bar{h}_2-\bar{\lambda}_1}=\left|\frac{a^{\{2,3\}}(3)}{a^{\{2,3\}}(2)}\right|. 
\end{align*}
$r_1$ and $r_2$ in Theorem \ref{th:two_station_stability} are given as 
\[
r_1 = \frac{\lambda_3+p \mu_2}{\mu_2-\lambda_1},\quad 
r_2 = \frac{\lambda_1}{\mu_4-\lambda_3}, 
\]
and inequality $r_1 r_2<1$ is equivalent to $\rho_2+\rho_3<1$. Hence, by Theorem \ref{th:two_station_stability}, we see that the two-station network with a non-preemptive priority service is stable in our sense if $\rho_2+\rho_3<1$ and it is unstable in our sense if $\rho_2+\rho_3>1$; this result is coincident with the existing results for two-station networks with a preemptive-resume priority service (see, for example, Bramson \cite{Bramson08} and Dai and Weiss \cite{Dai96b}).

\begin{remark}
It can be seen from the results above that the stability and instability conditions for the model are given only in terms of the mean arrival rates and the mean service times; hence, the stability region of the model does not depend on other features of the arrival processes and service time distributions. 
On the other hand, Dai et al.~\cite{Dai04} demonstrated that the stability region of a push-started Lu-Kumar network depended on the inter-arrival time distribution and service time distributions. Our two-station network is a simplified version of the push-started Lu-Kumar network, and it can be seen from the results in Dai et al.~\cite{Dai04} that, in the case of constant inter-arrival and service times and in the case of uniformly distributed inter-arrival and service times, the condition $\rho_2+\rho_4>1$ does not always imply instability of the two-station network with a non-preemptive priority service. 
Any arrival process with i.i.d.\ inter-arrival times can be approximated with any accuracy by a MAP and any service time distribution can also be approximated with any accuracy by a PH-distribution. However, any arrival process whose inter-arrival time distribution has bounded support cannot be represented as a MAP and any service time distribution with bounded support cannot also be represented as a PH-distribution. 
Therefore, as mentioned in Dai et al.~\cite{Dai04}, it can be considered that the boundedness of the support of the distributions affects the stability of the network. 
\end{remark}

%%%%%%%%%%%%%%%%%%%%%%%%%
{\it $(1,K)$-limited service model.}\quad 
Consider a two-station network in which customers in Q$_1$ (Q$_3$) are served according to a 1-limited service and those in Q$_4$ (resp.\ Q$_2$) according to a $K$-limited service, which means that the server in station 1 (resp.\ station 2) alternatively visits Q$_1$ and Q$_4$ (resp.\ Q$_3$ and Q$_2$), serves one customer upon a visit to Q$_1$ (resp.\ Q$_3$) and continuously serves customers upon a visit to Q$_4$ (resp.\ Q$_2$) until it completes serving just $K$ customers or Q$_4$ (resp.\ Q$_2$) becomes empty. 
We assume that switchover times of the servers are negligible. We call this service discipline a $(1,K)$-limited service. 
The $(1,K)$-limited service is equivalent to a non-preemptive priority service when the value of $K$ is unbounded (i.e., $K=\infty$). Using a $(1,K)$-limited service, we can control relative levels of priority between customers in different classes by varying the value of $K$.
For $i\in D$, we assume the service times of class-$i$ customers are subject to a PH-distribution with the same representation as that of the non-preemptive priority service model. 
The phase sets $S^s_{10}$, $S^s_{11}$, $S^s_{20}$ and $S^s_{23}$ are also the same as those of the non-preemptive service model. On the other hand, in order to represent how many class-4 customers (class-2 customers) have continuously been served in Q$_4$ (resp.\ Q$_2$), $S^s_{14}$ and $S^s_{22}$ are given as $S^s_{14}=\{ s^s_1+2, s^s_1+3, ..., s^s_1+K s^s_4+1 \}$ and $S^s_{22}=\{ s^s_3+2, s^s_3+3, ..., s^s_3+K s^s_2+1 \}$. 
The representation of the MSP in station~1 is given in block form, as follows (we omit several zeros in describing matrices): 
\[
\bar{T}_1^{00} = \begin{pmatrix} 
0 & & & &\cr 
\bone & -I & & & \cr
\bone & & -I & & \cr
\vdots & & & \ddots & \cr
\bone & & & & -I  \end{pmatrix},\quad 
\bar{T}_1^{+0} = \begin{pmatrix} 
-1 & \bbeta_1 & & &\cr 
 & \bar{H}_1 & & & \cr
 & \bone \bbeta_1 & -I & & \cr
 & \vdots & & \ddots & \cr
 & \bone \bbeta_1 & & & -I  \end{pmatrix}, 
\]
\[
\bar{T}_1^{1^*0} = \begin{pmatrix} 
0 & & & &\cr 
-\bar{H}_1 \bone & O & & & \cr
 & & O & & \cr
 & & & \ddots & \cr
 & & & & O  \end{pmatrix},\quad 
\bar{T}_1^{2^*0} = \begin{pmatrix} 
0 & & & &\cr 
 & -\bar{H}_1 \bone \bbeta_1 & & & \cr
 & & O & & \cr
 & & & \ddots & \cr
 & & & & O  \end{pmatrix}, 
\]
\[
\bar{T}_1^{0+} = \begin{pmatrix} 
-1 & & \bbeta_4 & &\cr 
 & -I & \bone \bbeta_4 & & \cr
 & & \bar{H}_4 & & \cr
 & & & \ddots & \cr
 & & & & \bar{H}_4  \end{pmatrix},\quad 
\bar{T}_1^{01^*} = \begin{pmatrix} 
0 & & & &\cr 
 & O & & & \cr
-\bar{H}_4 \bone & & O & & \cr
\vdots & & & \ddots & \cr
-\bar{H}_4 \bone & & & & O  \end{pmatrix}, 
\]
\[
\bar{T}_1^{02^*} = \begin{pmatrix} 
0 & & & &\cr 
 & O & & & \cr
 & & O & -\bar{H}_4 \bone \bbeta_4 & \cr
 & & & \ddots & \ddots \cr
 & & -\bar{H}_4 \bone \bbeta_4 & & O  \end{pmatrix},\quad 
\bar{T}_1^{++} = \begin{pmatrix} 
-1 & & \bbeta_4 & &\cr 
 & \bar{H}_1 & & & \cr
 & & \bar{H}_4 & & \cr
 & & & \ddots & \cr
 & & & & \bar{H}_4  \end{pmatrix},
\]
\[
\bar{T}_1^{1^*+} = \bar{T}_1^{2^*+} = \begin{pmatrix} 
0 & & & &\cr 
 & O & -\bar{H}_1 \bone \bbeta_4 & & \cr
 & & O & & \cr
 & & & \ddots &  \cr
 & & & & O  \end{pmatrix},\quad 
\bar{T}_1^{+1^*} = \begin{pmatrix} 
0 & & & &\cr 
 & O & & & \cr
 & -\bar{H}_4 \bone \bbeta_1 & O & & \cr
 & \vdots & & \ddots &  \cr
 & -\bar{H}_4 \bone \bbeta_1 & & & O  \end{pmatrix},
\]
\[
\bar{T}_1^{+2^*} = \begin{pmatrix} 
0 & & & &\cr 
 & O & & & \cr
 & & O & -\bar{H}_4 \bone \bbeta_4 & \cr
 & & & \ddots & \ddots  \cr
 & -\bar{H}_4 \bone \bbeta_1 & & & O  \end{pmatrix},\quad 
U_1^{0^*0} = \begin{pmatrix} 
0 & \bbeta_1 & & &\cr 
 & \bone \bbeta_1 & & & \cr
 & \bone \bbeta_1 & O & & \cr
 & \vdots & & \ddots & \cr
 & \bone \bbeta_1 & & & O \end{pmatrix},\quad 
\]
\[
U_1^{00^*} = \begin{pmatrix} 
0 & & \bbeta_4 & & \cr 
 & O & \bone \bbeta_4 & & \cr
 & & \bone \bbeta_4 & & \cr
 & & \vdots & \ddots & \cr
 & & \bone \bbeta_4 & & O \end{pmatrix}. 
\]
$U_1^{+^*0}$, $U_1^{+0^*}$, $U_1^{0^*+}$, $U_1^{0+^*}$, $U_1^{+^*+}$ and $U_1^{++^*}$ are identity matrices. The representation of the MSP in station~2 are analogously given. 
In the CTMC $\{\bar{\bY}(t)\}$ arising from the two-station network with a $(1,K)$-limited service, the state $((0,0,0,0),(1,1,1,1))$ is accessible from all other states; hence, by Remark \ref{re:semiirreducible}, the CTMC $\{\bar{\bY}(t)\}$ is semi-irreducible and the state $((0,0,0,0),(1,1,1,1))$ belongs to the irreducible class $\calS_0$. 

For simplicity, we assume the model parameters are symmetric, i.e., $\bar{\lambda}_1=\bar{\lambda}_3$, $\bar{h}_1=\bar{h}_3$, $\bar{h}_2=\bar{h}_4$ and $p=0$, and examine how the stability region of the two-station network is affected by the value of the parameter $K$. 
We assume the nominal condition and condition $\bar{h}_1<\bar{h}_2$. 
Furthermore, we assume $K$ satisfies 
\begin{align}
K>K^*=\max\!\left\{1,\,\frac{1-\rho_1}{\rho_2},\,\frac{\rho_1}{1-\rho_2} \right\}, 
\label{eq:Kcond}
\end{align}
where $\rho_1=\bar{\lambda}_1\bar{h}_1$ and $\rho_2=\bar{\lambda}_1\bar{h}_2$. In order for the model to satisfy the conditions in Assumption \ref{as:two_queue_spiral}, condition (\ref{eq:Kcond}) will be necessary. 
Supposing $\scrD_{stable}=\{D,\{1,2,3\},\{1,3,4\},\{1,4\},\{2,3\}\}$, we evaluate the mean drift vectors $\ba^A,\,A\in\scrD_{stable}$. 
\begin{itemize}
%%%%%
\item $\ba^D$. 
Consider the case where all the queues are saturated with customers. In station 1, after some finite period of time, the server serves one customer in Q$_1$ and then serves $K$ customers in Q$_4$; after that, the server repeats it. Hence, the departure rates of Q$_1$ and Q$_4$ are given by 
\[
\bar{\mu}^D(1)=\frac{1}{\bar{h}_1+K\bar{h}_2},\quad 
\bar{\mu}^D(4)=\frac{K}{\bar{h}_1+K\bar{h}_2}. 
\]
Since the model is symmetric, we have $\bar{\mu}^D(2)=\bar{\mu}^D(4)$ and $\bar{\mu}^D(3)=\bar{\mu}^D(1)$. 
By expression~(\ref{eq:aD_nw}), we obtain 
\begin{align*}
&a^D(1) = a^D(3) = \left(\bar{\lambda_1}-\frac{1}{\bar{h}_1+K\bar{h}_2}\right)/\nu = \frac{\rho_1+K\rho_2-1}{\bar{h}_1+K\bar{h}_2}\frac{1}{\nu} >0,\cr 
&a^D(2) = a^D(4) = \left(\frac{1}{\bar{h}_1+K\bar{h}_2}-\frac{K}{\bar{h}_1+K\bar{h}_2}\right)/\nu = -\frac{K-1}{\bar{h}_1+K\bar{h}_2}\frac{1}{\nu} <0, 
\end{align*}
where we use inequalities $K>(1-\rho_1)/\rho_2$ and $K>1$, which come from condition (\ref{eq:Kcond}); $\nu$ is the uniformization parameter. 
%%%%%
\item $\ba^{\{1,2,3\}}$ and $\ba^{\{1,3,4\}}$.  
Consider the case where queues Q$_1$, Q$_2$ and Q$_3$ are saturated with customers. In station 2, after some finite period of time, the server serves $K$ customers in Q$_2$ and then serves one customer in Q$_3$; after that, the server repeats it. Hence, the departure rates of Q$_2$ and Q$_3$ are given by 
\[
\bar{\mu}^{\{1,2,3\}}(2)=\frac{K}{\bar{h}_1+K\bar{h}_2},\quad 
\bar{\mu}^{\{1,2,3\}}(3)=\frac{1}{\bar{h}_1+K\bar{h}_2}. 
\]
In station 1, we have $a^{\{1,2,3\}}(4)=0$ and this implies $\bar{\mu}^{\{1,2,3\}}(4)=\bar{\mu}^{\{1,2,3\}}(3)$. Hence, we obtain 
\[
\bar{\mu}^{\{1,2,3\}}(1) = \left( 1-\bar{\mu}^{\{1,2,3\}}(3)\bar{h}_2 \right)/\bar{h}_1 = \frac{1+(K-1)\bar{h}_2/\bar{h}_1}{\bar{h}_1+K\bar{h}_2}. 
\]
By expression (\ref{eq:a123_nw}), we obtain 
\begin{align*}
&a^{\{1,2,3\}}(1) = \left(\bar{\lambda}_1 - \frac{1+(K-1)\bar{h}_2/\bar{h}_1}{\bar{h}_1+K\bar{h}_2}\right)/\nu = \frac{(K-1)(\rho_2-\bar{h}_2/\bar{h}_1)+\rho_1+\rho_2-1}{\bar{h}_1+K\bar{h}_2}\frac{1}{\nu} <0,\cr
&a^{\{1,2,3\}}(2) = \left(\frac{1+(K-1)\bar{h}_2/\bar{h}_1}{\bar{h}_1+K\bar{h}_2} - \frac{K}{\bar{h}_1+K\bar{h}_2}\right)/\nu = \frac{(K-1)(\bar{h}_2/\bar{h}_1-1)}{\bar{h}_1+K\bar{h}_2}\frac{1}{\nu} >0,\cr
&a^{\{1,2,3\}}(3) = \left(\bar{\lambda}_1 - \frac{1}{\bar{h}_1+K\bar{h}_2}\right)/\nu = \frac{\rho_1+K\rho_2-1}{\bar{h}_1+K\bar{h}_2}\frac{1}{\nu} > 0, 
\end{align*}
where we use the nominal condition, condition $\bar{h}_1<\bar{h}_2$ and inequality $K>(1-\rho_1)/\rho_2$. 
Since the model is symmetric, we obtain 
\begin{align*}
&a^{\{1,3,4\}}(1) = a^{\{1,2,3\}}(3) >0,\quad 
a^{\{1,3,4\}}(3) = a^{\{1,2,3\}}(1) <0,\quad
a^{\{1,3,4\}}(4) = a^{\{1,2,3\}}(2) >0.
\end{align*}
%%%%%
\item $\ba^{\{1,4\}}$ and $\ba^{\{2,3\}}$.  
Consider the case where queues Q$_1$ and Q$_4$ are saturated with customers. In station 1, after some finite period of time, the server serves one customer in Q$_1$ and then serves $K$ customers in Q$_4$; after that, the server repeats it. Hence, the output rates of Q$_1$ and Q$_4$ are given by 
\[
\bar{\mu}^{\{1,4\}}(1) = \frac{1}{\bar{h}_1+K\bar{h}_2},\quad 
\bar{\mu}^{\{1,4\}}(4) = \frac{K}{\bar{h}_1+K\bar{h}_2}. 
\]
By expression (\ref{eq:a14_nw}), we obtain  
\begin{align*}
&a^{\{1,4\}}(1) = \left(\bar{\lambda}_1 - \frac{1}{\bar{h}_1+K\bar{h}_2}\right)/\nu = \frac{\rho_1+K\rho_2-1}{\bar{h}_1+K\bar{h}_2}\frac{1}{\nu} > 0,\cr
&a^{\{1,4\}}(4) = \left(\bar{\lambda}_1 - \frac{K}{\bar{h}_1+K\bar{h}_2}\right)/\nu = \frac{\rho_1+K\rho_2-K}{\bar{h}_1+K\bar{h}_2}\frac{1}{\nu} < 0,
\end{align*}
where we use inequalities $K>(1-\rho_1)/\rho_2$ and $K>\rho_1/(1-\rho_2)$, which come from condition (\ref{eq:Kcond}). 
Since the model is symmetric, we obtain 
\begin{align*}
&a^{\{2,3\}}(2) = a^{\{1,4\}}(4) <0,\quad 
a^{\{2,3\}}(3) = a^{\{1,4\}}(1) >0.
\end{align*}
\end{itemize}

By the arguments above, we see that the model satisfies the conditions in Assumption \ref{as:two_queue_spiral}. Furthermore, the following conditions of Theorem \ref{th:two_station_stability} are also satisfied: 
\begin{align*}
&\left|\frac{a^D(1)}{a^D(4)}\right|=\left|\frac{a^D(3)}{a^D(2)}\right|
=\frac{\rho_1+K\rho_2-1}{K-1} 
\le \frac{\rho_1+K\rho_2-1}{K-\rho_1-K\rho_2}
= \left|\frac{a^{\{1,4\}}(1)}{a^{\{1,4\}}(4)}\right|=\left|\frac{a^{\{2,3\}}(3)}{a^{\{2,3\}}(2)}\right|, 
\end{align*}
where we use inequalities $K>(1-\rho_1)/\rho_2$. 
$r_1$ and $r_2$ in Theorem \ref{th:two_station_stability} are given as 
\[
r_1 = r_2 = \frac{\rho_1+K\rho_2-1}{-\rho_1+K(1-\rho_2)}, 
\]
and inequality $r_1 r_2<1$ is equivalent to 
\[
\left(\frac{\rho_1+K\rho_2}{1+K}\right)<\frac{1}{2}. 
\]
Therefore, by Theorem \ref{th:two_station_stability}, if $\rho_2<1/2$, then the two-station network with a $(1,K)$-limited service is stable in our sense for all $K>K^*$; if $\rho_2>1/2$, then it is stable in our sense for $K$ such that $K^*<K<\frac{1-2 \rho_1}{2 \rho_2-1}$ and unstable in our sense for $K$ such that $K>\max\!\left\{K^*,\,\frac{1-2\rho_1}{2\rho_2-1} \right\}$.
The condition $\rho_2>1/2$ corresponds to $\rho_2+\rho_4>1$ in the asymmetric parameter case. Hence, we see that even though a two-station network with a non-preemptive priority service is unstable, the corresponding two-station network with a $(1,K)$-limited service may be stable for some values of $K$. 
For example, when the parameters are set as $\lambda_1=1$, $\mu_1=5$ and $\mu_2=9/5$, we have $\rho_1+\rho_2=34/45<1$ and $\rho_2=5/9>1/2$. In this case, the two-station network with a $(1,K)$-limited service is stable if $2\le K\le 5$ and it is unstable if $K\ge 6$. 

\begin{remark}
In the symmetric parameter case, if $K=1$, then $a^D(2)=a^D(4)=0$ and hence we cannot apply our results to the $(1,K)$-limited service model. 
\end{remark}

%%%%%%%%%%%%%%%%%%%%%%%
%
% References
%
%%%%%%%%%%%%%%%%%%%%%%%
%


\begin{thebibliography}{99}
%
\bibitem{Bramson08}
Bramson, M.\ (2008). 
Stability of queueing networks. 
Probability Survey {\bf 5}, 169--345. 
%
\bibitem{Bramson10}
Bramson, M., Dai, J.G.\ and Harrison J.M.\  (2010). 
Positive recurrence of reflecting Brownian motion in three dimensions. 
The Annals of Applied Probability {\bf 20(2)}, 753--783.
%
\bibitem{Bremaud99}
Br\'emaud, P.\ (1999).
{\it Markov Chains: Gibbs Fields, Monte Carlo Simulation and Queues}. 
Springer, New York. 
%
\bibitem{Chen95} 
Chen, H.\ (1995). 
Fluid approximations and stability of multiclass queueing networks I: Work-conserving discipline. 
The Annals of Applied Probability 5(3), 637--665.
%
\bibitem{Dai95}
Dai, J.G.\  (1995). 
On positive Harris recurrence of multiclass queueing networks: a unified approach via fluid limit models. 
The Annals of Applied Probability {\bf 5(1)}, 49--77.
%
\bibitem{Dai96}
Dai, J.G.\ (1996).
A fluid limit model criterion for instability of multiclass queueing networks.
The Annals of Applied Probability 6(3), 751--757.
%
\bibitem{Dai96c}
Dai, J.G.\ and Vande Vate, J.\ (1996). 
Virtual stations and the capacity of two-station queueing networks. 
Under revision for Operations Research. 
%
\bibitem{Dai96b}
Dai, J.G.\ and Weiss, G.\ (1996). 
Stability and instability of fluid models for re-entrant lines. 
Mathematics of Operations Research {\bf 21}, 115--134.
%
\bibitem{Dai00}
Dai, J.G.\ and Vande Vate, J.\ (2000). 
The stability of two-station multi-type fluid networks. 
Operations Research {\bf 48}, 721--744. 
%
\bibitem{Dai04}
Dai, J.G., Hasenbein, J.J.\ and Vande Vate, J.H.\ (2004).
Stability and instability of a two-station queueing network.
The Annals of Applied Probability 14, 326--377.
%
\bibitem{Down97}
Down, D.\ and Meyn, S.P.\ (1997).
Piecewise linear test functions for stability and instability of queueing networks. 
Queueing Systems 27, 205--226.
%
\bibitem{Dumas97}
Dumas, V.\ (1997).
A multiclass network with non-linear, non-convex, non-monotonic stability condition. 
Queueing Systems 25, 1-43.
%
\bibitem{Fayolle89}
Fayolle, G.\  (1989). 
On random walks arising in queueing systems: ergodicity and transience via quadratic forms as Lyapounov functions -- Part I. 
Queueing Systems {\bf 5}, 167--184.
%
\bibitem{Fayolle95} 
Fayolle, G., Malyshev, V.A., and Menshikov, M.V.\  (1995). 
{\it Topics in the Constructive Theory of Countable Markov Chains}. 
Cambridge University Press, Cambridge.
%
\bibitem{Gamarnik05}
Gamarnik, D.\ and Hasenbein J.J.\ (2005). 
Instability in stochastic and fluid queueing networks, 
The Annals of Applied Probability 15(3), 1652--1690.
%
\bibitem{Kharroubi02}
Kharroubi, A.E., Tahar, A.B. and Yaacoubi, A.\  (2002). 
On the stability of the linear Skorohod problem in an orthant. 
Mathematical Methods of Operations Research {\bf 56}, 243--258. 
%
\bibitem{Kumar90}
Kumar, P.R.\ and Seidman, T.I.\ (1990). 
Dynamic instabilities and stabilization methods in distributed real-time scheduling of manufacturing systems. 
IEEE Trans.\ Automat.\ Control {\bf 35}, 289--298.
%
%\bibitem{Kumar93}
%Kumar P.R.\ (1993).
%Re-entrant lines. 
%Queueing Systems 13, 87--110.
%
\bibitem{Latouche99}
Latouche, G.\ and Ramaswami, V.\  (1999).  
{\it Introduction to Matrix Analytic Methods in Stochastic Modeling}. 
SIAM, Philadelphia.
%
\bibitem{Lu91}
Lu, S.H.\ and Kumar, P.R.\ (1991). 
Distributed scheduling based on due dates and buffer priorities. 
IEEE Trans.\ Automat.\ Control {\bf 36}, 1406--1416.
%
\bibitem{Malyshev81} 
Malyshev, V.A.\ and Menshikov, M.V.\ (1981). 
Ergodicity, continuity, and analyticity of countable Markov chains. 
Transactions of the Moscow Mathematical Society 1, 1--47.
%
\bibitem{Meyn94}
Meyn, S.P.\ and Down, D.\  (1994).  
Stability of generalized Jackson networks. 
The Annals of Applied Probability {\bf 4(1)}, 124--148.
%
\bibitem{Meyn95}
Meyn, S.P.\ (1995).
Transience of multiclass queueing networks via fluid limit models. 
The Annals of Applied Probability 5(4), 946--957.
%
\bibitem{Miyazawa11}
Miyazawa, M.\  (2011).
Light tail asymptotics in multidimensional reflecting processes for queueing networks. 
TOP {\bf 19(2)}, 233--299. 
%
\bibitem{Neuts94} 
Neuts, M.F.\  (1994).
{\it Matrix-Geometric Solutions in Stochastic Models}.   
Dover Publications, New York.
%
\bibitem{Ozawa04}
Ozawa, T.\  (2004). 
Analysis of queues with Markovian Service Processes. 
Stochastic Models {\bf 20(4)}, 391--413.
%
\bibitem{Ozawa12}
Ozawa, T.\  (2013). 
Asymptotics for the stationary distribution in a discrete-time two-dimensional quasi-birth-and-death process. 
Queueing Systems {\bf 74}, 109--149.
%
%\bibitem{Ozawa13b}
%Ozawa, T.\  (2013). 
%Positive recurrence and transience of  a two-station network with server states. 
%Submitted (arXiv: 1308.6104). 
%
\bibitem{Rybko92}
Rybko, A.\ and Stolyar, A.\ (1992). 
On the ergodicity of stochastic processes describing open queueing networks. 
Problemi Peredachi Informatsii {\bf 28}, 3--26.
%
\bibitem{Tezcan13}
Tezcan, T.\  (2013). 
Stability analysis of N-model systems under a static priority rule. 
Queueing Systems {\bf 73}, 235--259.
%
\end{thebibliography}
\end{document}